\documentclass [10pt,a4paper]{article}
\usepackage{amssymb}
\usepackage{enumerate,verbatim}
\usepackage{amsmath}
\usepackage{amsfonts,mathrsfs}
\usepackage{amsthm,wasysym}
\usepackage{epsfig,lpic,tikz}
\usetikzlibrary{calc}
\usepackage[applemac]{inputenc}
\usepackage[english]{babel}

\usepackage[width=0.9\textwidth]{caption}

\usepackage{pdfsync}
\usepackage{xcolor,hyperref}
\usepackage[T1]{fontenc}
\usepackage{lpic}
\usepackage{mathtools}

\hypersetup{
    colorlinks=true,
    linkcolor=blue,
    filecolor=magenta,      
    urlcolor=blue,%cyan,
    pdftitle={Overleaf Example},
    pdfpagemode=FullScreen,
    }

\definecolor{green}{rgb}{0,0.5,0}

\hypersetup{backref,colorlinks=true,linkcolor=blue,citecolor=blue}

\usepackage{marvosym}

\newcounter{teoremaganso}
\newcounter{appendix}

\newcounter{coryganso}

\renewenvironment{abstract}{\small\quotation\noindent
 {\bfseries \abstractname .}}{\endquotation \par}

\newenvironment{trivlistparm}[2]{\trivlist
\item[{#1 #2}]}{\endtrivlist}

\newenvironment{prooftext}[1]{\trivlistparm{\bfseries}{#1}}{\Qed\endtrivlistparm}
\newenvironment{prova}{\trivlistparm{\bfseries}{Proof.}}{\Qed\endtrivlistparm}

\catcode`\@=11

\def\resetthefootnote{\renewcommand{\thefootnote}{\@arabic\c@footnote} }
\def\@principiremex#1{\trivlist
 \item[\hskip \labelsep{\bfseries #1\ \thetheo.}]\ignorespaces}
\def\opar@principiremex#1[#2]{\trivlist
 \item[\hskip \labelsep{\bfseries #1\ \thetheo\ (#2).}]\ignorespaces}

\newcommand{\newTHEOremrom}[2]{\newenvironment{#1}{\refstepcounter{theo}\@ifnextchar[{\opar@principiremex{#2}}
{\@principiremex{#2}}}{\qedB\endtrivlist}} \catcode`\@=12
\DeclareMathSymbol{\square}{\mathord}{AMSa}{"03}
\newcommand{\qedB}{\nopagebreak\hspace*{\fill}$\square$\par}
\newcommand{\Qed}{\nopagebreak\hspace*{\fill}{\vrule width6pt height6pt depth0pt}\par}

\newtheorem {theo} {Theorem} [section]
\newtheorem {prop} [theo] {Proposition}
\newtheorem {cory} [theo] {Corollary}
\newtheorem {lem} [theo] {Lemma}
\newtheorem {bigtheo} [teoremaganso] {Theorem}

\newTHEOremrom {defi} {Definition}
\newTHEOremrom {obs} {Remark}
\newTHEOremrom {ex} {Example}

\newcommand{\refc}[1]{\mbox{$(\ref{#1})$}}
\newcommand{\secc}[1]{Section~\ref{#1}}

\newcommand{\teoc}[1]{Theorem~\ref{#1}}
\newcommand{\propc}[1]{Proposition~\ref{#1}}
\newcommand{\coryc}[1]{Corollary~\ref{#1}}
\newcommand{\lemc}[1]{Lemma~\ref{#1}}
\newcommand{\defic}[1]{Definition~\ref{#1}}
\newcommand{\obsc}[1]{Remark~\ref{#1}}

\newcommand{\figc}[1]{Figure~\ref{#1}}

\newcommand{\N}{\ensuremath{\mathbb{N}}}
\newcommand{\Z}{\ensuremath{\mathbb{Z}}}
\newcommand{\R}{\ensuremath{\mathbb{R}}}

\newcommand{\F}{\ensuremath{\mathcal{F}}}

\newcommand{\np}{{\ensuremath{{\mu}}}}

\newcommand{\h}{{\ensuremath{\ell_{n+1}}}}
\newcommand{\Sc}{\ensuremath{\mathbb{S}}}

\newcommand{\cc}{\ensuremath{\mathscr{C}}}

\newcommand{\mb}[1]{\ensuremath{\mathbb{#1}}}

\def\map#1#2#3{\mbox{${#1}\!:{#2}\longrightarrow{#3}$}}

\newcommand{\sist}[2]{
  \left\{\!
   \begin{array}{l}
    \dot x=#1 \\[2pt] \dot y=#2
   \end{array}
  \right.
}

\makeatletter
\newcommand{\lambdabar}{{\mathchoice
  {\smash@bar\textfont\displaystyle{0.25}{1.2}\lambda}
  {\smash@bar\textfont\textstyle{0.25}{1.2}\lambda}
  {\smash@bar\scriptfont\scriptstyle{0.25}{1.2}\lambda}
  {\smash@bar\scriptscriptfont\scriptscriptstyle{0.25}{1.2}\lambda}
}}
\newcommand{\smash@bar}[4]{%
  \smash{\rlap{\raisebox{-#3\fontdimen5#10}{$\m@th#2\mkern#4mu\mathchar'26$}}}%
}
\makeatother

\newcommand{\op}{\ensuremath{\mbox{\rm o}}}

\newcommand{\mc}[1]{\mathcal{#1}}
\newcommand{\mr}[1]{\mathrm{#1}}
\newcommand{\red}[1]{{\color{red} #1}}
\newcommand{\blue}[1]{{\color{blue}#1}}

\newcommand{\dsp}{\displaystyle}
\newcommand{\gorro}{\hat}
\title{\bf The cyclicity of hyperbolic hemicycles
\footnotetext{2010 {\it AMS Subject Classification}: 34C07; 34C20; 34C23.} 
\footnotetext{{\it Key words and phrases}: limit cycle, hemicycle, cyclicity, asymptotic expansion, Dulac map.}
\footnotetext{This work has been partially funded by the Ministry of Science, Innovation and Universities of Spain through the grants 
PID2020-118281GB-C33, PID2021-125625NB-I00 and PID2022-136613NB-I00, and by the Agency for Management of University and Research Grants of Catalonia through the grants 2021SGR00113 and 2021SGR01015.}
}

\author{D. Mar\'{\i}n and J. Villadelprat
\\*[.1truecm]
{\small \textsl{Departament de Matem{\`a}tiques, Edifici Cc,
Universitat Aut{\`o}noma de Barcelona,}}\\*[-.05truecm]
{\small\textsl{08193 Cerdanyola del Vallès (Barcelona), Spain}}
}
                  
\date{\today}

\begin{document}

\maketitle

\begin{abstract}
We consider families of planar polynomial vector fields of degree $n$ and study the cyclicity of a type of unbounded polycycle~$\Gamma$ called hemicycle. Compactified to the Poincar\'e disc,~$\Gamma$ consists of an affine straight line together with half of the line at infinity and has two singular points, which are hyperbolic saddles located at infinity. We prove four main results. \teoc{thmA} deals with the cyclicity of~$\Gamma$ when perturbed without breaking the saddle connections. For the other results we consider the case $n=2$. More concretely they are addressed to the quadratic integrable systems belonging to the class $Q_3^R$ and having two hemicycles, $\Gamma_u$ and $\Gamma_\ell$, surrounding each one a center. \teoc{thmB} gives the cyclicity of $\Gamma_u$ and $\Gamma_\ell$ when perturbed inside the whole family of quadratic systems. In \teoc{thmC} we study the number of limit cycles bifurcating simultaneously from $\Gamma_u$ and $\Gamma_\ell$ when perturbed as well inside the whole family of quadratic systems.  
Finally, in \teoc{ThmD} we show that for three specific cases there exists a simultaneous alien limit cycle bifurcation from $\Gamma_u$ and $\Gamma_\ell$.
\end{abstract}

%\tableofcontents

\section{Introduction and main results}

In this paper we consider families of planar polynomial vector fields $X_\mu$ depending on a parameter $\mu\in\R^N$ and we are interested in the number of limit cycles (i.e., isolated periodic orbits). More concretely, we study their bifurcations, which occur at the limit periodic sets of the family (where limit cycles accumulate as~$\mu$ varies). In this setting the first step is to obtain the sharp bound for the number of limit cycles that bifurcate from each limit periodic set $\Gamma.$ This bound is called the cyclicity of $\Gamma.$ The computation of the cyclicity is a crucial step to determine the bifurcation diagram for the number of limit cycles within the family. Before stating our main results we will give a precise definition of all these notions. The 
problems that we discuss in the present paper are related to questions surrounding Hilbert's 16th problem and its various weakened forms. We refer the interested reader to the monographs of  Il'yashenko~\cite{Ily}, Jibin~Li~\cite{Li}, or Roussarie~\cite{Roussarie} for more information on these issues.

We begin by recalling the notion of limit periodic set as introduced in {\cite[Definition 10]{Roussarie}}. This is the fundamental object that we aim to study and its definition is given in terms of the Hausdorff topology, which for reader's convenience we briefly explain next.
\begin{obs}\label{compact}
Let $S$ be a metrizable space and denote by $\mc C(S)$ the set of all compact non-empty subsets of~$S.$ Given any $K_1,K_2\in\mc C(S)$ we define
\[
 d_H(K_1,K_2)=\sup_{x_1\in K_1,x_2\in K_2}\left\{\inf_{x_2'\in K_2}d(x_1,x_2')\,,\inf_{x_1'\in K_1}d(x_1',x_2)\right\}.
\]
One can readily show that $d_H$ is a distance. It defines a topology on $\mc C(S)$, which is independent of the distance $d$ chosen, that is called the \emph{Hausdorff topology}. Moreover it turns out that
\[
 d_H(K_1,K_2)=\inf\big\{\varepsilon>0:K_1\subset N_\varepsilon(K_2)\text{ and }K_2\subset N_\varepsilon(K_1)\big\},
\]
where $N_\varepsilon(K)$ is the $\varepsilon$-neighbourhood of $K$. Finally, 
if $(S,d)$ is a compact metric space then so is $(\mc C(S),d_H)$. The interested reader is referred to \cite[p. 279]{Munkres} for both assertions.
\end{obs}

\begin{defi}\label{lps}
A  non-empty compact subset $\Gamma$ of a surface $S$  is a \emph{limit periodic set} for a germ of a family $\{X_\mu\}_{\mu\approx\mu_0}$ of vector fields on $S$ if there exists a sequence of parameters $\{\mu_n\}_{n\in\N}$ converging to $\mu_0$ such that each~$X_{\mu_n}$ has a limit cycle $\gamma_n$ and the sequence $\{\gamma_n\}_{n\in\N}$ converges to $\Gamma$ as $n\to\infty$ in the Hausdorff topology of the space $\mc C(S)$ of compact non-empty subsets of $S$. 
\end{defi}

It is well known, see \cite[Theorem 5]{Roussarie}, that  any limit periodic set of a germ of an analytic family $\{X_\mu\}_{\mu\approx\mu_0}$ such that $X_{\mu_0}$ has only isolated singular points is either a singular point, a period orbit or a graphic of $X_{\mu_0}$. We recall the notion of graphic and polycycle below:

\begin{defi}%\label{persistent}
Let $X$ be a vector field on $\R^2$ (or $\Sc^2$). A \emph{graphic} $\Gamma$ for $X$ is a compact, non-empty invariant subset which is a continuous image of $\Sc^1$ and consists of a finite number of isolated singularities $\{p_1,\ldots,p_m,p_{m+1}=p_1\}$ (not necessarily distinct) and compatibly oriented separatrices $\{s_1,\ldots,s_m\}$ connecting them (i.e., such that the $\alpha$-limit set of~$s_j$ is~$p_j$ and the $\omega$-limit set of~$s_j$ is~$p_{j+1}$). A graphic is said to be \emph{hyperbolic} if all its singular points are hyperbolic saddles. A \emph{polycycle} is a graphic with a return map defined on one of its sides. 
\end{defi}

The polycycles that we aim to study are unbounded and for this reason we need to compactify the vector field. 
Recall that to investigate the phase portrait of a polynomial vector field~$Y$ near infinity we can consider its Poincar\'e compactification $p(Y)$, see \cite[\S 5]{ADL} for details, which is an analytically equivalent vector field defined on the sphere $\Sc^2$. The points at infinity of~$\R^2$ are in bijective correspondence with the points of the equator of~$\Sc^2$, that we denote by~$\ell_{\infty}$. Moreover the trajectories of $p(Y)$ in $\Sc^2$ are symmetric with respect to the origin and so it suffices to draw its flow in the closed northern hemisphere only, the so called Poincar\'e disc. 

\begin{defi}\label{defcic}
Let $\Pi$ be an arbitrary collection of limit periodic sets for the germ of an analytic family $\{X_\mu\}_{\mu\approx\mu_0}$ of vector fields on $\Sc^2$. We define the \emph{cyclicity} of $\Pi$ with respect to  $\{X_\mu\}_{\mu\approx\mu_0}$ as
\[
\mr{Cycl}\big((\Pi,X_{\mu_0}),X_\mu\big)
%\mr{Cycl}_S(\Pi,\{X\}_\mu,\mu_0)
\!:=\inf\limits_{\varepsilon,\delta>0}\sup\limits_{\mu\in B_\delta(\mu_0)}\#\big\{\gamma \text{ limit cycle of $X_\mu$ such that }d_H(\gamma,\Gamma)<\varepsilon\text{ for some $\Gamma\in\Pi$}\big\},
\]
which may be infinite.
\end{defi}
\begin{obs}\label{cicl_sim}
Let us point out that if 
$\Pi=\{\Gamma\}$ then the cyclicity of $\Pi$ coincides with the usual cyclicity $\mr{Cycl}\big((\Gamma,X_{\mu_0}),X_\mu\big)$ of the limit periodic set $\Gamma$, cf. \cite[Definition~12]{Roussarie}. In contrast, if $\Pi$ consists of more than one limit periodic set then the cyclicity of $\Pi$ accounts for the limit cycles bifurcating \emph{simultaneously} from any of them. Finally, observe that if $\Pi\subset\Pi'$ then $\mr{Cycl}\big((\Pi,X_{\mu_0}),X_\mu\big)\leqslant\mr{Cycl}\big((\Pi',X_{\mu_0}),X_\mu\big)$.
\end{obs}
\begin{figure}[t]
\begin{center}
\begin{tikzpicture}
%\draw[step=1mm,gray,very thin] (-5,-8) grid (8,4);
\begin{scope}[shift={(-1.7,-3.5)},rotate=40,scale=0.5]
\path[fill=blue] (10.5,0.2) -- (10.5,-0.2) -- (11,0) -- cycle;
\end{scope}
\begin{scope}[shift={(6.8,-6.1)},rotate=150,scale=0.5]
\path[fill=blue] (10.5,0.2) -- (10.5,-0.2) -- (11,0) -- cycle;
\end{scope}
\begin{scope}[shift={(-1.7,-3.8)},rotate=30,scale=0.5]
\path[fill=red] (10.5,0.2) -- (10.5,-0.2) -- (11,0) -- cycle;
\end{scope}
\begin{scope}[shift={(7,-5.3)},rotate=150,scale=0.5]
\path[fill=red] (10.5,0.2) -- (10.5,-0.2) -- (11,0) -- cycle;
\end{scope}
\begin{scope}[shift={(-4.5,1.33)},rotate=-28,scale=0.5]
\path[fill=red] (10.5,0.2) -- (10.5,-0.2) -- (11,0) -- cycle;
\end{scope}
\begin{scope}[shift={(5,-.08)},rotate=210,scale=0.5]
\path[fill=red] (10.5,0.2) -- (10.5,-0.2) -- (11,0) -- cycle;
\end{scope}
\begin{scope}[shift={(3.7,-2.75)},rotate=180,scale=0.5]
\path[fill=blue] (10.5,0.2) -- (10.5,-0.2) -- (11,0) -- cycle;
\end{scope}
\begin{scope}[shift={(10,-2.75)},rotate=180,scale=0.5]
\path[fill=blue] (10.5,0.2) -- (10.5,-0.2) -- (11,0) -- cycle;
\end{scope}
\def\r{5}\def\s{1.5}
\draw [very thick,red] (\s,-2) to [out=135,in=90] (\s-\r,-2) to [out=-90, in =225] (\s,-2) to [out=45,in=90] (\s+\r,-2) to [out=-90,in=-45] (\s,-2);
 \draw [fill,red] (\s,-2) circle [radius=.1];
 \draw [] (\s+.75,-2) to [out=90,in=200] (\s+1.3,-1.4) to [out=20,in=90] (\s+\r-.5,-2) to [out=-90,in=0] (\s+3.,-3) to [out=180,in=-90] (\s+1.2,-2) to [out=90,in=200] (\s+1.5,-1.6);

 \draw [] (\s-.75,-2) to [out=90,in=-20] (\s-1.3,-1.4) to [out=160,in=90] (\s-\r+.5,-2) to [out=-90,in=180] (\s-3.,-3) to [out=0,in=-90] (\s-1.2,-2) to [out=90,in=-20] (\s-1.5,-1.6);

\draw [] (\s,-1.2) to [out=0,in=-145] (\s+1,-.75) to [out=35,in=90] (\s+\r+.6,-2) to [out=-90,in=-35] (\s+1,-3.25) to [out=145,in=0] (\s,-2.8) to [out=180,in=35] (\s-1,-3.25) to [out=-145,in=-90] (\s-\r-.6,-2) to [out=90,in=145] (\s-1,-0.6) to [out=-35, in=180] (\s,-.9) to [out=0,in=-145] (\s+1,-0.4);

\draw [thick,blue] (\s,-0.5) to [out=0,in=-145] (\s+1,0) to [out=35,in=90] (\s+\r+1,-2) to [out=-90,in=-35] (\s+1,-3.6) to [out=145,in=0] (\s,-3.3) to [out=180,in=35] (\s-1,-3.6) to [out=-145,in=-90] (\s-\r-1,-2) to [out=90,in=145] (\s-1,0) to [out=-35, in=180] (\s,-0.5);

 \draw [thick,blue] (4.5,-2) ellipse ({1.25cm} and {0.75cm}); 
 \draw [thick,blue] (-1.5,-2) ellipse ({1.25cm} and {0.75cm}); 
\node at (\s-0.5,-1.25) {$\Gamma_-$};
\node at (\s+0.5,-2.65) {$\Gamma_+$};
\node at (\s-2.1,-2) {$\gamma_-$};
\node at (\s+2.1,-2) {$\gamma_+$};
\node at (\s,-0.25) {$\gamma$};
\end{tikzpicture}
\end{center}
\caption{``Figure  eight-loop'' $\Gamma=\Gamma_-\cup\Gamma_+$ formed by two homoclinic connections $\Gamma_-$ and $\Gamma_+$. The limit cycles $\gamma_-$, $\gamma_+$ and $\gamma$ are close (with respect to the Hausdorff distance) to $\Gamma_-$, $\Gamma_+$ and $\Gamma$, respectively.}\label{dib8}
\end{figure}
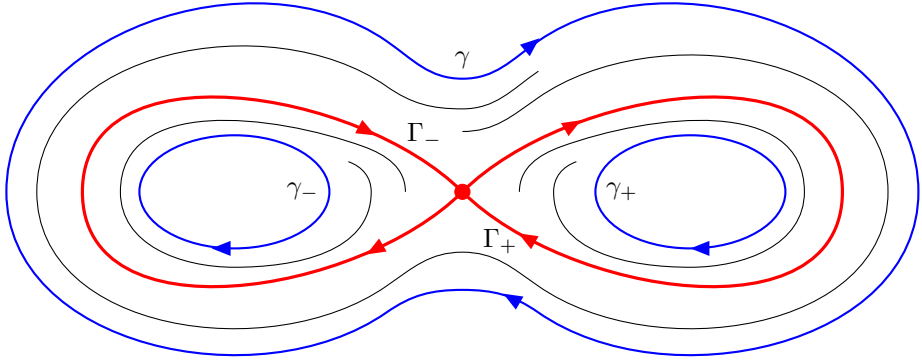

Note that the simultaneous cyclicity of $\{\Gamma_1,\ldots,\Gamma_r\}$ may not coincide with the cyclicity of $\Gamma_1\cup\cdots\cup\Gamma_r$, even in case that the latter is a limit periodic set. For instance, consider 
a germ $\{X_\mu\}_{\mu\approx\mu_0}$ such that $X_{\mu_0}$ has a saddle point  with two homoclinic loops $\Gamma_-$ and $\Gamma_+$ making up a  ``figure eight-loop'' $\Gamma=\Gamma_-\cup\Gamma_+$, see \figc{dib8}. Then the values of $\mathrm{Cycl}\big((\Pi,X_{\mu_0}),X_\mu\big)$ for 
\[
\Pi=\{\Gamma_+\},\; \Pi=\{\Gamma_-\},\; \Pi=\{\Gamma\} ,\; \Pi=\{\Gamma_+,\Gamma_-\},\; 
\Pi=\{\Gamma_+,\Gamma\},\;  \Pi=\{\Gamma_-,\Gamma\}\text{ and } \Pi=\{\Gamma_+,\Gamma_-,\Gamma\}
\]
%\[
%\begin{array}{llll}
%\mathrm{Cycl}\big((\Gamma_-,X_{\mu_0}),X_\mu\big) & \mathrm{Cycl}\big((\Gamma_+,X_{\mu_0}),X_\mu\big) & \mathrm{Cycl}\big((\Gamma,X_{\mu_0}),X_\mu\big) & \mathrm{Cycl}\big((\{\Gamma_-,\Gamma_+\},X_{\mu_0}),X_\mu\big)\\[10pt]
%\mathrm{Cycl}\big((\{\Gamma_+,\Gamma\},X_{\mu_0}),X_\mu\big) & \mathrm{Cycl}\big((\{\Gamma_-,\Gamma\},X_{\mu_0}),X_\mu\big)& \mathrm{Cycl}\big((\{\Gamma_-,\Gamma_+,\Gamma\},X_{\mu_0}),X_\mu\big) & \\
%\end{array}
%\]
may be all different. On the other hand, it is clear that 
\begin{equation*}%\label{cotasim}
\max_{j\in\{1,2,\ldots,r\}}\left\{\mr{Cycl}\big((\Gamma_j,X_{\mu_0}),X_\mu\big)\right\}\leqslant
\mathrm{Cycl}\big((\{\Gamma_1,\ldots,\Gamma_r\},X_{\mu_0}),X_\mu\big)\leqslant \sum_{j=1}^r\mathrm{Cycl}\big((\Gamma_j,X_{\mu_0}),X_\mu\big).
\end{equation*}

In this paper we study the cyclicity problem for perturbations of planar polynomial vector fields with an invariant straight line (see \cite{DRR2,GMM02} and references therein for previous results on the issue). After a suitable rotation we can assume that this invariant straight line is $\{y=0\}.$ In the first part of the paper this line is assumed to be invariant throughout all the perturbation. Any such family $\{X_\mu\}_{\mu\in\Lambda}$ can be written as
\begin{equation}\label{DS}
X_\mu\quad\sist{yf(x,y;\mu)+g(x;\mu),}{yq(x,y;\mu),}
\end{equation}
where $\Lambda$ is an open subset of $\R^N$ and $f$, $g$ and $q$ are polynomials with the coefficients depending analytically on $\mu.$ We assume that $\deg(f)=\deg(q)=n$ and $\deg(g)=n+1$ and that the following hypothesis hold:
\begin{enumerate}

\item[\bf{H1}] $g(x;\mu)<0$ for all $x\in\R$ and $\mu\in\Lambda$, which implies that $n$ is odd, and

\item[\bf{H2}] $\h(x,y;\mu)\!:=yf_n(x,y;\mu)-xq_n(x,y;\mu)+g_{n+1}(\mu)x^{n+1}>0$
for all $(x,y)\neq (0,0)$ and  $\mu\in\Lambda$.

\end{enumerate}
Here, and in what follows, $f_n(x,y;\mu)$ and $q_n(x,y;\mu)$ denote, respectively, the homogeneous part of degree~$n$ of $f(x,y;\mu)$ and $q(x,y;\mu)$, whereas $g_{n+1}(\mu)$ is the leading coefficient of $g(x;\mu).$ The second hypothesis is related with the angle variation $\theta$ of the solutions of \refc{DS} near the infinity because one can verify that 
\[
 r^2\dot\theta=y\big(xq(x,y)-yf(x,y)-g(x)\big).
\]  
Since $\h$ is a homogeneous polynomial of even degree, {\bf H2} is equivalent to 
$zf_n(1,z)-q_n(1,z)+g_{n+1}>0$ and $f_n(z,1)-zq_n(z,1)+g_{n+1}z^{n+1}>0$  for all $z\in\R$ and $\mu\in\Lambda.$

\begin{figure}[t]
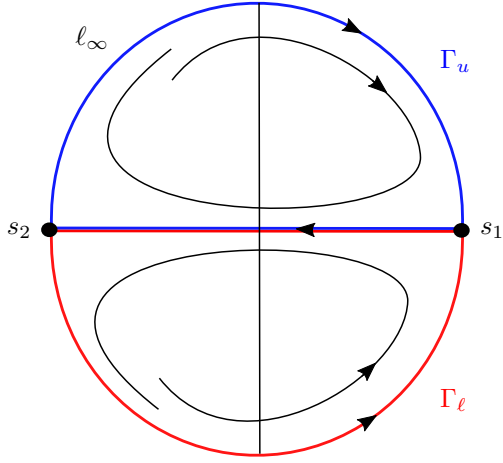

   \centering
  \begin{lpic}[l(0mm),r(0mm),t(0mm),b(5mm)]{dib3(0.75)}
    \lbl[l]{-6,40;$s_2$}   
    \lbl[l]{6,74;$\ell_\infty$}       
    \lbl[l]{77,40;$s_1$}     
    \lbl[l]{70,70;$\blue{\Gamma_u}$}           
    \lbl[l]{70,10;$\red{\Gamma_\ell}$}              
   \end{lpic}
  \caption{Placement of the hyperbolic saddles and the polycycles $\Gamma_u$ and $\Gamma_\ell$ in the Poincar\'e disc for the polynomial vector field \refc{DS}.}
  \label{dib3}
 \end{figure}

Conditions {\bf H1} and {\bf H2} guarantee that, after compactifying the polynomial vector field $X_\mu$ to the Poincar\'e disc, the boundary of the upper (respectively, lower) half-plane is a polycycle $\Gamma_u$ (respectively, $\Gamma_\ell$) with two hyperbolic saddles, see \figc{dib3},
\[
 \text{$s_1\!:=\{y=0,x>0\}\cap \ell_\infty$ and $s_2\!:=\{y=0,x<0\}\cap \ell_\infty.$}
\] 
This type of polycycle, formed by an invariant line and half of the equator of $\Sc^2$, is called \emph{hemicycle} in~\cite{DRR2}. Moreover the vector fields of the form \refc{DS} verifying {\bf H1} and {\bf H2} are called \emph{D-systems} by the authors in \cite{GMM02}.

Our first main result is addressed to the cyclicity of $\Gamma_u$ when perturbed \emph{inside} the family of $D$-systems. This result will be given in terms of two functions $d_0(\mu)$ and $d_1(\mu)$. In order to define them we first need to introduce several other functions. For the sake of shortness we shall omit the dependence of $\mu$ in these functions when there is no risk of confusion. We define
\begin{equation}\label{K}
 K(x_1,x_2;\mu)\!:=\left. 1-\frac{xq(x,y)}{yf(x,y)+g(x)}\right|_{(x,y)=\left(\frac{1}{x_1},\frac{x_2}{x_1}\right)}
 \text{ and }\lambda(\mu)\!:=-K(0,0;\mu)=-1+\frac{q_n(1,0)}{g_{n+1}}>0.
\end{equation}
The function $K$ is related to the projective compactification of $X_\mu$, whereas $\lambda(\mu)$ is the hyperbolicity ratio of its saddle at infinity. Let us remark that, on account of {\bf H1} and {\bf H2}, the functions $K$ and $1/K$ are well defined 
in a neighbourhood of $\{x_1=0\}$ and $\{x_2=0\}.$ Then, setting
\begin{align*}
&M_1(u)=\exp\left(\int_0^u\left(\frac{1}{K(0,z)}+\frac{1}{\lambda}\right)\frac{dz}{z}\right)\partial_1\Big(\frac{1}{K}\Big)(0,u)\\
\intertext{and}
&M_2(u)=\exp\left(\int_0^u\big(K(z,0)+\lambda\big)\frac{dz}{z}\right)\partial_2K(u,0),
\end{align*}
we define
\begin{align}\label{F1}
F_{1}(\mu)&=-\int_0^{+\infty} \Big(M_1(z)-M_1(0)
   +\exp(G_{1})\big(M_1(-z)-M_1(0)\big)\Big)\frac{dz}{z^{1+1/\lambda}},\\[3pt]
F_{2}(\mu)&=\int_0^{+\infty} \Big(M_2(-z)-M_2(0)+\exp(G_{2})\big(M_2(z)-M_2(0)\big)\Big)\frac{dz}{z^{1+\lambda}}\label{F2}
   \intertext{and}\label{F3}
F_{3}(\mu)&=-G_2\big(\partial_1K\partial_2K+\partial_{12}K\big)(0,0),   
\end{align}
where
\begin{equation*}
G_{1}=\int_{-1}^1\left(\frac{q_n(1,z)}{\h(1,z)}+1+\frac{1}{\lambda}
        +\frac{zq_n(z,1)}{\h(z,1)}\right)\frac{dz}{z}\\
  \text{ and }
G_{2}=\int_{0}^{{+\infty}}\left(\frac{q(z,0)}{g(z)}+\frac{q(-z,0)}{g(-z)}\right)dz.
\end{equation*}
Taking this notation into account, the functions that determine the cyclicity (and stability) of the polycyle~$\Gamma_u$ at first and second order are the following:
\begin{equation}\label{d0}
d_{0}(\mu)\!:=-\int_{-{\infty}}^{+\infty}\left(\frac{q(z,0)}{g(z)}
    +\lambda \frac{q_n(z,1)}{\h(z,1)}\right)dz
 \;\text{ and }\;
 d_1(\mu)\!:=
 \left\{
  \begin{array}{ll}
   F_1(\mu) & \text{ if $\lambda(\mu)>1,$}\\[5pt]
   F_2(\mu) & \text{ if $\lambda(\mu)<1,$}\\[5pt]
   F_3(\mu) & \text{ if $\lambda(\mu)=1.$}
  \end{array} 
 \right.
\end{equation}
Let us advance that the function $d_0$, together with the functions $F_1,$ $F_2$ and $F_3$ defining $d_1,$ generate the ideal of coefficients at order one and two of the asymptotic expansions of the displacement function studied in \teoc{teo1}.
This result also shows that $d_0$ is analytic on $\Lambda$ and $d_1$ is analytic on 
$\Lambda\setminus\Lambda_1$, where $\Lambda_1\!:=\{\mu\in\Lambda:\lambda(\mu)= 1\}.$ In the next statement $\mathscr R_u(\,\cdot\,;\mu)$ stands for the return map of the vector field~$X_\mu$ around the polycycle~$\Gamma_u$, see \figc{dib3}, and we use the notion of functional independence given in \defic{indepe}.

\begin{bigtheo}\label{thmA}
Let us consider the family of polynomial vector fields $\{X_\mu\}_{\mu\in\Lambda}$ given in \refc{DS} and verifying the assumptions {\bf H1} and {\bf H2}. Then the following assertions hold for any $\mu_0\in\Lambda$ such that $\mathscr R_u(\,\cdot\,;\mu_0)\not\equiv \text{Id}:$

\begin{enumerate}
\item[$(a)$] If $d_0(\mu_0)\neq0$ then 
         $\mathrm{Cycl}\big((\Gamma_u,X_{\mu_0}),X_\mu\big)=0.$
         
\item[$(b)$] If $d_0$ vanishes and is independent at $\mu_0$ then        
         $\mathrm{Cycl}\big((\Gamma_u,X_{\mu_0}),X_\mu\big)\geqslant 1.$
         
\item[$(c)$] If $d_1(\mu_0)\neq 0$ 
         then $\mathrm{Cycl}\big((\Gamma_u,X_{\mu_0}),X_\mu\big)\leqslant 1.$

\item[$(d)$] If $d_0$ and $d_1$ vanish and are independent at $\mu_0$ and 
         $\lambda(\mu_0)\neq 1$ then 
         $\mathrm{Cycl}\big((\Gamma_u,X_{\mu_0}),X_\mu\big)\geqslant 2.$ Moreover the same lower bound 
         holds in case that $\lambda(\mu_0)=1$ and the restrictions
         $\left.d_0\right|_{\Lambda_1}$ and $\left.d_1\right|_{\Lambda_1}$ vanish and are independent at $\mu_0$.
        
\end{enumerate}
\end{bigtheo}
With regard to the application of \teoc{thmA} it is worth noting that if $d_0(\mu_0)\neq 0$, or $d_1(\mu_0)\neq 0$, then $\mathscr R_u(\,\cdot\,;\mu_0)\not\equiv \text{Id}.$ This is a consequence of \teoc{teo1}, which is a fundamental result to prove \teoc{thmA}.

The stability of this kind of hemicycle was previously studied in \cite[Theorem 7]{GMM02}. Indeed, using our notation, the authors prove that $\mathscr R_u(s;\mu)=e^{d_0(\mu)}s+\op(s)$, so that if $d_0(\mu_0)<0$ (respectively, $d_0(\mu_0)>0$) then the polycycle $\Gamma_u$ of the vector field $X_{\mu_0}$ is asymptotically stable (respectively, unstable). In this paper, by performing a second order analysis we also obtain the stability in case that $d_0(\mu_0)=0$ and $d_1(\mu_0)\neq 0$ (see \obsc{estable}). That being said, the goal of the present paper is not to study the stability of the hemicycle but its cyclicity. The first notion concerns single vector fields, whereas the second one is addressed to families of vector fields (i.e., depending on parameters). This is the reason why we need the remainder in the asymptotic expansion of $\mathscr R_u(s;\mu)$ at $s=0$ to be uniform with respect to the parameters. Let us also note that similar results (for both, stability and cyclicity) can be obtained for the hemicycle $\Gamma_\ell$ by performing the change of variables $(x,y)\mapsto (x,-y).$ 

\teoc{thmA} is a general result for the cyclicity of the polycycle $\Gamma_u$ of a D-system $X_{\mu_0}$ with $\mathscr R_u(\,\cdot\,;\mu_0)\not\equiv\text{Id}$ when perturbed inside the family of D-systems~\refc{DS}. Note that in doing so the polycycle $\Gamma_u$ is persistent (i.e., the connections between the two vertices remain unbroken through the perturbation). In contrast, the rest of our main results concern the cyclicity of quadratic D-systems $X_{\mu_0}$ with $\mathscr R_u(\,\cdot\,;\mu_0)\equiv\text{Id}$ when perturbed inside the whole family of quadratic systems. This means in particular that the connection breaks, see \figc{dib5}. More concretely, in Theorems~\ref{thmB}, \ref{thmC} and~\ref{ThmD}, for each $(a_0,b_0)\in (-2,0)\times (0,2)$, we perturb the quadratic D-system 
 \begin{figure}[t]
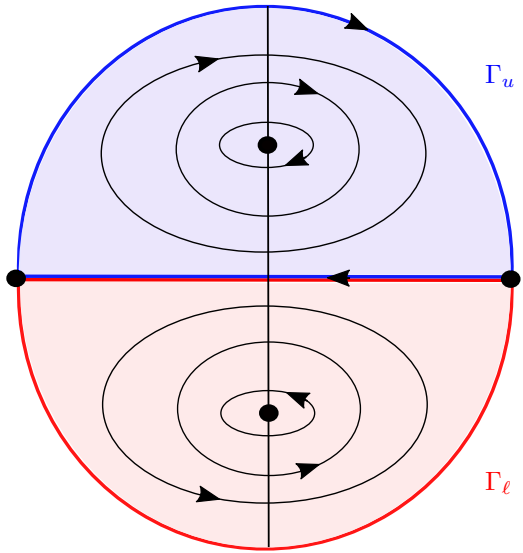

   \centering
  \begin{lpic}[l(0mm),r(0mm),t(0mm),b(5mm)]{hemi0(0.9)}
   \lbl[l]{70,70;$\blue{\Gamma_u}$}          
   \lbl[l]{70,10;$\red{\Gamma_\ell}$}        
   \end{lpic}
  \caption{Phase portrait in the Poincar\'e disc of the quadratic differential system \refc{DSq} for each $(a_0,b_0)\in (-2,0)\times (0,2)$.}\label{dib0}
 \end{figure}
\begin{equation}\label{DSq}
 \sist{\frac{b_0-2}{4}+(1-b_0)y+a_0x^2+b_0y^2,}{-2xy,}
\end{equation}
that one can show it verifies assumptions {\bf H1} and {\bf H2}. Moreover 
it has two centers, located at the points $(0,\frac{1}{2})$ and  $(0,\frac{b_0-2}{2b_0})$ whose period annulus foliate, respectively, the upper and lower half-planes, see 
\figc{dib0}. 

%\begin{bigtheo}\label{thmB}
%If $(a_0,b_0)\in (-2,0)\times (0,2)$ and $a_0\neq -1$ then the cyclicities of $\Gamma_u$ and $\Gamma_\ell$ when we perturb~\refc{DSq} inside the whole family of quadratic differential systems are exactly 2. Furthermore both cyclicities are at least 2 for $(a_0,b_0)\in \{-1\}\times (0,2)$.
%\end{bigtheo}

\begin{bigtheo}\label{thmB}
Let us take any $(a_0,b_0)\in (-2,0)\times (0,2)$. Then
the cyclicity of $\Gamma_u$ when we perturb~\refc{DSq} inside the whole family of quadratic differential systems is exactly 2 if $a_0\neq -1$ and at least 2 if $a_0=-1.$ Moreover the same statement is true for $\Gamma_\ell.$
\end{bigtheo}

We point out that this result does not imply that the number of limit cycles bifurcating simultaneously from $\Gamma_u$ and $\Gamma_\ell$ is four. As a matter of fact this number is at most three by the forthcoming \teoc{thmC}. Using the terminology from \cite{Iliev}, both centers of the unperturbed system \refc{DSq} are inside the reversible component $Q_3^R$ of the \emph{center manifold} of the quadratic systems.
There are three other components: Hamiltonian~$Q_3^H$, codimension four~$Q_4$ and generalized Lotka-Volterra~$Q_3^{LV}$. It turns out (cf. \lemc{centreLV}) that the centers of the unperturbed system belong also to the $Q_3^{LV}$ component in case that $(a_0+b_0)(a_0-b_0+2)=0$, and when this occurs the proof of \teoc{thmB} is a little more difficult.  

\begin{figure}[t]
\begin{center}
\hspace{-13mm}
\begin{tikzpicture}[scale=1.2]
%\draw[step=1mm,gray,very thin] (-6,-1) grid (1,4);
\draw[fill,color=gray]  (-4,2) to [out=-80,in=135]  (-3.3,0.65) to [out=-45,in=180] (-2,0) to (-2/3,0) to (-4,2);
\draw[fill,color=gray]  (-4.8,4) to [out=-60,in=100] (-4,2) to (-6,3.2) to (-6,4) to (-4.8,4);
\draw[fill,color=gray] (-2/3,0) to (0,0) to (1,-1) to (-2/3,0);
\draw[fill,color=gray] (1,-1) to (1.5,-1.3) to (1.5,-1.5) to (1,-1);
\draw[] (0,-1.5) to (0,4);
\draw[] (-6,0) to (1.5,0);
\draw[thick] (1.5,-1.3)  to  (-6,3.2);
\draw[dashed] (-2,0) to (-2,2) to (0,2);
\draw[blue] (-3.5,-1.5) to (1.5,3.5);
\draw[blue] (0,0) to (-4,4);
\draw[thick] (-2,0) to (0,0) to (1.5,-1.5) ;
\draw[thick] (-4.8,4) to [out=-60,in=100] (-4,2) to [out=-80,in=135]  (-3.3,0.65) to [out=-45,in=180] (-2,0);
\draw [fill] (-4,2) circle [radius=.05];
\node at (-3.7,2.1) {$Q_4^+$};
\draw [fill] (-2,0) circle [radius=.05];
\node at (-2,0.2) {$S_2$};
\node at (-2,-0.2) {$-2$};
\node at (0.2,2) {$2$};
\node at (0.15,0.15) {$0$};
\draw [fill] (-4,1) circle [radius=.05];
\node at (-4.2,1) {$S_4$};
\draw [fill] (-1/2,0) circle [radius=.05];
\node at (-1/2,0.2) {$S_3$};
\draw [fill] (-1,1) circle [radius=.05];
\node at (-1,.8) {$S_1$};
\draw [fill] (-2/3,0) circle [radius=.05];
\node at (-2/3,-0.2) {$Q_4^-$};
\node at (-6,2.4) {$3a_0+5b_0+2=0$};
\node at (-4.3,3.5) {$\mc C$};
\begin{scope}[shift={(-1.24,0)},scale=0.25]
\path[fill] (10.5,0.2) -- (10.5,-0.2) -- (11,0) -- cycle;
\end{scope}
\node at (1.4,.18) {$a_0$};
\begin{scope}[shift={(0,1.26)},rotate=90,scale=0.25]
\path[fill] (10.5,0.2) -- (10.5,-0.2) -- (11,0) -- cycle;
\end{scope}
\node at (.2,3.9) {$b_0$};
\draw (-0.05,1) to (0.05,1);
\node at (0.2,1) {$1$};
\draw (-1,-0.05) to (-1,0.05);
\node at (-1.1,-0.2) {$-1$};
\draw (-4,-0.05) to (-4,0.05);
\node at (-4,-.2) {$-4$};
\draw (1,.05) to (1,-0.05);
\node at (1,-.2) {$1$};
\draw (-0.05,-1) to (0.05,-1);
\node at (0.25,-1) {$-1$};
\end{tikzpicture}
\end{center}
\caption{According to Illiev's conjecture, the shaded area corresponds to those parameters $(a_0,b_0)$ for which
the period annulus $\mathscr P$ of the center at $(0,\frac{1}{2})$ of system \refc{DSq} has cyclicity~3. Its boundary has two components: the straight line $3a_0+5b_0+2=0$ and a piecewise curve $\mc C.$ The straight line 
corresponds to parameters for which the center itself has cyclicity 3. The curve $\mc C$ corresponds to parameters for which the polycycle at $\partial\mathscr P$ has cyclicity 3. The parameters $S_1=(-1,1)$, $S_2=(-2,0)$, $S_3=(-\frac{1}2,0)$ and $S_4=(-4,1)$  are the four isochronous quadratic centers. The blue straight lines are the intersection points with the component $Q_3^{LV}$ of the center manifold. The parameters $Q_4^+=(-4,2)$ and $Q_4^-=(-\frac{2}3,0)$ are the intersection points with the component $Q_4$.}\label{FigIliev}
\end{figure}
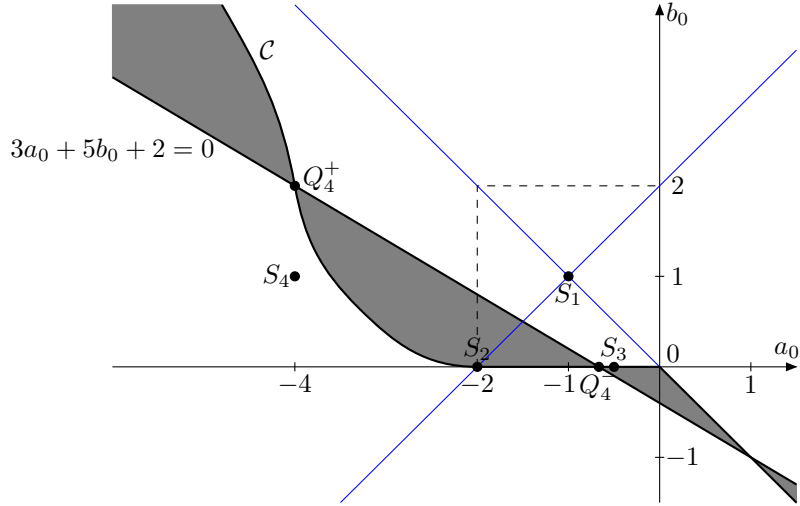

Closely related to \teoc{thmB},  a result due to Swirszcz (see~\cite[Theorem 1]{Swirszcz}) is worth to be quoted. Indeed, in that paper the author also studies the cyclicity of a polycycle of a quadratic reversible system when perturbed inside the whole quadratic family. More concretely, he perturbs the differential system~\refc{DSq} but taking $(a_0,b_0)\in\mc S\!:=\{0<b_0<-a_0\}\cap\{a_0<-2\}$. For these parameters the singular point $(0,\frac{1}{2})$ is also a center but the polycycle at the boundary of its period annulus is not an hemicycle. It is a bicycle~$\Gamma_b$ with the two vertices at infinity, and consisting of a branch of a hyperbola together with a segment of $\ell_\infty$. Recall that the \emph{period annulus} of a center~$p$ is its largest punctured neighbourhood~$\mathscr P$ which is entirely covered by periodic orbits, and that its boundary~$\partial\mathscr P$ has two connected components: the center itself and a polycycle. By using a completely different approach than ours, and with a lower level of detail in the proofs,
Swirszcz identifies a curve $\mc C$ (see \figc{FigIliev}) such that the cyclicity of $\Gamma_b$ is 3 if $(a_0,b_0)\in\mc S\cap\mc C$ and 2 if $(a_0,b_0)\in\mc S\setminus\mc C$. It is to be noted that the only parameter value in $\mc S$ which intersects another center component is $(a_0,b_0)=(-4,2),$ that belongs also to the $Q_4$ component.  

In another vein, Iliev studies in his seminal paper~\cite{Iliev} the cyclicity of the period annulus $\mathscr P$ of the quadratic centers. We stress that the definition of cyclicity for $\mathscr P$ is different than the one for a polycycle because the former is open (see \defic{cycl-G}). Among other results Iliev proves that the cyclicity 
of the period annulus $\mathscr P$ of the center at $(0,\frac{1}{2})$ of the differential system~\refc{DSq} is $3$ for $(a_0,b_0)=(-4,2)$ and $2$ for $(a_0,b_0)=(-1,1)$. These two parameters are denoted, respectively, by $Q_4^+$ and $S_1$ in \figc{FigIliev}. Moreover he conjectures that the cyclicity of $\mathscr P$ is equal to 3 if $(a_0,b_0)$ is inside the shaded area in \figc{FigIliev} and equal to 2 if $(a_0,b_0)$ is outside. 
Previous to Iliev's conjecture, there is a result by Shafer and Zegeling (see~\cite[Theorem 3.2]{ShaZeg94}) that determines some regions where the cyclicity of $\mathscr P$ is equal to 3. They also give a numerical approximation to the curve $\mc C.$ In this setting \teoc{thmB} reinforces Iliev's conjecture because it shows that the curve $\mc C$ does not enter the square $(a_0,b_0)\in (-2,0)\times(0,2)$.
%
%In \cite[Theorems 1.1 and 1.2]{Liu12} it is claimed that the cyclicity of $\mc A$ at $\{a=-3/2\}$ is $3$ if $b\in(0,1/2)$, $2$ if $b\in(1/2,3/2)$ and $1$ if $b\in(3/2,1)$. This last value contradicts Iliev's conjecture and we think that the correct value is $2$ which is the cyclicity of the outer boundary $\Gamma_u$ of $\mc A$  according to our \teoc{thmB}.
%The problem would come from the statement Theorem~1.2 of \cite{Liu12} where the cyclicities of each of the two period annuli and the simultaneous cyclicity of both of them are expressed in a somewhat messy way. 

Let us recall at this point that \emph{Hilbert's 16th problem} asks for the maximum number $H(n)$ of limit cycles of a planar polynomial differential system of degree $\leqslant n.$ It is still open for any $n\geqslant 2.$ In 1994 Dumortier, Roussarie and Rousseau conceived a program (see \cite{DRR2}) to prove that $H(2)$ is finite. In short, they reduced this problem to prove the finite cyclicity for only 121 (different classes of) graphics occurring in quadratic systems. According to the notation in that paper, the quadratic system~\refc{DSq} with $(a_0,b_0)\in (-2,0)\times(0,2)$ is inside the class~$H_2^1$ of hyperbolic hemicycles surrounding a center (see \cite[Figure 7]{DRR2}). Thus,
%\footnote{\blue{La ciclicidad finita del policiclo~$H_2^1$ es consecuencia del Teorema 2 de MorMouRou94. Este p\'arrafo deber\'{\i}a redactarse de nuevo. Curioso. En la p\'agina 313 asumen que pueden poner rectas las separatrices de cada silla de manera anal\'{\i}tica, incluso con respecto de los par\'ametros. Me gustar\'{\i}a saber como.}} 
\teoc{thmB} can be viewed as a small contribution to the completion of the program to prove that $H(2)<\infty.$ Nevertheless some authors (e.g. \cite{Rou}) attribute to Mourtada the proof of the finite cyclicity of any hyperbolic polycycle in an unpublished series of manuscripts (see~\cite[Theorem 0]{Mou5} and references therein). For other results about the cyclicity of quadratic hemicycles in this context the interested reader is referred to \cite{DGR02,RouRou08}.  

Note that \teoc{thmB} provides the cyclicity of $\Gamma_u$ and $\Gamma_\ell$ individually, i.e., taking $\Pi=\{\Gamma_u\}$ and $\Pi=\{\Gamma_\ell\}$ in \defic{defcic}. In our third main result we study the cyclicity of $\Pi=\{\Gamma_u,\Gamma_\ell\}$, cf. \obsc{cicl_sim}, when we perturb \refc{DSq} inside the family of quadratic differential systems. In its statement we use the following parameter subsets:
\color{black}
\begin{align*}
&\mathcal K_1\!:=\{(a_0,b_0)\in (-2,0)\times (0,2): a_0+b_0\leqslant 0\text{ or }a_0-b_0+2\leqslant 0\}
\intertext{and}
&\mathcal K_2\!:=\{(a_0,b_0)\in (-2,0)\times (0,2): a_0+b_0> 0\text{ and }a_0-b_0+2> 0\}.
\end{align*}
\begin{bigtheo}\label{thmC}
If $(a_0,b_0)$ belongs to $\mathcal K_1\setminus\{a_0=-1\}$ $($respectively, $\mathcal K_2)$ then the cyclicity of $\Pi=\{\Gamma_u,\Gamma_\ell\}$ when we perturb \refc{DSq} inside the whole family of quadratic differential systems is exactly 3 $($respectively, 2$)$. Moreover it is at least 3 for $(a_0,b_0)\in\{-1\}\times (0,2).$ 
\end{bigtheo}

We stress that \teoc{thmC} deals with the simultaneous bifurcation of limit cycles from $\Gamma_u$ and $\Gamma_\ell$, which are the outer boundaries of two period annuli. Note that if this simultaneous cyclicity is 3 then, as a consequence of \teoc{thmB}, two limit cycles bifurcate from $\Gamma_u$ and one from $\Gamma_\ell$, or vice versa. The simultaneous bifurcation of limit cycles from the two period annuli has been studied for $a_0=-\frac{3}{2}$ and $a_0=-\frac{1}{2}$ in \cite{Liu12} and \cite{CLP09}, respectively, and also for $(a_0,b_0)=(-\frac{1}2,\frac{1}2)$ and  $(a_0,b_0)=(-1,1)$ in \cite{PFL14} and \cite{FG}, respectively. The authors do not know of any previous work dealing with the simultaneous bifurcation from two polycycles. 

We turn now to the statement of our last main result, \teoc{ThmD}, which deals with alien limit cycles. This notion was introduced by Dumortier and Roussarie in \cite{alien2}, where the authors bring to light that there are limit cycles bifurcating from a polycycle which cannot be detected as a zero of the first Melnikov function (see also~\cite{alien1,alien3,alien8,alien4}). Our aim in this paper with regard to this issue is twofold. On the one hand, to propose a definition of this phenomenon more intrinsic and geometric than the one used in the literature and not depending on the computation of Melnikov functions. On the other hand we want to show that there exist alien limit cycles in the context of \emph{simultaneous bifurcations}. To this end, following Gavrilov \cite{G08} we first introduce the notion of cyclicity of an open subset $U$ as follows. (He considers the case when $U$ is a period annulus and here we extend it slightly.)

\begin{defi}\label{cycl-G}
Let $\{X_\mu\}_{\mu\approx\mu_0}$ be a germ of an analytic family of vector fields on $\Sc^2$ and let  $K$ be a compact subset of $\Sc^2$. We define the cyclicity of $K$ with respect to the germ $\{X_\mu\}_{\mu\approx\mu_0}$ as
\[
\mr{Cycl}_G\big((K,X_{\mu_0}),X_\mu\big)=\inf\limits_{\varepsilon,\delta>0}\sup\limits_{\mu\in B_\delta(\mu_0)}
\#\big\{\gamma\subset N_{\varepsilon}(K)\text{ limit cycle of $X_\mu$}\big\}\in\Z_{\ge 0}\cup\{\infty\},
\]
where $N_\varepsilon(K)$ is the tubular $\varepsilon$-neighbourhood of $K$.
If $U\subset\Sc^2$ is open we define 
\[
\mr{Cycl}_G\big((U,X_{\mu_0}),X_\mu\big)=\sup\left\{\mr{Cycl}_G\big((U,X_{\mu_0}),X_\mu\big):K\subset U\text{ compact}\right\},
\]
which may also be infinite.
\end{defi}
In case that $U$ is a period annulus with finite cyclicity in the above sense, Gavrilov proves in \cite[Theorem~1]{G08} 
that $\mr{Cycl}_G\big((U,X_{\mu_0}),X_\mu\big)$ is the same as in an appropriate one-parameter analytic deformation. 
This is related with the notion of essential perturbation introduced by Illiev \cite{Iliev} and enables to tackle the problem by computing Melnikov functions. This well-known approach allows to bound the number of limit cycles bifurcating from any compact set $K\subset U$ by means of the Weierstrass Preparation Theorem, however it gives not enough information on $U\setminus K.$ This motivates the following definition. 

\begin{defi}\label{inside}
Let $\{X_\mu\}_{\mu\approx\mu_0}$ be a germ of an analytic family of vector fields on $\Sc^2$ and consider an open subset $U$ of~$\Sc^2$. We define the \emph{boundary cyclicity of $U$ from \underline{inside}} as 
\[
\underline{\mr{Cycl}}^{\,U}_{\,G}\big((\partial U,X_{\mu_0}),X_\mu\big)=\inf\left\{\mr{Cycl}_G\big((U\setminus K,X_{\mu_0}),X_\mu\big):K\subset U\text{ compact}\right\}\in\Z_{\ge 0}\cup\{\infty\}.
\]
\end{defi}

If $\partial U$ is a polycycle with a return map which is not the identity then it can be shown by a compactness and continuity argument that $\underline{\mr{Cycl}}^{\,U}_{\,G}\big((\partial U,X_{\mu_0}),X_\mu\big)=0$. On the other hand, we prove in \lemc{lema_inside} that 
\[
 \underline{\mr{Cycl}}^{\,U}_{\,G}\big((\partial U,X_{\mu_0}),X_\mu\big)\leqslant\mr{Cycl}_G\big((\partial U,X_{\mu_0}),X_\mu\big).
 \] 
These two facts lead to the following definition:

\begin{defi}\label{alien_defi}
Let $\{X_\mu\}_{\mu\approx\mu_0}$ be a germ of an analytic family of vector fields on $\Sc^2$ such that $X_{\mu_0}$ is a D-system satisfying hypothesis $\mathbf{H1}$ and $\mathbf{H2}$. Assume additionally that the return maps $\mathscr R_u(\,\cdot\,;\mu_0)$ and $\mathscr R_\ell(\,\cdot\,;\mu_0)$ of the hemicycles $\Gamma_u$ and $\Gamma_\ell$ are both the identity. Taking $U=\R^2\setminus\{y=0\}$, if
\begin{equation*}
\underline{\mr{Cycl}}^{\,U}_{\,G}\big((\partial U,X_{\mu_0}),X_\mu\big)<\mr{Cycl}_G\big((\partial U,X_{\mu_0}),X_\mu\big)
\end{equation*}
then we say that an \emph{alien limit cycle bifurcation} occurs at $\partial U=\Gamma_u\cup\Gamma_\ell$ from inside $U$ for $\{X_\mu\}_{\mu\approx\mu_0}$.
\end{defi}

We have not given the notion of alien limit cycle bifurcation for an arbitrary collection of limit periodic sets because the involved casuistry would make the definition more complicated than it should be. This is already evident for the case of the ``figure eight-loop'' in \figc{dib8}. That being said, we do give the definition of alien limit cycle bifurcation for any unfolding of a polycycle satisfying rather natural hypothesis, which is the case of those 2-saddle cycles studied in \cite{alien1,alien3,alien2,alien8,alien4}. This will be done in \secc{provateoD}, see \defic{robert1}. 
Our definition differs from the one used by Dumortier and Roussarie in \cite{alien2} because we account only for limit cycles which cannot be detected as zeroes of any Melnikov function of \emph{any} order, cf. \lemc{robert2}. 

Under the hypothesis in \defic{alien_defi}, the vertices of $\Gamma_u$ and $\Gamma_\ell$ are hyperbolic saddles. In this case it follows from \lemc{L} that  
\[
 \mr{Cycl}_G\big((\partial U,X_{\mu_0}),X_\mu\big)=\mathrm{Cycl}\big((\{\Gamma_u,\Gamma_\ell\},X_{\mu_0}),X_\mu\big).
\]
Hence \defic{alien_defi} takes into account the \emph{simultaneous} bifurcation of limit cycles from $\Gamma_u$ and $\Gamma_\ell.$ In this regard we obtain the following result about alien bifurcations in the quadratic family:

\begin{bigtheo}\label{ThmD}
If  $(a_0,b_0)\in\{(-1,1),(-\frac{1}2,\frac{1}2),(-\frac{1}2,\frac{3}2)\}$ then an alien limit cycle bifurcation occurs at $\Gamma_u\cup\Gamma_\ell$ when we perturb \refc{DSq} inside the whole family of quadratic differential systems. 
\end{bigtheo}

\begin{figure}[t]
\begin{center}
\includegraphics[width=9cm]{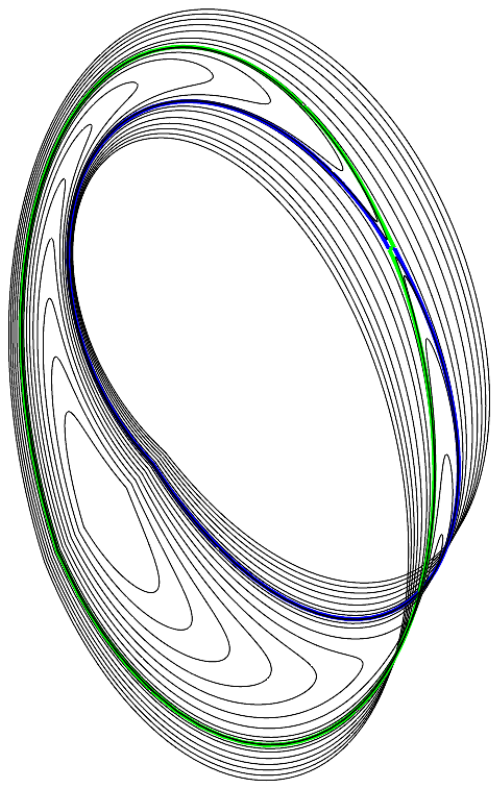}
\end{center}
\caption{Quadratic reversible double centers in \refc{DSq} compactified to the Moebius strip $\mathbb{R P}^2\setminus\mathbb D$. One of the two centers is depicted at the front of the drawing, while we place the other one in the removed invariant disk $\mathbb D$ for convenience. The polycycle $\bar\Gamma_u=\bar\Gamma_\ell$ is represented by the two circles in blue and green intersecting at the saddle point at the back.}\label{mobius}
\end{figure}

Let us remark that in the present paper we consider families of planar polynomial vector fields $\{X_\mu\}_{\mu\in\Lambda}$ and that the statements of our main results should more formally be addressed to the compactified family $\{p(X_\mu)\}_{\mu\in\Lambda}$ of analytic vector fields on the Poincar\'e sphere $\Sc^2.$ For simplicity in the exposition we commit an abuse of language by identifying both families. It is clear that the number of limit cycles of $X_\mu$ and $p(X_\mu)$ is the same because the line at infinity $\ell_\infty$ is invariant in all the cases under consideration. Related with this we note that, although the corresponding analytic extension of the polynomial vector field to~$\Sc^2$ does not descend to the quotient $\mathbb{RP}^2$ of $\Sc^2$ by the central symmetry with respect to the origin, the induced foliation does. Since limit cycles depend on the foliation, and not on the specific way in which the orbits are parametrized, one could consider the notion of cyclicity in the real projective plane $\mathbb{RP}^2$ instead of the sphere~$\Sc^2$. It is worth to point out that these two notions are not equivalent. Indeed, the two hemicycles $\Gamma_u$ and $\Gamma_\ell$ in $\Sc^2$ project to the same polycycle $\bar\Gamma_u=\bar\Gamma_\ell$ on $\mathbb{RP}^2$  (see Figure~\ref{mobius})
and by applying Theorems~\ref{thmB} and~\ref{thmC}, respectively,
\[
\mathrm{Cycl}\big((\Gamma_u,X_{\mu_0}),X_\mu\big)=2
\text{ and }
\mr{Cycl}_{\mathbb{RP}^2}\big((\bar\Gamma_u,X_{\mu_0}),X_\mu\big)=3
\]
for any $(a_0,b_0)\in\mathcal K_1\setminus\{a_0=-1\}$.

The paper is organized as follows. Sections~\ref{provateoA} and~\ref{provateoB} are devoted to prove Theorems~\ref{thmA} and~\ref{thmB}, respectively. Both results strongly rely on the asymptotic development of the difference map $\mathscr D(s;\mu)$ given in \teoc{teo1}. This is a rather technical result that follows by applying the tools developed in \cite{MV19,MV20,MV21} to study the Dulac map and its proof is deferred to Appendix~\ref{B} for reader's convenience. Another important ingredient in the proof of \teoc{thmB} is \teoc{division}, which provides a very useful division of the difference map in the ideal generated by its coefficients. The proofs of Theorems~\ref{thmC} and~\ref{ThmD} are given in sections~\ref{provateoC} and~\ref{provateoD}, respectively. Appendix~\ref{A} gathers the essential definitions and results from \cite{MV19,MV20,MV21} that we use in the present paper, together with some other auxiliary results. Finally, in Appendix~\ref{B} we demonstrate \teoc{teo1} and \propc{prop2}, which have the longest and most technical proofs. 

\section{Proof of \teoc{thmA}}\label{provateoA}

\begin{figure}[t]
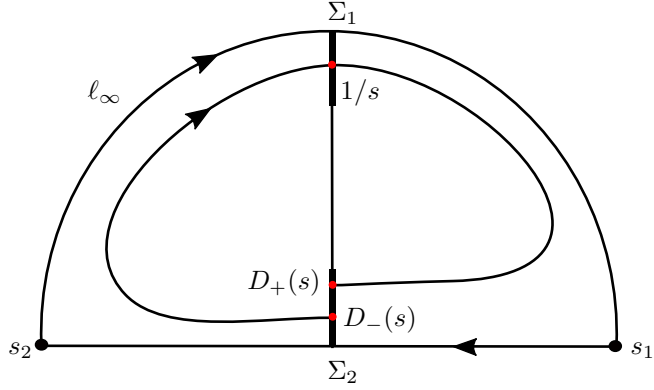

   \centering
  % \vspace{5mm}
  \begin{lpic}[l(0mm),r(0mm),t(0mm),b(5mm)]{dib2b(0.7)}
    \lbl[l]{55,-3;$\Sigma_2$}
    \lbl[l]{55,65;$\Sigma_1$}      
    \lbl[l]{10,50;$\ell_\infty$}           
    \lbl[l]{112,1;$s_1$}
    \lbl[l]{-5,1;$s_2$}               
    \lbl[l]{58,7;${D_-(s)}$}      
    \lbl[l]{57.5,50;${1/s}$}      
    \lbl[l]{40,14;${D_+(s)}$}    
   \end{lpic}
  \caption{
  Dulac maps for the definition of $\mathscr D=D_+-D_-$ in \teoc{teo1}}\label{dib2}
 \end{figure}

In this section we consider the family of vector fields $\{X_\mu\}_{\mu\in\Lambda}$ given by \refc{DS}
and satisfying the hypothesis {\bf H1} and {\bf H2}. We take two local transverse sections, $\Sigma_1$ and $\Sigma_2$ parametrised, respectively, by $s\mapsto (0,\frac{1}{s})$ and $s\mapsto (0,s)$ with $s>0.$ We also define $D_+(s;\mu)$ to be the Dulac map of $X_\mu$ from $\Sigma_1$ to $\Sigma_2$ and $D_-(s;\mu)$ to be the Dulac map of $-X_\mu$ from $\Sigma_1$ to $\Sigma_2,$ see \figc{dib2}. The limit cycles of $X_\mu$ that are close to 
$\Gamma_u$ in Hausdorff sense are in one to one correspondence with the isolated positive zeroes of the \emph{difference map}
 \[
  \mathscr D(s;\mu)\!:=D_+(s;\mu)-D_-(s;\mu)
 \]
near $s=0.$ The following result gives the asymptotic development of $\mathscr D(s;\mu)$ at $s=0$ and the functions $\lambda,$ $F_1$, $F_2$, $F_3$ and $d_0$ in its statement are the ones defined in \refc{K}, \refc{F1}, \refc{F2}, \refc{F3} and \refc{d0}, respectively. In the statement we use the Ecalle-Roussarie comensator $\omega(s;\alpha)$, see \defic{defi_comp}, and $\F_{\ell}^\infty(\mu_0)$ stands for a function $\ell$-flat with respect to $s$ at $\mu_0$, see \defic{defi2}. 

\begin{theo}\label{teo1} 
Let us fix any $\mu_0\in\Lambda$ and set $\lambda_0\!:=\lambda(\mu_0).$ Then $\mathscr D(s;{\mu})=\Delta_{0}({\mu})s^\lambda +\F_{\ell}^\infty(\mu_0)$ for any ${\ell}\in \big[\lambda_0,\min(2\lambda_0,\lambda_0+1)\big),$ where $\Delta_0$ is an analytic function at $\mu_0$ that can be written as $\Delta_{0}=\kappa_0d_0$, with $\kappa_0$ analytic at~$\mu_0$ and $\kappa_0(\mu_0)>0.$ In addition,

\begin{enumerate}[$(1)$]

\item If $\lambda_0>1$ then $\mathscr D(s;{\mu})=\Delta_{0}({\mu})s^\lambda +\Delta_{1}({\mu})s^{\lambda+1} +\F_{\ell}^\infty(\mu_0)$ for any ${\ell}\in \big[\lambda_0+1,\min(2\lambda_0,\lambda_0+2)\big)$. Furthermore $\Delta_{1}$ is an analytic function at $\mu_0$ that can be written as $\Delta_{1}=\kappa_1F_1+\bar\kappa_1\Delta_0$, where $\kappa_1$ and $\bar\kappa_1$ are analytic at $\mu_0$ and $\kappa_1(\mu_0)>0.$

\item If $\lambda_0<1$ then $\mathscr D(s;{\mu})=\Delta_{0}({\mu})s^\lambda +\Delta_{2}({\mu})s^{2\lambda}+\F_{\ell}^\infty(\mu_0)$ for any ${\ell}\in\big[2\lambda_0,\min(3\lambda_0,\lambda_0+1)\big)$. Moreover $\Delta_{2}$ is an analytic function at $\mu_0$ that can be written as $\Delta_{2}=\kappa_2F_2+\bar\kappa_2\Delta_0$, where $\kappa_2$ and $\bar\kappa_4$ are analytic at $\mu_0$ and $\kappa_2(\mu_0)>0.$ 

\item If $\lambda_0=1$ then $\mathscr D(s;{\mu})=\Delta_{0}({\mu})s^\lambda +\Delta_{3}({\mu})s^{\lambda+1}\omega(s;1-\lambda)+\Delta_4(\mu)s^{\lambda+1}+\F_{\ell}^\infty(\mu_0)$ for any ${\ell}\in [2,3)$
and where $\Delta_{3}$ and $\Delta_{4}$ are analytic functions at $\mu_0$. Moreover there exist analytic functions $\kappa_3$ and $\bar\kappa_3$ at $\mu_0$ with $\kappa_3(\mu_0)>0$ such that the equality $\Delta_{3}=\kappa_3F_3+\bar\kappa_3\Delta_0$ holds on $\{\mu\in\Lambda: \lambda(\mu)=1\}.$

\end{enumerate}
\end{theo}
 
 Since the proof of \teoc{teo1} is rather long and technical and also requires several results from previous papers, we postpone it to Appendix~\ref{B} for reader's convenience.
 
\begin{obs}\label{estable}
On account of the definition of $d_1$ given in \refc{d0}, \teoc{teo1} provides the following information about the stability of the polycycle $\Gamma_u$ for the vector field $X_{\mu_0}$:
\begin{enumerate}[$(a)$]
\item If $d_0(\mu_0)<0$ (respectively, $d_0(\mu_0)>0$) then $\Gamma_u$ is asymptotically stable (respectively, unstable).
\item If $d_0(\mu_0)=0$ and $d_1(\mu_0)<0$ (respectively, $d_1(\mu_0)>0$) then $\Gamma_u$ is asymptotically  stable (respectively, unstable).
\end{enumerate}
The key point for this observation is that the functions $\kappa_i$ 
in the statement of \teoc{teo1} are strictly positive at $\mu_0.$ 
\end{obs}

For simplicity in the exposition, from now on we will use the following definition.

\begin{defi}\label{Z_0}
Let $h(s;\mu)$ be a function in $\cc^{\infty}_{s>0}(U)$ for some open set $U\subset\R^N.$ Given any $\mu_0\in U$ we define $\mathcal Z_0(h(\,\cdot\,;\mu),\mu_0)$ to be the smallest integer $\ell$ having the property that there exist $\delta>0$ and a neighbourhood~$V$ of $\mu_0$ such that for every $\mu\in V$ the function $h(s;\mu)$ has no more than $\ell$ isolated zeros on $(0,\delta)$ counted with multiplicities. 
\end{defi}

\begin{prooftext}{Proof of \teoc{thmA}.} 
Recall (see \figc{dib2}) that the limit cycles of the vector field~\refc{DS} that are close to~$\Gamma_u$ in Hausdorff sense are in one to one correspondence with the isolated positive zeroes of the difference map
 \[
  \mathscr D(s;\mu)=D_+(s;\mu)-D_-(s;\mu)
 \]
near $s=0.$ Hence, see \defic{Z_0}, we have that $\mathrm{Cycl}\big((\Gamma_u,X_{\mu_0}),X_\mu\big)\leqslant\mathcal Z_0\big(\mathscr D(\,\cdot\,;\mu),\mu_0\big)$. Note moreover that, by \teoc{teo1}, 
\begin{equation}\label{teoAeq1}
 \mathscr D(s;{\mu})=\Delta_{0}({\mu}) s^\lambda+\F_{\ell}^\infty(\mu_0)
\end{equation}
for any ${\ell}\in \big[\lambda_0,\min(2\lambda_0,\lambda_0+1)\big),$ where $\lambda_0\!:=\lambda(\mu_0)$ and $\Delta_{0}=\kappa_0d_0$ with $\kappa_0(\mu_0)>0.$ If $d_0(\mu_0)\neq 0$ then, taking any $\ell>\lambda_0$ (see \defic{defi2}), 
\[
 \lim_{(s,\mu)\to (0^+,\mu_0)}s^{-\lambda}\mathscr D(s;{\mu})=\Delta_0(\mu_0)\neq 0,
\] 
which implies $\mathcal Z_0\big(\mathscr D(\,\cdot\,;\mu),\mu_0\big)=0$ and proves $(a)$. 

On the other hand, since $\mathscr D(\,\cdot\,;\mu_0)\equiv 0$ if, and only if, $\mathscr R_u(\,\cdot\,;\mu_0)\equiv\text{Id},$
the assertion in $(b)$ follows from the equality in \refc{teoAeq1} by applying \propc{cota_inf} with $n=1.$ 

We turn next to the proof of $(c)$ and $(d)$. To this end we shall use that, by applying \teoc{teo1},
\begin{equation}\label{teoAeq3}
 \mathscr D(s;\mu)=\Delta_{0}({\mu})s^\lambda+
\left\{
\begin{array}{ll}
\Delta_{1}({\mu})s^{\lambda+1}+ \F_{\ell_1}^\infty(\mu_0)&\text{ if $\lambda_0>1,$}\\[10pt]
\Delta_{2}({\mu})s^{2\lambda}+ \F_{\ell_2}^\infty(\mu_0)&\text{ if $\lambda_0<1,$}\\[10pt]
\Delta_{3}({\mu})s^{\lambda+1}\omega(s;1-\lambda)+\Delta_4(\mu)s^{\lambda+1}+ \F_{\ell_3}^\infty(\mu_0)&\text{ if $\lambda_0=1,$}
\end{array}
\right.
\end{equation}
for any ${\ell_1}\in \big[\lambda_0+1,\min(2\lambda_0,\lambda_0+2)\big)$, ${\ell_2}\in  \big[2\lambda_0,\min(3\lambda_0,\lambda_0+1)\big)$ and $\ell_3\in [2,3)$, respectively. Moreover, in its respective case, the coefficient $\Delta_i$ is an analytic function at $\mu_0$. In addition, for $i\in\{0,1,2,3\},$ there exist analytic functions $\kappa_i$ and $\bar\kappa_i$ at $\mu_0$ with $\kappa_i(\mu_0)>0$ such that we can write
\begin{equation}\label{teoAeq5}
 \Delta_0=\kappa_0d_0,\quad\Delta_1=\kappa_1 F_1+\bar\kappa_1\Delta_0,\quad\Delta_2=\kappa_2 F_2+\bar\kappa_2\Delta_0\,\text{ and }\,\left.\Delta_3\right|_{\Lambda_1}=\left.(\kappa_3 F_3+\bar\kappa_3\Delta_0)\right|_{\Lambda_1}
\end{equation}
where recall that $\Lambda_1\!:=\{\mu\in\Lambda:\lambda(\mu)=1\}.$

In order to show $(c)$ we can suppose that $\Delta_0=\kappa_0d_0$ vanishes at $\mu_0$ because otherwise we have already proved that $\mathrm{Cycl}\big((\Gamma_u,X_{\mu_0}),X_\mu\big)=0.$ On account of this the assumption $d_1(\mu_0)\neq 0$ implies, see the definition given in \refc{d0}, that $\Delta_1(\mu_0)\neq 0$ if $\lambda_0>1,$ $\Delta_2(\mu_0)\neq 0$ if $\lambda_0<1$ and $\Delta_3(\mu_0)\neq 0$ if $\lambda_0=1.$ In the first case, from \refc{teoAeq3} and applying \lemc{FLK},
\begin{align*}
 \partial_s\big(s^{-\lambda}\mathscr D(s;\mu)\big)&=\partial_s\big(\Delta_0(\mu)+\Delta_1(\mu)s+s^{-\lambda}\F_{\ell_1}^\infty(\mu_0)\big)\\&=\Delta_1(\mu)-\lambda s^{-\lambda-1}\F_{\ell_1}^\infty(\mu_0)+s^{-\lambda}\F_{\ell_1-1}^\infty(\mu_0)\\
 &=\Delta_1(\mu)+\F_{\varepsilon}^\infty(\mu_0)
\end{align*}
for some $\varepsilon>0$ small enough since we can take $\ell_1>\lambda_0+1$. Therefore, see \defic{defi2}, the derivative $\partial_s\big(s^{-\lambda}\mathscr D(s;\mu)\big)$ tends to $\Delta_1(\mu_0)\neq 0$ as $(s,\mu)\to (0^+,\mu_0).$ Thus, by applying Rolle's Theorem,  
\[
\mathrm{Cycl}\big((\Gamma_u,X_{\mu_0}),X_\mu\big)\leqslant
\mathcal Z_0\big(\mathscr D(\,\cdot\,;\mu),\mu_0\big)\leqslant 1,
\]  
as desired. Similarly, in the second case (i.e., $\lambda_0<1$) we have that 
\[
\partial_s\big(s^{-\lambda}\mathscr D(s;\mu)\big)=\lambda\Delta_2(\mu)s^{\lambda-1}+\F_{\varepsilon}^\infty(\mu_0)
=s^{\lambda-1}\big(\lambda\Delta_2(\mu)+s^{1-\lambda}\F_{\varepsilon}^\infty(\mu_0)\big)
\]
for some $\varepsilon>0$ small enough. Then, due to $\Delta_2(\mu_0)\neq 0$, we conclude by Rolle's Theorem as before that $\mathcal Z_0\big(\mathscr D(\,\cdot\,;\mu),\mu_0\big)\leqslant 1$. If $\lambda_0=1$ then, from \refc{teoAeq3} once again and taking $\ell_3\in [2,3)$ into account, the application of \lemc{FLK} yields
\begin{align*}
 \partial_s\big(s^{-\lambda}\mathscr D(s;\mu)\big)&=
 \partial_s\big(
 \Delta_0(\mu)+\Delta_3(\mu)s\omega(s;1-\lambda)+\Delta_4(\mu)s+s^{-\lambda}\F_{\ell_3}^\infty(\mu_0)\big)\\
 &=\Delta_3(\mu)\lambda\omega(s;1-\lambda)+\Delta_4(\mu)-\Delta_3(\mu)+\F_{\varepsilon}^\infty(\mu_0) 
\end{align*}
for $\varepsilon>0$ small enough.
Here we use that $\partial_ss\omega(s;\alpha)=(1-\alpha)\omega(s;\alpha)-1,$ see \defic{defi_comp}. Consequently, after dividing the above asymptotic expansion by its leading monomial, one can show that if $(s,\mu)\to (0^+,\mu_0)$ then
\[
 \frac{\partial_s\big(s^{-\lambda}\mathscr D(s;\mu)\big)}{\omega(s;1-\lambda)}=
\lambda \Delta_3(\mu)+\frac{\Delta_4(\mu)-\Delta_3(\mu)}{\omega(s;1-\lambda)}+\frac{\F_{\varepsilon}^\infty(\mu_0)}{\omega(s;1-\lambda)}\to \lambda_0\Delta_3(\mu_0)\neq 0,
\]
since $\lim_{(s,\alpha)\to (0^+,0)}\frac{1}{\omega(s;\alpha)}=0$ by $(a)$ in \cite[Lemma A.4]{MV19}. By Rolle's Theorem again, this implies that $\mathcal Z_0\big(\mathscr D(\,\cdot\,;\mu),\mu_0\big)\leqslant 1$ in the case $\lambda_0=1$ as well and completes the proof of assertion $(c)$. 

Let us show finally the validity of the two assertions in $(d)$. The first one concerns the case $\mu_0\notin\Lambda_1,$ i.e., $\lambda_0\neq 1.$ If $\lambda_0>1$ then, from \refc{teoAeq3}, 
\begin{align*}
s^{-\lambda}\mathscr D(s;\mu)&=\Delta_{0}(\mu)+\Delta_{1}(\mu)s+ f_2(s;\mu)\\
&=\kappa_0d_0+(\kappa_1F_1+\bar\kappa_1\kappa_0d_0)s+f_2(s;\mu)\\
&=d_0\kappa_0(1+\bar\kappa_1s)+d_1\kappa_1s+f_2(s;\mu),
\end{align*}
where in the first equality $f_2\in s^{-\lambda}\F_{\ell_1}^\infty(\mu_0)\subset\F_{1+\varepsilon}^\infty(\mu_0)$ for $\varepsilon>0$ small enough by \lemc{FLK} due to $\ell_1>\lambda_0+1$, in the second one we take \refc{teoAeq5} into account, and in the third one that $d_1(\mu)=F_1(\mu)$ if $\lambda(\mu)>1.$ 
Thus, setting $f_0(s;\mu)=\kappa_0(1+\bar\kappa_1s)$ and $f_1(s;\mu)=\kappa_1s$, we can write
\begin{equation}\label{teoAeq4}
s^{-\lambda}\mathscr D(s;\mu)=d_0(\mu)f_0(s;\mu)+d_1(\mu)f_1(s;\mu)+f_2(s;\mu).
\end{equation}
By assumption we have that $d_0$ and $d_1$ vanish and are independent at $\mu_0$ and that $\mathscr D(\,\cdot\,;\mu_0)\not\equiv 0$ due to $\mathscr R_u(\,\cdot\,;\mu_0)\not\equiv\text{Id}$. Accordingly, since $\frac{f_1(s;\mu)}{f_0(s;\mu)}=\frac{\kappa_1s}{\kappa_0(1+\bar\kappa_1s)}$ and $\frac{f_2(s;\mu)}{f_1(s;\mu)}\in s^{-1}\F_{1+\varepsilon}^\infty(\mu_0)$ tend to zero as $s\to 0^+$, we can apply \propc{cota_inf} with $n=2$ to conclude that 
$\mathrm{Cycl}\big((\Gamma_u,X_{\mu_0}),X_\mu\big)\geqslant 2.$ If $\lambda_0<1$ then following verbatim from  \refc{teoAeq3} and \refc{teoAeq5} we get the equality in \refc{teoAeq4} with $f_0(s;\mu)=\kappa_0(1+\bar\kappa_2s^\lambda)$, $f_1(s;\mu)=\kappa_2s^\lambda$ and $f_2\in s^{-\lambda}\F_{\ell_2}^\infty(\mu_0)\subset\F_{\lambda_0+\varepsilon}^\infty(\mu_0)$. Thus the assumptions in \propc{cota_inf} are also verified and so the lower bound $\mathrm{Cycl}\big((\Gamma_u,X_{\mu_0}),X_\mu\big)\geqslant 2$ is true for the case
$\lambda_0>1$ as well. Let us consider finally the case $\lambda_0=1,$ which is slightly different. In this case, from \refc{teoAeq3} and taking \defic{defi_comp} into account, if $\mu\in\Lambda_1$ then
\begin{align*} 
s^{-\lambda}\mathscr D(s;\mu)&=\Delta_{0}(\mu)-\Delta_3(\mu) s\log s+\Delta_4(\mu) s+\hat f_2(s;\mu)\\
&=d_0\kappa_0(1-\bar\kappa_3s\log s)-F_3\kappa_3s\log s+\Delta_4s+\hat f_2(s;\mu),
\end{align*}
where in the first equality $\hat f_2\in s^{-1}\F_{\ell_3}^\infty(\mu_0)\subset\F_{1+\varepsilon}^\infty(\mu_0)$ and the second one follows from \refc{teoAeq5} due to $\mu\in\Lambda_1$. Hence, since $d_1=F_3$ on $\Lambda_1,$ we can write
\[
\left.s^{-\lambda}\mathscr D(s;\mu)\right|_{\mu\in\Lambda_1}=\left.d_0\right|_{\Lambda_1}f_0(s;\mu)
+\left.d_1\right|_{\Lambda_1}f_1(s;\mu)+f_2(s;\mu)
\]
taking the functions $f_0(s;\mu)=\kappa_0(1-\bar\kappa_3s\log s),$ $f_1(s;\mu)=-\kappa_3s\log s$ and $f_2(s;\mu)=\Delta_4s+\hat f(s;\mu).$ Once again, $\frac{f_1(s;\mu)}{f_0(s;\mu)}=-\frac{\kappa_3s\log s}{\kappa_0(1-\bar\kappa_3s\log s)}$ and 
$\frac{f_2(s;\mu)}{f_1(s;\mu)}=-\frac{\Delta_4+s^{-1}\hat f(s;\mu)}{\kappa_3\log s}$ tend to zero as $s\to 0^+$ and, on the other hand,
$\left.d_1\right|_{\Lambda_1}$ and $\left.d_1\right|_{\Lambda_1}$ vanish and are independent at $\mu_0$ by assumption. Consequently, by applying \propc{cota_inf} with $W=\Lambda_1$ and $n=2$ we get that $\mathrm{Cycl}\big((\Gamma_u,X_{\mu_0}),X_\mu\big)\geqslant 2$ in case that $\lambda_0=1,$ as desired. This proves the second assertion in $(d)$ and concludes the proof of the result. 
\end{prooftext}

\section{Proof of \teoc{thmB}}\label{provateoB}

The following result shows that to prove \teoc{thmB} it suffices to consider a 5-dimensional perturbation.

\begin{lem}\label{5param}
Any quadratic differential system which is close $($in the topology of coefficients$)$ to \refc{DSq} for some $(a_0,b_0)\in\R^2$ with $a_0\neq -2$ can be brought by means of an affine change of coordinates and a constant rescaling of time to 
\begin{equation}\label{pert}
X_\mu\quad\sist{\frac{b-2}{4}+\varepsilon_1x+(1-b)y+ax^2+\varepsilon_2xy+by^2,}{\varepsilon_0-2xy,}
\end{equation}
with $(a,b,\varepsilon_0,\varepsilon_1,\varepsilon_2)\approx (a_0,b_0,0,0,0).$
\end{lem}

\begin{prova}
We consider the group $\mathrm{Aff}(2,\R)$ of affine transformations 
\[
g(x,y)=(g_{11}x+g_{12}y+g_{13},g_{21}x+g_{22}y+g_{23})
\]
and the pull-back $g^\star(Y_{a,b})=\big(Dg^{-1}\big)(Y_{a,b}\circ g)$ of 
\[
Y_{a,b}\!:=\left((1-b)y+by^2+\frac{b-2}{4}+ax^2\right)\!\partial_x-2xy\partial_y.
\]
Note that $Y_{a,b}=w_0+aw_1+bw_2$ with $w_0\!:=(y-\frac{1}{2})\partial_x-2xy\partial_y$, $w_1\!:=x^2\partial_x$ and $w_2\!:=(-y+y^2+\frac{1}{4})\partial_x$.
An easy computation performed with \texttt{Maple} shows that if $a_0\neq -2$ then the vector fields $v_0=\partial_y$, $v_1=x\partial_x$ and $v_2=xy\partial_x$ span a complementary to the tangent space at the point $(\lambda,g,a,b)=(1,\mathrm{id},a_0,b_0)$ of the orbit
\[
\{\lambda g^\star(Y_{a,b}):\,\lambda\in\R^*,\ g\in\mathrm{Aff}(2,\R),\ a,b\in\R\}
\]
in the $12$-dimensional space $\mathcal P_2$  of all polynomial vector fields of degree 2.
In other words, if $a_0\neq -2$ then the map $F:U\!:=\R^*\!\times\mathrm{Aff}(2,\R)\times\R^5\to\mathcal P_2$ defined by
\[
F(\lambda,g,a,b,\varepsilon_0,\varepsilon_1,\varepsilon_2)=\lambda g^\star \big(Y_{a,b}+\varepsilon_0v_0+\varepsilon_1 v_1+\varepsilon_2 v_2\big)
\]
is a local diffeomorphism between neighbourhoods of $(1,\mathrm{id},a_0,b_0,0,0,0)$ in $U$ and $Y_{a_0,b_0}$ in $\mathcal P_2$. This proves the result.
%\begin{verbatim}
%z:=[[0,0],[1,0],[0,1],[2,0],[1,1],[0,2]]:
%v:=[add(a[op(i)]*x^i[1]*y^i[2],i in z),add(b[op(i)]*x^i[1]*y^i[2],i in z)];
%f:=w->[seq(coeff(coeff(w[1],x,i[1]),y,i[2]), i in z),
%seq(coeff(coeff(w[2],x,i[1]),y,i[2]), i in z)];f(v);
%
%G:=w->[lambda*(g[2,2]*eval(w[1],{x=g[1,1]*x+g[1,2]*y+g[1,3],
%y=g[2,1]*x+g[2,2]*y+g[2,3]})-g[1,2]*eval(w[2],{x=g[1,1]*x+g[1,2]*y+g[1,3]
%y=g[2,1]*x+g[2,2]*y+g[2,3]})),
%lambda*(-g[2,1]*eval(w[1],{x=g[1,1]*x+g[1,2]*y+g[1,3],
%y=g[2,1]*x+g[2,2]*y+g[2,3]})+g[1,1]*eval(w[2],
%{x=g[1,1]*x+g[1,2]*y+g[1,3],y=g[2,1]*x+g[2,2]*y+g[2,3]}))]:
%
%glist:=[g[1,1],g[1,2],g[1,3],g[2,1],g[2,2],g[2,3],lambda,a,b]:
%id:={g[1,1]=1,g[1,2]=0,g[1,3]=0,g[2,1]=0,g[2,2]=1,g[2,3]=0,lambda=1}:
%A:=[seq(eval(diff(f(G([y*((1-b)+b*y)+(b-2)/4+a*x^2,-2*x*y])),gg),id), gg in glist),
%f([0,1]),f([x,0]),f([x*y,0])];
%\end{verbatim}
%The determinant is $-16(a+2)\neq 0$ if $a\neq -2$.
%This shows that the tangent space to the slice $(\varepsilon_0,\varepsilon_1,\varepsilon_2)$ is complementary to the tangent space of the set of orbits by the affine group of the Loud family.
\end{prova}

We stress that henceforth $X_\mu$ refers to the differential system in \refc{pert}. 
That being said, the key point for our purposes is that $X_\mu$ writes as 
\[
\sist{yf(x,y;\mu)+g(x;\mu),}{\varepsilon_0+yq(x,y;\mu),}
\]
with $f(x,y)=1-b+\varepsilon_2 x+by,$ $g(x)=\frac{b-2}{4}+\varepsilon_1x+ax^2$ and $q(x,y)=-2x,$ so that $X_\mu$  is a D-system for $\varepsilon_0=0$. Moreover one can easily check that $X_\mu$ with $a\in (-2,0)$, $b\in (0,2)$, $\varepsilon_0=0$, $\varepsilon_1\approx 0$ and $\varepsilon_2\approx 0$ verifies assumptions \textbf{H1} and \textbf{H2}. Accordingly, for these parameter values, $X_\mu$ has a polycycle $\Gamma_u$ at the boundary of the upper half-plane with two hyperbolic saddles, $s_1=\{y=0,x>0\}\cap \ell_\infty$ and $s_2=\{y=0,x<0\}\cap \ell_\infty$. Since~$\varepsilon_0$ does not affect the homogenous part of higher degree of $X_\mu$, the location and character of these two singular points remains unaltered taking $\varepsilon_i\approx 0$ for $i=0,1,2.$

\begin{figure}[t]
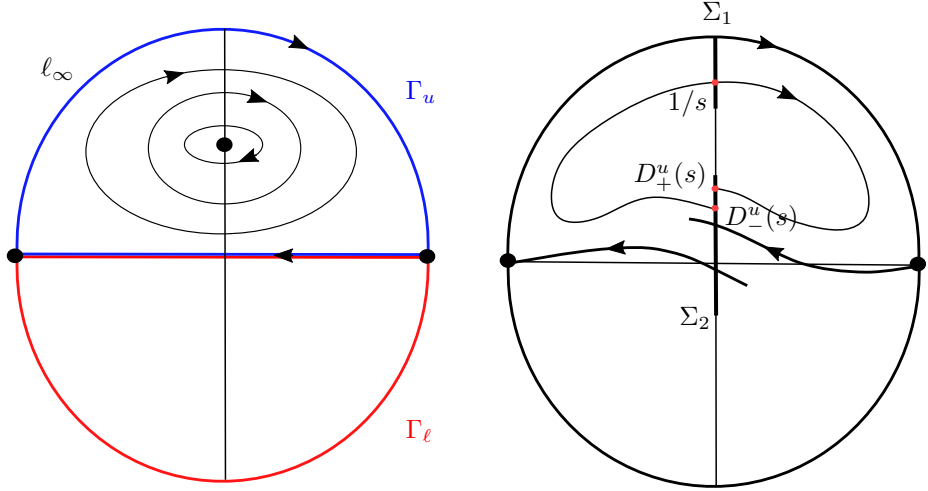

   \centering
  \begin{lpic}[l(0mm),r(0mm),t(0mm),b(5mm)]{dib4(0.75)}
    \lbl[l]{6,74;$\ell_\infty$}       
    \lbl[l]{110,55;$D^u_+(s)$} 
    \lbl[l]{126,47.5;$D^u_-(s)$}  
    \lbl[l]{116,68;$1/s$}      
    \lbl[l]{118,30;$\Sigma_2$}  
    \lbl[l]{122,84;$\Sigma_1$}                  
    \lbl[l]{70,70;$\blue{\Gamma_u}$}           
    \lbl[l]{70,10;$\red{\Gamma_\ell}$}              
   \end{lpic}
  \caption{Phase portrait in the Poincar\'e disc of the vector field $X_\mu$ in~\refc{pert} for $\varepsilon_0=\varepsilon_1=\varepsilon_2=0$ (left) and $\varepsilon_0\neq 0$ (right). On the right, Dulac maps $D_\pm$ to define the function $\mathscr D_u(s;\mu)=D^u_+(s;\mu)-D^u_-(s;\mu)$ studied in \propc{prop2}. The points in red are $(0,D^u_\pm(s))$ and $(0,1/s)$.}
  \label{dib5}
 \end{figure}

Let us fix any $\mu_0=(a_0,b_0,\varepsilon_0,\varepsilon_1,\varepsilon_2)$ with $(a_0,b_0)\in (-2,0)\times (0,2)$ and $\varepsilon_i\approx 0$ for $i=0,1,2.$ We take two transverse sections on $x=0$: $\Sigma_1$, parametrized by $s\mapsto (0,1/s)$ with $s\in (0,\delta),$ and $\Sigma_2$, parametrized by $s\mapsto (0,s)$ with $s\in (-\delta,\delta).$ 
For $\mu\approx\mu_0$ and $\delta>0$ small enough, we have a well defined Dulac map $D^u_+(\,\cdot\,;\mu)$ for $X_\mu$ from $\Sigma_1$ to $\Sigma_2$ and a well defined Dulac map $D^u_-(\,\cdot\,;\mu)$ for $-X_\mu$ from $\Sigma_1$ to $\Sigma_2,$ see \figc{dib5}. This follows by first applying the local center-stable manifold theorem (see \cite[Theorem 1]{Kelley} for instance) and then appealing to the smooth dependence of the solutions of $X_\mu$ on initial conditions and parameters.
In our next result we study the asymptotic development of the difference map
\[
 \mathscr D_u(s;\mu)\!:=D^u_+(s;\mu)-D^u_-(s;\mu).
\]
It is clear that the positive zeros of this function are in one-to-one correspondence with the limit cycles of~$X_\mu$ bifurcating from $\Gamma_u$ to the upper half-plane.

\begin{prop}\label{prop2} 
Fix any $\mu_0=(a_0,b_0,0,0,0)$ with $(a_0,b_0)\in (-2,0)\times (0,2)$. Then 
\[
 \mathscr D_u(s;{\mu})=\delta_{u}+ \Delta^{u}_{0} s^\lambda+ \F_{{L}}^\infty(\mu_0),
 \text{ for any ${{L}}\in \big[\lambda_0,\min(2\lambda_0,\lambda_0+1)\big)$,}
\]
where $\lambda$, $\delta_{u}$ and $\Delta^{u}_0$ are smooth functions in a neighbourhood of $\mu_0$ and $\lambda_0\!:=\lambda(\mu_0)=-\frac{a_0+2}{a_0}.$ In addition $\mathscr D_u(s;{\mu_{0}})\equiv 0$, 
$\partial_{\varepsilon_0}\delta_{u}(\mu_{0})>0$, $\partial_{\varepsilon_1}\delta_{u}(\mu_{0})=\partial_{\varepsilon_2}\delta_{u}(\mu_{0})=0$ and
\[
\textstyle\Delta^{u}_0(\mu)=-\kappa_{01}(\mu)\left(2\frac{\sqrt{b(a+2)}}{\sqrt{a(b-2)}}\,\varepsilon_1+\varepsilon_2\right)
 +\kappa_{02}(\mu)\delta_{u}(\mu),
\] 
where $\kappa_{0i}$ are smooth functions at $\mu=\mu_0$ for $i=1,2$ and $\kappa_{01}(\mu_0)>0.$
Furthermore the following assertions are also true in case that $a_0\neq -1:$
\begin{enumerate}[$(1)$]

\item If $a_0>-1$ then $\mathscr D_u(s;{\mu})=\delta_u+ \Delta^{u}_{0} s^\lambda +\Delta^{u}_{1} s^{\lambda+1} + \F_{{L}}^\infty(\mu_0)$ for any ${{L}}\in \big[\lambda_0+1,\min(2\lambda_0,\lambda_0+2)\big)$, where
$\Delta^{u}_1$ is a smooth function in a neighbourhood of $\mu_0$ satisfying that 
\[
\textstyle\Delta^{u}_1(\mu)=\kappa_{11}(\mu)\left(\varepsilon_1+\frac{a(b-1)}{2(a+1)b}\varepsilon_2
+\op(\|(\varepsilon_1,\varepsilon_2)\|)\right)+\kappa_{12}(\mu)\Delta^{u}_0(\mu)+\kappa_{13}(\mu)\delta_u(\mu)
\]
where $\kappa_{1i}$ are smooth functions at $\mu=\mu_0$ for $i=1,2,3$ and $\kappa_{11}(\mu_0)>0.$

\item If $a_0<-1$ then $\mathscr D_u(s;{\mu})=\delta_{u}+\Delta^{u}_{0} s^\lambda +\Delta^{u}_{2} s^{2\lambda}+\F_{{L}}^\infty(\mu_0)$ for any ${{L}}\in\big[2\lambda_0,\min(3\lambda_0,\lambda_0+1)\big)$, where
$\Delta^{u}_2$ is a smooth function in a neighbourhood of $\mu_0$ satisfying that 
\[
\textstyle \Delta^{u}_2(\mu) =\kappa_{21}(\mu)\left( \frac{2(a+2)(b-1)}{(a+1)(b-2)}\varepsilon_1+\varepsilon_2+\op(\|(\varepsilon_1,\varepsilon_2)\|)\right)+\kappa_{22}(\mu)\Delta^{u}_0(\mu)+\kappa_{23}(\mu)\delta_{u}(\mu),
\] 
where $\kappa_{2i}$ are smooth functions at $\mu=\mu_0$ for $i=1,2,3$ and $\kappa_{21}(\mu_0)>0.$

\end{enumerate}
\end{prop}

For reader's convenience the proof of \propc{prop2} is deferred to Subsection~\ref{ApB2}.

Let us fix $\mu_0=(a_0,b_0,0,0,0)$ with $a_0\in (-2,0)$ and $b_0\in (0,2)$. The differential system \refc{pert} has only two finite singularities for $\mu\approx\mu_0$, which are of focus type and close to the points $(0,\frac{1}{2})$ and $(0,\frac{b_0-2}{2b_0})$. Let us denote them by $c_u(\mu)$ and $c_\ell(\mu)$, respectively. 
We also define the parameter subset
\[
\mc Z_u\!:=\{\mu\approx\mu_0:\mathscr D_u(\,\cdot\,;\mu)\equiv 0\}.
\]
The next result shows that $\mc Z_u$ is precisely the center manifold for the focus at $(0,\frac{1}{2})$.  We remark, in connection with our discussion in \figc{FigIliev}, that the subsets $Z_0$ and $Z_1$ correspond to the components $Q_3^R$ and $Q_3^{LV}$, respectively. For completeness, we note that the combination of this result with \lemc{involucio} provides also the description of the center manifold for the focus at $(0,\frac{b_0-2}{2b_0})$.

\begin{lem}\label{centreLV}
$\mc Z_u=\textstyle\{\mu\approx\mu_0:c_u(\mu)\text{ is a center of }X_\mu\}$ and $\mc Z_u=Z_0\cup Z_1$, where
\[
Z_0\!:=\{\mu\approx\mu_0:\varepsilon_0=\varepsilon_1=\varepsilon_2=0\}
\text{ and }
Z_1\!:=\{\mu\approx\mu_0:a+b=\varepsilon_0=2\varepsilon_1+\varepsilon_2=0\}.
\]
Moreover, if $\mu\in\mc Z_u$ then the period annulus of the center at $(0,\frac{1}{2})$ is $\big\{(x,y)\in\R^2:y>0\big\}\setminus\{(0,\frac{1}{2})\}$.
\end{lem}

\begin{prova}
Let us fix~$\hat\mu\approx \mu_0$ and consider the straight line $L$ passing through the singularities $c_u(\hat\mu)$ and $c_\ell(\hat\mu)$. These two points split $L$ into three open segments where $X_{\hat\mu}$ is transverse because the vector field is quadratic. Let us denote the unbounded segment having $c_u(\hat\mu)$ as endpoint by $\Sigma_1$ and the bounded segment by $\Sigma_2.$ We parametrize them analytically by \map{\sigma_1}{(0,1)}{\Sigma_1} and \map{\sigma_2}{(0,1)}{\Sigma_2}, respectively, such that $\lim_{s\to 0}\|\sigma_1(s)\|=+\infty$, $\lim_{s\to 1}\sigma_1(s)=c_u(\hat\mu)$, $\lim_{s\to 0}\sigma_2(s)=c_\ell(\hat\mu)$ and $\lim_{s\to 1}\sigma_2(s)=c_u(\hat\mu).$ By transversality and the fact that $c_u(\hat\mu)$ and $c_\ell(\hat\mu)$ are the only finite singularities of $X_{\hat\mu}$, the application of the Poincar\'e-Bendixson Theorem shows that there is a well defined Poincar\'e map for $X_{\hat\mu}$ from $\Sigma_1$ to~$\Sigma_2.$ Taking the parametrizations previously introduced, we denote it by \map{\mathcal P_+}{(0,1)}{(0,1)}, which is an analytic function by applying the Implicit Function Theorem. 
Similarly, we denote by \map{\mathcal P_-}{(0,1)}{(0,1)} the Poincar\'e map for $-X_{\hat\mu}$ from $\Sigma_1$ to~$\Sigma_2$, which is analytic as well. Observe that, by construction, the periodic orbits surrounding $c_u({\hat\mu})$ correspond to zeros of $\mc D\!:=\mathcal P_+-\mathcal P_-.$   
Moreover $\hat\mu\in\mc Z_u$ if, and only if, $\mc D\equiv 0$ on $(0,\delta_1)$ and, on the other hand, $c_u(\hat\mu)$ is a center if, and only if, $\mc D\equiv 0$ on $(1-\delta_2,1).$ Accordingly, since $\mc D$ is analytic on $(0,1)$, this proves that 
$\hat\mu\in\mc Z_u$ if, and only if, $c_u(\hat\mu)$ is a center. So far we have proved that 
\[
\mc Z_u=\textstyle\{\mu\approx\mu_0:c_u(\mu)\text{ is a center of }X_\mu\}=:U.
\]
Our next task is to show that $U=Z_0\cup Z_1.$ To prove the inclusion $U\subset Z_0\cup Z_1$ we take any $\mu\in U$ and, due to $U=\mc Z_u$, by applying \propc{prop2} we get that $\delta_u(\mu)=0$ and $\Delta_0^u(\mu)=0$, which imply 
\[
 \varepsilon_0=0\text{ and }2\frac{\sqrt{b(a+2)}}{\sqrt{a(b-2)}}\,\varepsilon_1+\varepsilon_2=0.
\] 
Here the first equality follows by the Implicit Function Theorem using that $\delta_u|_{\varepsilon_0}\equiv 0$ and $\partial_{\varepsilon_0}\delta_u(\mu_0)\neq 0.$ 
Recall on the other hand that trace equal to zero is a necessary condition for a singular point to be a center. 
One can verify that if $\varepsilon_0=0$ then $c_u(\mu)=(0,\frac{1}{2})$ and that its trace is equal to $\varepsilon_1+\frac{1}{2}\varepsilon_2$. The vanishing of this quantity, together with the two equalities above, yields to either $\{\varepsilon_0=\varepsilon_1=\varepsilon_2=0\}$ or  
$\{a+b=\varepsilon_0=2\varepsilon_1+\varepsilon_2=0\}.$ Therefore $U\subset Z_0\cup Z_1$. To prove the reverse inclusion we note first that if $\mu\in Z_0$ then the function
\[
 H_0(x,y)=|y|^a(x^2+l y^2+m y+n),
 \]
with $ l=\frac{b}{a+2}$,  $m=-\frac{b-1}{a+1}$ and $n=\frac{b-2}{4a}$, is a global first integral of $X_\mu$. The continuity of $H_0$ at $c_u(\mu)$ implies that it must be a center, so that $\mu\in U.$ Finally, if $\mu\in Z_1$ then one can verify that 
\begin{align*}
H_1(x,y)&=|y|^a(r_1(x,y)+i\alpha_1x)^{1-i\frac{\varepsilon_2}{\alpha_1}}(r_1(x,y)-i\alpha_1x)^{1+i\frac{\varepsilon_2}{\alpha_1}}\\
&=|y|^a(r_1(x,y)^2+\alpha_1^2x^2)e^{\frac{2\varepsilon_2}{\alpha_1}\mathrm{arg}(r_1(x,y)+i\alpha_1x)},
\end{align*}
with $r_1(x,y)=2by+(2-b)+\varepsilon_2 x$ and $\alpha_1=\sqrt{4b(2-b)-\varepsilon_2^2}$, 
is a well defined first integral of $X_\mu$ outside any ray from $\{r_1(x,y)=0,x=0\}=\{c_\ell(\mu)\}$ to infinity. In particular it is continuous at $c_u(\mu)$, so that again it must be a center and $\mu\in U.$ This proves the result. 
\end{prova}

\begin{lem}\label{seusa}
Suppose that $F(u_1,u_2,v)$ is a smooth function on a neighbourhood $U$ of $(0,0,v_0)\in\R^2\times\R^n$ 
verifying $F=\mr{o}(\|(u_1,u_2)\|)$. Then there exist smooth functions $F_1(u_1,u_2,v)$ and $F_2(u_2,v)$ on $U$  such that $F(u_1,u_2,v)=u_1F_1(u_1,u_2,v)+u_2^2F_2(u_2,v)$.
\end{lem}
\begin{prova}
The hypothesis implies that $F(0,0,v)\equiv 0$ and $\partial_{u_i}F(0,0,v)\equiv 0$. Then
\begin{align*}
F(u_1,u_2,v)&=F(u_1,u_2,v)-F(0,u_2,v)+F(0,u_2,v)-F(0,0,v)\\
&=u_1\underbrace{\int_0^1\partial_{u_1}F(tu_1,u_2,v)dt}_{F_1(u_1,u_2,v)}+u_2\underbrace{\int_0^1\partial_{u_2}F(0,tu_2,v)dt}_{G(u_2,v)}
%&=u_1F_1(u_1,u_2,v)+u_2G(u_2,v)
\end{align*}
where $F_1$ and $G$ are smooth functions on $U$.
Since $G(0,v)=\partial_{u_2}F(0,0,v)=0$, we also deduce that $G(u_2,v)=u_2F_2(u_2,v)$
where
\[F_2(u_2,v)=\int_0^1\partial_{u_2}G(tu_2,v)dt\] is also smooth on $U$. Hence we can write $F=u_1F_1+u_2^2F_2$ and the result follows.
\end{prova}

In the statement of our next result, and in what follows, we denote
\begin{equation}\label{c_pm}
\varepsilon_\pm(\mu)=-\varepsilon_2\mp 2\sqrt{\frac{b(a+2)}{a(b-2)}}\varepsilon_1
\text{ and }
c_\pm(\mu)=(a+1)\pm (1-b).
\end{equation}

\begin{theo}\label{division}
Given any $\mu_0=(a_0,b_0,0,0,0)$ with $a_0\in (-2,0)\setminus\{-1\}$ and $b_0\in (0,2),$ there exist a neighbourhood $U$ of $\mu_0$ in $\R^5$ and $\delta>0$ such that $\nu=\Phi(\mu)\!:= (\varepsilon_0,\varepsilon_+,\varepsilon_-,c_+,c_-)$ is a local change of coordinates in $U$ and we can write 
\begin{equation}\label{thmdiv}
\left.\mathscr D_u(s;\mu)\right|_{\mu=\Phi^{-1}(\nu)}=\nu_1g_1(s;\nu)+\nu_2g_2(s;\nu)+\nu_3\nu_5g_3(s;\nu),
\end{equation}
where, setting $\nu_0=\Phi(\mu_0)=(0,0,0,\nu_4^0,\nu_5^0),$ 
\begin{enumerate}[$(a)$]
\item $g_1(s;\nu)=\kappa_1(\nu)+\F_\delta^\infty(\nu_0)$,

\item $g_2(s;\nu)=s^{\underline\lambda(\nu)}\big(\kappa_2(\nu)+\F_\delta^\infty(\nu_0)\big)$ where $\underline\lambda(\nu)|_{\nu=\Phi(\mu)}=-\frac{a+2}{a},$ and

\item $g_3(s;\nu)=s^{\underline{\lambda}'(\nu)}\big(\kappa_3(\nu)+\F_\delta^\infty(\nu_0)\big)$ where 
$\underline\lambda'(\nu)=\underline\lambda(\nu)+\min\big(\underline\lambda(\nu),1\big).$

\end{enumerate}
Moreover $\kappa_1,$ $\kappa_2$ and $\kappa_3$ are smooth strictly positive functions on $\Phi(U).$
\end{theo}

\begin{prova} 
The result is a consequence of \propc{prop2}. Note first that, since $\partial_{\varepsilon_0}\delta_{u}(\mu_{0})>0$ and $\delta_u|_{\varepsilon_0=0}\equiv 0,$ we can write $\delta_u=\rho_0\varepsilon_0$ with $\rho_0$ a smooth positive function. Thus, setting $\lambda'(\mu)\!:=\lambda(\mu)+\min\big(\lambda(\mu),1\big),$
\[
\alpha_1\!:=\left\{\begin{array}{ll} -\frac{2(a+2)(b-1)}{(a+1)(b-2)} & \text{if $a<-1$,}\\ -1 & \text{if $a>-1$,}\end{array}\right.\quad
\alpha_2\!:=\left\{\begin{array}{ll} -1 & \text{if $a<-1$,}\\ -\frac{a(b-1)}{2(a+1)b} & \text{if $a>-1$,}\end{array}\right.
\text{ and }\rho_1\!:=\left\{\begin{array}{ll} \kappa_{11} & \text{if $a>-1$,}\\ \kappa_{12} & \text{if $a<-1$,}\end{array}\right.
\]
we can recap the whole statement of \propc{prop2} as 
\begin{equation}\label{diveq0}
\mathscr D_u(s;\mu)=\varepsilon_0(\rho_0+\star s^\lambda+\star s^{\lambda'})+\varepsilon_+(\kappa_{01}s^\lambda+\star s^{\lambda '})+(\alpha_1\varepsilon_1+\alpha_2\varepsilon_2+\rho_2)\rho_1s^{\lambda'}+\F_{L}^\infty(\mu_0),
\end{equation}
where $\star$ are unspecified smooth functions on $\mu$, $\rho_2=\rho_2(a,b,\varepsilon_1,\varepsilon_2)=\op(\|(\varepsilon_1,\varepsilon_2)\|)$ and $L=\lambda'(\mu_0)+\delta'$ for some $\delta'>0$ small enough. We remark that $\kappa_{01},$ $\kappa_{11}$ and $\kappa_{21}$ are smooth strictly positive functions given in 
\propc{prop2}. Thus $\rho_1$ is a smooth strictly positive function as well. 

On the other hand, from \refc{c_pm} we get that $\alpha_1\varepsilon_1+\alpha_2\varepsilon_2=\alpha_+\varepsilon_++\alpha_-\varepsilon_-$ with
\begin{equation}\label{diveq3}
 \alpha_{\pm}\!:=\frac{1}{2}\left(-\alpha_2\mp\frac{\alpha_1}{2}\sqrt{\frac{a(b-2)}{b(a+2)}}\right)=
 \left\{
 \begin{array}{ll}
  \frac{1}{2}\left(
   1\pm\frac{(a+2)(b-1)}{(a+1)(b-2)}\sqrt{\frac{a(b-2)}{b(a+2)}}
  \right)
  & \text{if $a<-1$,}\\[7pt] 
  \frac{1}{4}\left(
   \frac{a(b-1)}{(a+1)b}\pm\sqrt{\frac{a(b-2)}{b(a+2)}}
  \right)
  & \text{if $a>-1$.}
 \end{array}
 \right.
\end{equation}
Hence, since $\nu=\Phi(a,b,\varepsilon_0,\varepsilon_1,\varepsilon_2)\!:= (\varepsilon_0,\varepsilon_+,\varepsilon_-,c_+,c_-)$ is a smooth change of coordinates in a neighbourhood $U$ of $\mu_0$
and $\big(\rho_2\circ\Phi^{-1}\big)(\nu)=\underline\rho_2(\varepsilon_+,\varepsilon_-,c_+,c_-)=\op(\|(\varepsilon_+,\varepsilon_-)\|)$, the application of \lemc{seusa} yields
\begin{equation}\label{diveq1}
\left.\big(\alpha_1\varepsilon_1+\alpha_2\varepsilon_2+\rho_2(\mu)\big)\right|_{\mu=\Phi^{-1}(\nu)}=
 (\underline\alpha_++\star)\varepsilon_++(\underline\alpha_-+\varepsilon_-\eta_1)\varepsilon_-=
 \star\varepsilon_++(\underline\alpha_-+\varepsilon_-\eta_1)\varepsilon_-
\end{equation}
with $\eta_1=\eta_1(\varepsilon_-,c_+,c_-)$. Here, and in what follows, for the sake of shortness, given a function $h=h(\mu)$ we denote $\underline h=\underline h(\nu)=h(\mu)|_{\mu=\Phi^{-1}(\nu)}.$ Following this convention, from 
\refc{diveq0} and \refc{diveq1} we get
\[
\mathscr D_u(s;\mu)|_{\mu=\Phi^{-1}(\nu)}=\varepsilon_0(\underline\rho_0+\star s^{\underline\lambda}+\star s^{\underline\lambda'})+\varepsilon_+(\underline\kappa_{01}s^{\underline\lambda}+\star s^{\underline\lambda '})
+\varepsilon_-(\underline\alpha_-+\varepsilon_-\eta_1)\underline\rho_1s^{\underline\lambda'}+r(s;\nu),
\]
where, setting $\nu_0\!:=\Phi(\mu_0)$ and applying assertion $(h)$ in \lemc{FLK}, $r\in\F_{L}^\infty(\nu_0)$. Note that, by \lemc{centreLV}, if $\mu\in Z_0=\{\varepsilon_0=\varepsilon_1=\varepsilon_2=0\}$ then $\mathscr D_u(s;\mu)\equiv 0.$ Thus, since $\Phi(Z_0)=\{\varepsilon_0=\varepsilon_+=\varepsilon_-=0\}$, we get that $r(s;\nu)|_{\varepsilon_0=\varepsilon_+=\varepsilon_-=0}\equiv 0$. By applying \lemc{Fdiv} this implies that the remainder can be written as $r=\varepsilon_0r_0+\varepsilon_+r_1+\varepsilon_-r_2$ with $r_i\in\F_{L}^\infty(\nu_0).$ Consequently
\[
\mathscr D_u(s;\Phi^{-1}(\nu))=\varepsilon_0\big(\underline\rho_0+\star s^{\underline\lambda}+\star s^{\underline\lambda'}+r_0(s;\nu)\big)
+\varepsilon_+\big(\underline\kappa_{01} s^{\underline\lambda}+\star s^{\underline\lambda'}+r_1(s;\nu)\big)
+\varepsilon_-\big((\underline\alpha_-+\varepsilon_-\eta_1)\underline\rho_1s^{\underline\lambda'}+r_2(s;\nu)\big).
\]
Furthermore, by \lemc{centreLV} again, if $\mu\in Z_1=\{a+b=\varepsilon_0=2\varepsilon_1+\varepsilon_2=0\}$
then $\mathscr D_u(s;\mu)\equiv 0.$ Thus, since one can easily check that $\Phi(Z_1)=\{\varepsilon_0=\varepsilon_+=c_-=0\}$, we can assert that 
\[
 \left.(\underline\alpha_-+\varepsilon_-\eta_1)\underline\rho_1s^{\underline\lambda'}+r_2(s;\nu)\right|_{\varepsilon_0=\varepsilon_+=c_-=0}\equiv 0.
 \]
Since $\rho_1(\mu_0)> 0$ and one can verify using \refc{diveq3} that $\underline\alpha_-=c_-\eta_2$ with $\eta_2(\nu_0)>0$, the above identity implies $\eta_1(\varepsilon_-,c_+,c_-)|_{c_-=0}\equiv 0$ and $r_2(s;\nu)|_{\varepsilon_0=\varepsilon_+=c_-=0}\equiv 0.$ Accordingly $\eta_1(\varepsilon_-,c_+,c_-)=c_-\eta_3(\varepsilon_-,c_+,c_-)$ and, by \lemc{Fdiv} once again, $r_2=\varepsilon_0r_3+\varepsilon_+r_4+c_-r_5$ with $r_i\in\F_{L}^\infty(\nu_0).$ Consequently
\[
\mathscr D_u(s;\Phi^{-1}(\nu))=\varepsilon_0\big(\underline\rho_0+\star s^{\underline\lambda}+\star s^{\underline\lambda'}+\bar r_0(s;\nu)\big)
+\varepsilon_+\big(\underline\kappa_{01} s^{\underline\lambda}+\star s^{\underline\lambda'}+\bar r_1(s;\nu)\big)
+c_-\varepsilon_-\big(\eta_4s^{\underline\lambda'}+r_5(s;\nu)\big).
\]
where the new remainders $\bar r_0=r_0+r_3$ and $\bar r_1=r_1+r_4$ also belong to $\F_L(\nu_0)$ and $\eta_4\!:=(\eta_2+\varepsilon_-\eta_3)\underline\rho_1$ satisfies $\eta_4(\nu_0)=(\eta_2\underline\rho_1)(\nu_0)>0.$ By applying \lemc{FLK} we can take 
$\delta>0$ small enough in order that the functions $s^{\underline\lambda}$, $s^{\underline\lambda'}$, $s^{\underline\lambda'-\underline\lambda}$, $s^{-\underline\lambda}\bar r_1$ and $s^{-\underline\lambda'}r_5$ belong to $\F_\delta(\nu_0)$. In doing so we obtain
\[
\mathscr D_u(s;\Phi^{-1}(\nu))=\varepsilon_0\big(\underline\rho_0+\F_\delta(\nu_0)\big)
+\varepsilon_+s^{\underline\lambda}\big(\underline\kappa_{01}+\F_\delta(\nu_0)\big)
+c_-\varepsilon_-s^{\underline\lambda'}\big(\eta_4+\F_\delta(\nu_0)\big).
\]
Since $\nu=(\varepsilon_0,\varepsilon_+,\varepsilon_-,c_+,c_-)$, from this expression we obtain \refc{thmdiv} by renaming the unit functions. This completes the proof.
\end{prova}

\begin{prooftext}{Proof of \teoc{thmB}.}
We prove first the assertion with regard to the hemicycle $\Gamma_u$. By \lemc{5param} it suffices to consider the quadratic 5-parameter perturbation given in \refc{pert}. We set $\mu_0=(a_0,b_0,0,0,0)$ and note that the limit cycles of $X_\mu$ that are close to 
$\Gamma_u$ in Hausdorff sense are in one to one correspondence with the isolated positive zeroes of 
\[
 \mathscr D_u(s;\mu)=D^u_+(s;\mu)-D^u_-(s;\mu),
\]  
see \figc{dib5}. That being said, by applying \teoc{division} we know that there exist a neighbourhood $U$ of~$\mu_0$ and $\delta>0$ small enough such that 
$\nu\!:=\Phi(\mu)=(\varepsilon_0,\varepsilon_+,\varepsilon_-,c_+,c_-)$ is a local change of coordinates in $U$ and
\begin{equation}\label{proBeq1}
\left.\mathscr D_u(s;\mu)\right|_{\mu=\Phi^{-1}(\nu)}=\nu_1\big(\kappa_1+\F_\delta^\infty(\nu_0)\big)+\nu_2s^{\underline\lambda(\nu)}\big(\kappa_2+\F_\delta^\infty(\nu_0)\big)+\nu_3\nu_5s^{\underline{\lambda}'(\nu)}\big(\kappa_3+\F_\delta^\infty(\nu_0)\big),
\end{equation}
where $\nu_0=\Phi(\mu_0)$, $\kappa_i(\nu_0)>0$ and $\underline\lambda'=\underline\lambda+\min(\underline\lambda,1).$  

Recall on the other hand, see \lemc{centreLV}, that $\mathscr D_u(s;\mu)\equiv 0$ if, and only if, $\mu\in Z_0\cup Z_1$ where 
\[
Z_0=\{\varepsilon_0=\varepsilon_1=\varepsilon_2=0\}
\text{ and }
Z_1=\{a+b=\varepsilon_0=2\varepsilon_1+\varepsilon_2=0\}.
\]
One can check in this respect that $\Phi(Z_0\cup Z_1)=\{\nu_1=\nu_2=\nu_3\nu_5=0\}.$ Taking this into account, and the fact that $\Phi(\mu_0)=\nu_0,$ we claim that  
there exist $s_0>0$ and an open ball $B_r(\nu_0)$ of radius $r>0$ centered $\nu_0$ such that 
\refc{proBeq1} has at most two zeros on $(0,s_0)$, counted with multiplicities, for all $\nu$ inside $V\!:=B_r(\nu_0)\cap\{\nu_1^2+\nu_2^2+(\nu_3\nu_5)^2\neq 0\}.$ This will imply, see \defic{Z_0}, that 
\[
\mathrm{Cycl}\big((\Gamma_u,X_{\mu_0}),X_\mu\big)\leqslant
\mathcal Z_0\big(\mathscr D_u(\,\cdot\,;\mu),\mu_0\big)=
\mathcal Z_0\big(\mathscr D_u(\,\cdot\,;\Phi^{-1}(\nu)),\nu_0\big)\leqslant 2.
\]  
In order to prove the claim we note first that, due to $\lim_{s\to 0}\big(\kappa_1(\nu)+\F_\delta^\infty(\nu_0)\big)=\kappa_1(\nu)\neq 0$ uniformly for $\nu\approx\nu_0$, we can take $r>0$ and $s_0>0$ small enough such that 
\[
\mathscr R_0(s;\nu)\!:=
\frac{ \left.\mathscr D_u(s;\mu)\right|_{\mu=\Phi^{-1}(\nu)}}{\kappa_1+\F_\delta^\infty(\nu_0)}=\nu_1+\nu_2s^{\underline\lambda(\nu)}\big(\kappa_4+\F_\delta^\infty(\nu_0)\big)+\nu_3\nu_5s^{\underline{\lambda}'(\nu)}\big(\kappa_5+\F_\delta^\infty(\nu_0)\big)
\]
is well defined for all $s\in (0,s_0)$ and $\nu\in B_r(\nu_0)$ and has exactly the same number of zeros, counted with multiplicities, as $\mathscr D_u(s;\Phi^{-1}(\nu))$. Accordingly $\mathcal Z_0\big(\mathscr R_0(\,\cdot\,;\nu),\nu_0\big)=\mathcal Z_0\big(\mathscr D_u(\,\cdot\,;\mu),\mu_0\big).$ We note that the second equality above follows from \refc{proBeq1} by applying \lemc{FLK} and that $\kappa_4\!:=\kappa_2/\kappa_1$ and $\kappa_5\!:\kappa_3/\kappa_1$ are strictly positive smooth functions.  
If $\nu\in V$ verifies $\nu_2=\nu_3\nu_5=0$ then $\nu_1\neq 0$ and, consequently, $\mathscr R_0(s;\nu)\neq 0$. This remark shows the validity of the claim for all $\nu\in V$ such that $\nu_2=\nu_3\nu_5=0$. To study the other cases we apply the so-called derivation-division algorithm. To this end we first observe that, by \lemc{FLK} again,
\begin{align*}
 &\partial_s\mathscr R_0(s;\nu)=\nu_2s^{\underline\lambda -1}\big(\underline\lambda \kappa_4+\F_\delta^\infty(\nu_0)\big)+\nu_3\nu_5s^{\underline{\lambda}'-1}\big(\underline\lambda'\kappa_5+\F_\delta^\infty(\nu_0)\big)\\
 \intertext{and}
 &\mathscr R_1(s;\nu)\!:=\frac{\partial_s\mathscr R_0(s;\nu)}{s^{\underline\lambda -1}\big(\underline\lambda\kappa_4+\F_\delta^\infty(\nu_0)\big)}=\nu_2+\nu_3\nu_5s^{\underline{\lambda}'-\underline{\lambda}}\big(\kappa_6+\F_\delta^\infty(\nu_0)\big),
\end{align*}
where $\kappa_6(\nu_0)>0$. Note that $\lim_{s\to 0^+}\big(\underline\lambda\kappa_4(\nu)+\F_\delta^\infty(\nu_0)\big)=\underline\lambda \kappa_4(\nu)\neq 0$ uniformly for $\nu\approx\nu_0.$
Therefore, by reducing $r>0$ and $s_0>0$ if necessary, $\mathscr R_1(s;\nu)$ is well defined  for all $s\in (0,s_0)$ and $\nu\in B_r(\nu_0)$ and has exactly the same number of zeros, counted with multiplicities, as $\partial_s\mathscr R_0(s;\nu)$. If $\nu\in V$ verifies 
$\nu_3\nu_5=0$ then we can suppose that $\nu_2\neq 0$ (otherwise we end up in the previous case) and, consequently, $\mathscr R_1(s;\nu)\neq 0$. Hence, by applying Rolle's Theorem, the claim follows in this case.
So far we have proved the validity of the claim for all $\nu\in V$ such that 
$\nu_3\nu_5= 0.$ To study the case $\nu_3\nu_5\neq 0$ we apply \lemc{FLK} once again to obtain
\[
 \mathscr R_2(s;\nu)\!:=\partial_s\mathscr R_1(s;\nu)=\nu_3\nu_5s^{\underline{\lambda}'-\underline{\lambda}-1}\big(\kappa_7+\F_\delta^\infty(\nu_0)\big)
\] 
with $\kappa_5=(\underline{\lambda}'-\underline{\lambda})\kappa_6\neq 0$ for $\nu\approx\nu_0.$ Exactly as before, by reducing $r>0$ and $s_0>0$ if necessary, we have that $\kappa_7(\nu)+\F_\delta^\infty(\nu_0)\neq 0$ for all $\nu\in B_r(\nu_0)$ and $s\in (0,s_0).$ Therefore $\mathscr R_2(s;\nu)\neq 0$ for all $\nu\in V$ with $\nu_3\nu_5\neq 0$ and $s\in (0,s_0)$ and the claim follows in this case by applying twice Rolle's Theorem. This exhausts all the possible cases for $\nu\in V$ and completes the proof of the claim. Accordingly $\mathrm{Cycl}\big((\Gamma_u,X_{\mu_0}),X_\mu\big)\leqslant 2$.

The fact that $\mathrm{Cycl}\big((\Gamma_u,X_{\mu_0}),X_\mu\big)\geqslant 2$ is also a consequence of \refc{proBeq1}. Indeed, by applying \propc{cota_inf} we can take a sequence $\lim_{n\to\infty}\hat\nu_n=\nu_0$
with $\hat\nu_n\in\Phi(U)\cap\{\nu_1=\nu_2=0\text{ and }\nu_3\nu_5\neq 0\}$ 
such that, setting $\hat\mu_n\!:=\Phi^{-1}(\hat\nu_n)$, we have $\mathrm{Cycl}\big((\Gamma_u,X_{\hat\mu_n}),X_\mu\big)\geqslant 2$ for all $n\in\N.$ Since $\lim_{n\to\infty}\hat\mu_n=\mu_0$ this clearly implies that $\mathrm{Cycl}\big((\Gamma_u,X_{\mu_0}),X_\mu\big)\geqslant 2$, as desired.

So far we have proved that the cyclicity of $\Gamma_u$ when we perturb \refc{DSq} inside the whole family of quadratic differential systems is exactly 2. In order to show this for $\Gamma_\ell$ we use an orbital symmetry that preserves the two-parameter family \refc{DSq} and interchanges $\Gamma_\ell$ with $\Gamma_u$. More concretely, we take $\phi(x,y)=(\eta x,-\eta^2y)$ with $\eta\!:=\sqrt{\frac{b_0}{2-b_0}}$. Then one can verify that the coordinate change $(\bar x,\bar y)=\phi(x,y)$, together with the time reparametrization $\bar t=\eta^{-1}t$, induce the parameter change $(\bar a_0,\bar b_0)=(a_0,2-b_0)$ in the family \refc{DSq}. Due to $\phi(\Gamma_\ell)=\Gamma_u$, the result follows because we have already proved its validity for $\Gamma_u$. This completes the proof of the result.
\end{prooftext}

\color{black}

\section{Proof of Theorem~\ref{thmC}}\label{provateoC}

\begin{lem}\label{involucio}
For each $b\in (0,2)$, define the linear map $\phi(x,y)=(\eta_bx,-\eta_b^2y)$ with $\eta_b\!:=\sqrt{\frac{b}{2-b}}$ and consider the vector field $X_\mu$ in \refc{pert}.
Then $\phi_\star(X_\mu)=\eta_b^{-1}X_{\sigma(\mu)}$ with $\sigma(a,b,\varepsilon_0,\varepsilon_1,\varepsilon_2)=
(a,2-b,-\eta_b^3\varepsilon_0,\eta_b\varepsilon_1,-\varepsilon_2/\eta_b)$.  \end{lem}

\begin{prova}
This follows by an easy computation and it is left to the reader.
\end{prova}
 \begin{figure}[t]
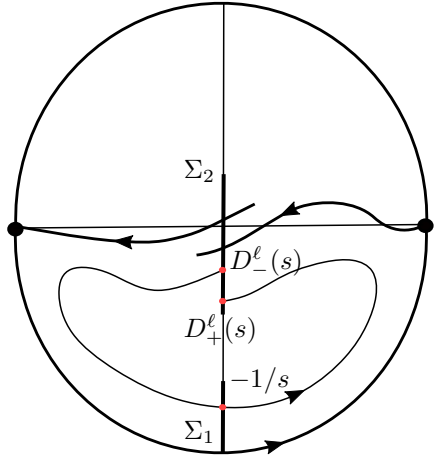

   \centering
  \begin{lpic}[l(0mm),r(0mm),t(0mm),b(5mm)]{dib5(0.75)}
    \lbl[l]{39,34;$D^\ell_-(s)$} 
    \lbl[l]{31,22;$D^\ell_+(s)$}  
    \lbl[l]{39,13;$-1/s$}      
    \lbl[l]{31,50;$\Sigma_2$}  
    \lbl[l]{31,4;$\Sigma_1$}                        
   \end{lpic}
  \caption{Dulac maps $D^\ell_\pm$ to define $\mathscr D_\ell(s;\mu)=D^\ell_+(s;\mu)-D^\ell_-(s;\mu)$. The points in red are $(0,D^\ell_\pm(s))$ and $(0,-1/s)$.}\label{dib4}
 \end{figure}
The previous result will enable us to study the limit cycles bifurcating from $\Gamma_\ell$ by taking 
advantage of \teoc{division}, which is addressed to the ones bifurcating from $\Gamma_u$. To this end 
we take two transverse sections on $x=0$, $\Sigma_1$ and $\Sigma_2$, parametrized by $s\mapsto (0,-1/s)$ with $s\in (0,\delta)$ and $s\mapsto (0,s)$ with $s\in (-\delta,\delta)$, respectively. Then, see \figc{dib4}, we consider the Dulac map $D^\ell_+(\,\cdot\,;\mu)$ for~$X_\mu$ from $\Sigma_1$ to $\Sigma_2$ and the Dulac map $D^\ell_-(\,\cdot\,;\mu)$ for~$-X_\mu$ from $\Sigma_1$ to $\Sigma_2$ and define
\[
 \mathscr D_\ell(s;\mu)\!:=D^\ell_+(s;\mu)-D^\ell_-(s;\mu).
\]
We remark that, according to the parametrization of $\Sigma_1$, the function $\mathscr D_\ell(s;\mu)$ is defined for positive $s.$ Taking these definitions into account we now prove the following result. With regard to its statement we stress that the change of parameters $\nu=\Phi(\mu)$ is the same as the one given in \teoc{division}, cf. \refc{c_pm}.

\begin{cory}\label{division2}
Given any $\mu_0=(a_0,b_0,0,0,0)$ with $a_0\in (-2,0)\setminus\{-1\}$ and $b_0\in (0,2),$ there exist a neighbourhood $U$ of $\mu_0$ in $\R^5$ and $\delta>0$ such $\nu=\Phi(\mu)\!:= (\varepsilon_0,\varepsilon_+,\varepsilon_-,c_+,c_-)$ is a local change of coordinates in $U$ and we can write 
\begin{equation*}
\left.\mathscr D_\ell(s;\mu)\right|_{\mu=\Phi^{-1}(\nu)}=\nu_1\hat g_1(s;\nu)+\nu_3\hat g_2(s;\nu)+\nu_2\nu_4\hat g_3(s;\nu),
\end{equation*}
where, setting $\nu_0=\Phi(\mu_0)=(0,0,0,\nu_4^0,\nu_5^0),$ 
\begin{enumerate}[$(a)$]
\item $\hat g_1(s;\nu)=\hat\kappa_1(\nu)+\F_\delta^\infty(\nu_0)$,

\item $\hat g_2(s;\nu)=s^{\underline\lambda(\nu)}\big(\hat\kappa_2(\nu)+\F_\delta^\infty(\nu_0)\big)$ where $\underline\lambda(\nu)|_{\nu=\Phi(\mu)}=-\frac{a+2}{a},$ and

\item $\hat g_3(s;\nu)=s^{\underline{\lambda}'(\nu)}\big(\hat\kappa_3(\nu)+\F_\delta^\infty(\nu_0)\big)$ where 
$\underline\lambda'(\nu)=\underline\lambda(\nu)+\min\big(\underline\lambda(\nu),1\big).$

\end{enumerate}
Moreover $\hat\kappa_1,$ $\hat\kappa_2$ and $\hat\kappa_3$ are smooth strictly positive functions on $\Phi(U).$

\end{cory}

\begin{prova}
By applying \lemc{involucio} (and following the notation given in its statement) one can easily show that
$D_{\pm}^\ell(s;\mu)=-\eta_b^{-2}D_{\pm}^u\big(\eta_b^{-2}s;\sigma(\mu)\big).$ Thus $\mathscr D_{\ell}(s;\mu)=-\eta_b^{-2}\mathscr D_{u}\big(\eta_b^{-2}s;\sigma(\mu)\big)$ and, consequently,
\begin{align*}
 \mathscr D_{\ell}(s;\mu)|_{\mu=\Phi^{-1}(\nu)}
 &=-\left.\eta_b^{-2}\mathscr D_{u}\big(\eta_b^{-2}s;\sigma(\mu)\big)\right|_{\mu=\Phi^{-1}(\nu)}\\[5pt]
 &=-\hat\eta_\nu^{-2}\mathscr D_{u}\big(\hat\eta_\nu^{-2}s;\sigma(\Phi^{-1}(\nu))\big)\\[5pt]
 &=-\hat\eta_\nu^{-2}\mathscr D_{u}\big(\hat\eta_\nu^{-2}s;\mu\big)|_{\mu=\Phi^{-1}(\hat\sigma(\nu))}
\end{align*}
where in the second equality we set $\hat\eta_\nu\!:=\eta_b|_{\mu=\Phi^{-1}(\nu)}=\sqrt{\frac{2-\nu_4+\nu_5}{2+\nu_4-\nu_5}}$ and in the third one $\hat\sigma\!:=\Phi\circ\sigma\circ\Phi^{-1}.$ Some computations show that
\[
 \hat\sigma(\nu)=(-\hat\eta_\nu^3\nu_1,-\nu_3/\hat\eta_\nu,-\nu_2/\hat\eta_\nu,\nu_5,\nu_4).
\]
Therefore, from the equality \refc{thmdiv} in \teoc{division}, we obtain that
\begin{align*}
 \mathscr D_{\ell}(s;\mu)|_{\mu=\Phi^{-1}(\nu)}
 &=-\hat\eta_\nu^{-2}\mathscr D_{u}\big(\hat\eta_\nu^{-2}s;\mu\big)|_{\mu=\Phi^{-1}(\hat\sigma(\nu))}\\[3pt]
 &=\hat\eta_\nu^{-2}\Big(
  \hat\eta_\nu^3\nu_1g_1\big(\hat\eta_\nu^{-2}s;\hat\sigma(\nu)\big)+\nu_3/\hat\eta_\nu g_2\big(\hat\eta_\nu^{-2}s;\hat\sigma(\nu)\big)+\nu_2\nu_4/{\hat\eta_\nu} g_3\big(\hat\eta_\nu^{-2}s;\hat\sigma(\nu)\big)
 \Big),
\end{align*}
and so the result follows setting 
\[
 \hat g_1(s;\nu)\!:=\hat\eta_\nu g_1\big(\hat\eta_\nu^{-2}s;\hat\sigma(\nu)\big),\;
 \hat g_2(s;\nu)\!:=\hat\eta_\nu^{-3} g_2\big(\hat\eta_\nu^{-2}s;\hat\sigma(\nu)\big)\text{ and } 
 \hat g_3(s;\nu)\!:=\hat\eta_\nu^{-3} g_3\big(\hat\eta_\nu^{-2}s;\hat\sigma(\nu)\big),
\]
which satisfy conditions $(a)$, $(b)$ and $(c)$ in the statement due to $\underline\lambda\circ\hat\sigma=\underline\lambda$, $\hat\eta_\nu>0$ and assertion $(h)$ of \lemc{FLK}. This concludes the proof of the result.
\end{prova}

\begin{prooftext}{Proof of \teoc{thmC}.}
By applying \lemc{5param} it suffices to consider the quadratic 5-perturbation $\{X_\mu\}$ given in \refc{pert}. To begin with let us take $\mu_0=(a_0,b_0,0,0,0)$ with $a_0\neq -1$ and note that then by \teoc{division} and \coryc{division2}, respectively, we obtain that
\begin{align}\notag
\mathscr R_u(s;\nu)\!:=\frac{\left.\mathscr D_u(s;\mu)\right|_{\mu=\Phi^{-1}(\nu)}}{g_1(s;\nu)}&=\nu_1+\nu_2\frac{g_2(s;\nu)}{g_1(s;\nu)}+\nu_3\nu_5\frac{g_3(s;\nu)}{g_1(s;\nu)}\\\label{thCeq1}
&=\nu_1+\nu_2h_2(s;\nu)+\nu_3\nu_5h_3(s;\nu)\\
\intertext{and}\notag
\mathscr R_\ell(s;\nu)\!:=\frac{\left.\mathscr D_\ell(s;\mu)\right|_{\mu=\Phi^{-1}(\nu)}}{\hat g_1(s;\nu)}&=\nu_1 
+\nu_3\frac{\hat g_2(s;\nu)}{\hat g_1(s;\nu)}+\nu_2\nu_4\frac{\hat g_3(s;\nu)}{\hat g_1(s;\nu)}\\\label{thCeq2}
&=\nu_1+\nu_3\hat h_2(s;\nu)+\nu_2\nu_4\hat h_3(s;\nu),
\end{align}
where by applying \lemc{FLK} we have that
\begin{equation}\label{thCeq3}
 h_2\!:=g_2/g_1=s^{\underline\lambda(\nu)}\big(\kappa_4(\nu)+\F_\delta^\infty(\nu_0)\big)\text{ and }
 h_3\!:=g_3/g_1=s^{\underline\lambda'(\nu)}\big(\kappa_5(\nu)+\F_\delta^\infty(\nu_0)\big)
\end{equation}
with $\kappa_i(\nu_0)>0$ and 
\begin{equation}\label{thCeq4}
 \hat h_2\!:=\hat g_2/\hat g_1=s^{\underline\lambda(\nu)}\big(\hat\kappa_4(\nu)+\F_\delta^\infty(\nu_0)\big)\text{ and }
 \hat h_3\!:=\hat g_3/\hat g_1=s^{\underline\lambda'(\nu)}\big(\hat\kappa_5(\nu)+\F_\delta^\infty(\nu_0)\big)
\end{equation}
with $\hat\kappa_i(\nu_0)>0$. Note that the limit cycles of $X_\mu$ that are close to 
$\Gamma_u$ (respectively, $\Gamma_\ell$) in Hausdorff sense are in one to one correspondence with the isolated positive zeroes of $\mathscr D_u(\,\cdot\,;\mu)$ (respectively, $\mathscr D_\ell(\,\cdot\,;\mu)$). In turn, those zeroes are in one to one correspondence with the ones of $\mathscr R_u(\,\cdot\,;\nu)$ and $\mathscr R_\ell(\,\cdot\,;\nu)$, respectively, where $\nu=\Phi(\mu).$

We claim first that $\mathrm{Cycl}\big((\{\Gamma_u,\Gamma_\ell\},X_{\mu_0}),X_\mu\big)\leqslant 3.$ We prove it by contradiction. If the claim is false then, since we know by \teoc{thmB} that  $\mr{Cycl}((\Gamma_u,X_{\mu_0}),X_\mu)=
\mr{Cycl}((\Gamma_\ell,X_{\mu_0}),X_\mu)=2$, by applying Rolle's Theorem there would exist three sequences $s_n\to 0^+$, $s_n'\to 0^+$ and $\nu_n\to\nu_0\!:=\Phi(\mu_0)$ such that $\partial_s\mathscr R_u(s_n;\nu_n)=\partial_s\mathscr R_\ell(s'_n;\nu_n)=0$ for all $n$. On the other hand, by \lemc{FLK} again, from \refc{thCeq3} we get that
\[
 \lim_{s\to 0^+}\frac{\partial_sh_3(s;\nu)}{\partial_sh_2(s;\nu)}=0\text{ uniformly on $\nu\approx\nu_0.$}
\] 
Then from  \refc{thCeq1} we obtain that $\partial_s\mathscr R_u(s_n;\nu_n)=\nu_2\partial_sh_2(s_n;\nu)+\nu_3\nu_5\partial_sh_3(s_n;\nu)\big|_{\nu=\nu_n}=0$ for all $n$ and, consequently,
\[
 \left.\frac{\nu_2}{\nu_3\nu_5}\right|_{\nu=\nu_n}=-\frac{\partial_sh_3(s_n;\nu_n)}{\partial_sh_2(s_n;\nu_n)}\to 0
 \text{ as $n\to\infty.$}
\] 
Therefore $\lim_{n\to\infty}\frac{\nu_2}{\nu_3\nu_5}\big|_{\nu=\nu_n}=0.$ Exactly the same way, but using \refc{thCeq4} and that
$\partial_s\mathscr R_\ell(s'_n;\nu_n)=0$ for all~$n,$ we get that $\lim_{n\to\infty}\frac{\nu_3}{\nu_2\nu_4}\big|_{\nu=\nu_n}=0$. The combination of both limits implies that $\frac{1}{\nu_4\nu_5}\big|_{\nu=\nu_n}$ tends to $0$ as $n\to \infty,$ which is a contradiction because $\lim_{n\to\infty}\nu_n=\nu_0\in\R^5.$ This proves the claim. 

In order to proceed we take $\varepsilon>0$ and $s_0>0$ small enough such that the functions $h_i(s;\nu)$ and $\hat h_i(s;\nu)$ for $i=1,2$ are strictly positive for all $s\in (0,s_0)$ and $\nu\in B_\varepsilon(\nu_0).$

We claim next that $\mr{Cycl}((\{\Gamma_u,\Gamma_\ell\},X_{\mu_0}),X_\mu)\geqslant 3$ for all $(a_0,b_0)\in (-2,0)\times (0,2)$ with $a_0\neq -1$ verifying that $a_0+b_0\leqslant 0$ or $a_0+2-b_0\leqslant 0.$ Let us assume for instance that $a_0+b_0\leqslant 0$ (the other case follows verbatim). To this end recall, see \refc{c_pm}, that $\nu_0=\Phi(\mu_0)=
(0,0,0,a_0+2-b_0,a_0+b_0)$ and so the fifth component of $\nu_0$ is not positive. That being said we take 
$\bar\nu\in B_\varepsilon(\nu_0)\cap\{\nu_1=\nu_2=0,\nu_3\neq 0,\nu_5<0\}$ and $s_1\in (0,s_0)$ in order that $\bar\nu_3\mathscr R_u(s_1;\bar\nu)<0$ and $\bar\nu_3\mathscr R_\ell(s_1;\bar\nu)>0$, see \refc{thCeq1} and \refc{thCeq2}, respectively. Next, by continuity, we can take $\hat\nu\in B_\varepsilon(\nu_0)\cap\{\nu_1=0,\nu_2\nu_3>0\}$ close enough to $\bar\nu$ in order to have
\[
 \hat\nu_3 \bar\nu_3>0,\quad 
 \mathscr R_u(s_1;\hat\nu) \mathscr R_u(s_1;\bar\nu)>0
 \text{ and }
 \mathscr R_\ell(s_1;\hat\nu)\mathscr R_\ell(s_1;\bar\nu)>0.
\]
We take then $s_2\in (0,s_1)$ small enough such that, on account of \refc{thCeq1} and \refc{thCeq3}, 
$\hat\nu_2\mathscr R_u(s_2;\hat\nu)>0.$ Finally, by continuity again, we choose $\nu^\star\in B_\varepsilon(\nu_0)\cap\{\nu_1\nu_2<0\}$ close enough to $\hat\nu$ such that 
\[
 \begin{array}{ll}
    \mathscr R_u(s_1;\nu^\star) \mathscr R_u(s_1;\hat\nu)>0 
      & \qquad \nu_2^\star \hat\nu_2>0\\[7pt]
    \mathscr R_\ell(s_1;\nu^\star) \mathscr R_\ell(s_1;\hat\nu)>0
     &  \qquad \nu_3^\star \hat\nu_3>0\\[7pt]
    \mathscr R_u(s_2;\nu^\star) \mathscr R_u(s_2;\hat\nu)>0 & 
 \end{array}
\]
Observe that we can also take $s_3\in (0,s_2)$ small enough such that, thanks to \refc{thCeq1} and \refc{thCeq2}, 
\[
\nu_1^\star\mathscr R_u(s_3;\nu^\star)>0\text{ and }\nu_1^\star\mathscr R_\ell(s_3;\nu^\star)>0.
\]
Then $\mathscr R_\ell(s_1;\nu^\star)\mathscr R_\ell(s_3;\nu^\star)<0$ due to $\nu_2^\star\nu_3^\star>0$ and $\nu_1^\star\nu_2^\star<0.$ Therefore, by Bolzano's Theorem, there exists $s_\ell\in (s_3,s_1)$ such that $\mathscr R_\ell(s_\ell;\nu^\star)=0.$ On the other hand,
\[
\mathscr R_u(s_3;\nu^\star)\mathscr R_u(s_2;\nu^\star)<0
\text{ and } 
\mathscr R_u(s_2;\nu^\star)\mathscr R_u(s_1;\nu^\star)<0
\]
due to $\nu_1^\star\nu_2^\star<0$ and $\nu_2^\star\nu_3^\star>0$, respectively.
Consequently, by applying Bolzano's Theorem again, there exist $s_u^1\in (s_3,s_2)$ and $s_u^2\in (s_2,s_1)$ such that 
$\mathscr R_u(s_u^1;\nu^\star)=\mathscr R_u(s_u^2;\nu^\star)=0.$ Summing-up, we have proved that there exist $\nu^\star\in B_\varepsilon(\nu_0)$ and $s_\ell,s_u^1,s_u^2\in (0,s_0)$ with $s_u^1\neq s_u^2$ such that 
\[ 
 \mathscr R_\ell(s_\ell;\nu^\star)=\mathscr R_u(s_u^1;\nu^\star)=\mathscr R_u(s_u^2;\nu^\star)=0.
\]  
Accordingly $\mr{Cycl}((\{\Gamma_u,\Gamma_\ell\},X_{\mu_0}),X_\mu)\geqslant 3$ because $\nu_0=\Phi(\mu_0)$ and we can take $\varepsilon>0$ and $s_0>0$ arbitrarily small. This proves the claim. (For completeness let us note that the case $a_0+2-b_0\leqslant 0$ leads to the simultaneous bifurcation of one limit cycle from $\Gamma_u$ and two from $\Gamma_\ell.$)

Thanks to the claim we also have that $\mr{Cycl}((\{\Gamma_u,\Gamma_\ell\},X_{\mu_0}),X_\mu)\geqslant 3$ for each $\mu_0=(a_0,b_0,0,0,0)$ with $(a_0,b_0)\in \{-1\}\times (0,2)$ because in any neighbourhood of such $\mu_0$ there exist a parameter $\mu_\star$, not in $\{a=-1\},$ verifying that $\mr{Cycl}((\{\Gamma_u,\Gamma_\ell\},X_{\mu_\star}),X_\mu)\geqslant 3$.

Our last task is to show that if $(a_0,b_0)\in\mathcal K_2,$ i.e., $a_0+b_0>0$ and $a_0+2-b_0>0$ is verified, then 
$\mr{Cycl}((\{\Gamma_u,\Gamma_\ell\},X_{\mu_0}),X_\mu)=2$. To this end, on account of
\[
 \mr{Cycl}\big((\{\Gamma_u,\Gamma_\ell\},X_{\mu_0}),X_\mu\big)\geqslant \max\left\{\mr{Cycl}\big((\Gamma_u,X_{\mu_0}),X_\mu\big),\mr{Cycl}\big((\Gamma_\ell,X_{\mu_0}),X_\mu\big)\right\}=2,
\]
it is clear that it suffices to prove that $\mr{Cycl}\big((\{\Gamma_u,\Gamma_\ell\},X_{\mu_0}),X_\mu\big)\leqslant 2.$ We shall bound this number by studying the positive zeros of $\mathscr R_u(s;\nu)$ and $\mathscr R_\ell(s;\nu)$, see \refc{thCeq1} and \refc{thCeq2}, bifurcating from $s=0$ when $\nu$ tends to $\nu_0\in\{\nu_1=\nu_2=\nu_3=0, \nu_4>0\text{ and }\nu_5>0\}.$ Recall here that $\nu_0=\Phi(\mu_0)=
(0,0,0,a_0+2-b_0,a_0+b_0)$. On account of \refc{thCeq3} and \refc{thCeq4}, respectively, the application of \lemc{FLK} yields
\[
 h_2(s;\nu)=\kappa_4(\nu)\left(s\left(1+\F_\delta^\infty(\nu_0)\right)\right)^{\underline{\lambda}(\nu)} 
\text{ and }
\hat h_2(s;\nu)=\hat\kappa_4(\nu)\left(s\left(1+\F_\delta^\infty(\nu_0)\right)\right)^{\underline{\lambda}(\nu)}.
\]
Accordingly, by applying twice \lemc{difeoF} we deduce that 
\[
(t,\nu)=\Psi(s,\nu)\!:=\big(h_2(s;\nu),\nu\big)\text{ and }(t,\nu)=\hat\Psi\big(s;\nu)\!:=\big(\hat h_2(s;\nu),\nu\big)
\]
are well defined changes of variables satisfying
\[
\Psi^{-1}(t,\nu)=\big(\sigma((t/\kappa_4(\nu))^{1/\underline{\lambda}(\nu)};\nu),\nu\big)\text{ and }
\hat\Psi^{-1}(t,\nu)=\big(\hat\sigma((t/\hat\kappa_4(\nu))^{1/\underline{\lambda}(\nu)};\nu),\nu\big),
\]
where $\sigma(u;\nu)\!:=u(1+\F_\delta^\infty(\nu_0))$ and $\hat\sigma(u;\nu)\!:=u(1+\F_\delta^\infty(\nu_0))$. 
Our aim is to apply these changes of variables in \refc{thCeq1} and~\refc{thCeq2}, respectively. To this end note that, by  \lemc{FLK} once again, from \refc{thCeq3} and \refc{thCeq4} we get 
\[
\big(h_3\circ\Psi^{-1}\big)(t;\nu)=t^{\vartheta(\nu)}(\kappa(\nu)+f(t;\nu))
\text{ and }
\big(\hat h_3\circ\hat\Psi^{-1}\big)(t;\nu)=t^{\vartheta(\nu)}(\hat\kappa(\nu)+\hat f(t;\nu))
\]
with $\vartheta(\nu)\!:=\underline{\lambda}'(\nu)/\underline{\lambda}(\nu)=1+\min(1,1/\underline{\lambda}(\nu))>1$,  $\kappa$ and $\hat\kappa$ smooth positive functions, and $f,\hat f\in\F^\infty_{\delta_1}(\nu_0)$ for some $\delta_1>0$ small enough. Accordingly, from \refc{thCeq1} and~\refc{thCeq2},
\begin{align*}
\bar{\mathscr R}_u(t;\nu)\!:=(\mathscr R_u\circ\Psi^{-1})(t,\nu)&=\nu_1+\nu_2t+\nu_3\nu_5 t^{\vartheta(\nu)}(\kappa(\nu)+f(t;\nu))
\intertext{and}
\bar{\mathscr R}_\ell(t;\nu)\!:=(\mathscr R_\ell\circ\hat \Psi^{-1})(t,\nu)&=\nu_1+\nu_3t+\nu_2\nu_4 t^{\vartheta(\nu)}(\hat\kappa(\nu)+\hat f(t;\nu)).
\end{align*}
We are now in position to prove that $\mr{Cycl}((\{\Gamma_u,\Gamma_\ell\},X_{\mu_0}),X_\mu)\leqslant 2$. By contradiction, if this number is greater than 2 then (by exchanging the subindices $u$ and $\ell$ if necessary) for all $\varepsilon>0$ there would exist 
\[
 (t_1,t_2,t_3,\nu)\in W_\varepsilon\!:=(0,\varepsilon)^3\times B_{\varepsilon}(\nu_0)\setminus\big\{t_1=t_2\text{ or }\nu_1=\nu_2=\nu_3=0\big\}
\]
verifying that
\[
\bar{\mathscr R}_u(t_1;\nu)=\bar{\mathscr R}_u(t_2;\nu)=\bar{\mathscr R}_\ell(t_3;\nu)=0.
\]
These three equalities can be written as
\[
\left(
\begin{array}{ccc}
1 & t_1 & \nu_5t_1^{\vartheta(\nu)}(\kappa(\nu)+f(t_1;\nu))\\
1 & t_2 & \nu_5t_2^{\vartheta(\nu)}(\kappa(\nu)+f(t_2;\nu))\\
1 & \nu_4t_3^{\vartheta(\nu)}(\hat\kappa(\nu)+\hat f(t_3;\nu)) & t_3
\end{array}
\right)
\left(\begin{array}{c}\nu_1\\ \nu_2\\ \nu_3\end{array}
\right)
=\left(\begin{array}{c}0\\ 0\\ 0\end{array}\right).
\]
A necessary condition for this to hold is that the determinant
\[
D(t_1,t_2,t_3;\nu)\!:=\left|\begin{array}{ccc}
1 & t_1 & \nu_5t_1^{\vartheta(\nu)}(\kappa(\nu)+f(t_1;\nu))\\
1 & t_2 & \nu_5t_2^{\vartheta(\nu)}(\kappa(\nu)+f(t_2;\nu))\\
1 & \nu_4t_3^{\vartheta(\nu)}(\hat\kappa(\nu)+\hat f(t_3;\nu)) & t_3
\end{array}
\right|
\]
is equal to zero because $(t_1,t_2,t_3,\nu)\in W_\varepsilon$. 
An easy computation shows that we can write
\begin{equation}\label{thCeq5}
\frac{D(t_1,t_2,t_3;\nu)}{t_2-t_1}=t_3\left(1-\nu_5t_3^{\vartheta(\nu)-1}\left(\hat\kappa(\nu)+\hat f(t_3;\nu)\right)A_0(t_1,t_2;\nu)\right)+\nu_5t_1t_2A_1(t_1,t_2;\nu),
\end{equation}
where, for $i=0,1$,
\begin{align*}
A_i(t_1,t_2;\nu)\!:=\,&\frac{t_2^{\vartheta(\nu)-i}(\kappa(\nu)+f(t_2;\nu))-t_1^{\vartheta(\nu)-i}(\kappa(\nu)+f(t_1;\nu))}{t_2-t_1}\\
=\,&\frac{t_2^{\vartheta(\nu)-i}-t_1^{\vartheta(\nu)-i}}{t_2-t_1}\left(\kappa(\nu)+\frac{f_i\big(t_2^{\vartheta(\nu)-i};\nu\big)-f_i\big(t_1^{\vartheta(\nu)-i};\nu\big)}{t_2^{\vartheta(\nu)-i}- t_1^{\vartheta(\nu)-i}}\right),
\end{align*}
with $f_i(r;\nu)\!:=rf\left(r^{\frac{1}{\vartheta(\nu)-i}};\nu\right)\in\F_{1+\delta_2}^\infty(\nu_0)$ for some $\delta_2>0$ small enough. By applying (twice) the Mean Value Theorem there exist $\alpha_i>0$ between $t_1$ and $t_2$, together with $\beta_i>0$ between $t_1^{\vartheta(\nu)-i}$ and $t_2^{\vartheta(\nu)-i}$, (depending both on $t_1$, $t_2$ and $\nu$) such that 
\[
 A_i(t_1,t_2;\nu)=\big(\vartheta(\nu)-i\big)\alpha_i^{\vartheta(\nu)-i-1}\big(\kappa(\nu)+\partial_rf_i(\beta_i;\nu)\big)
 \text{ for each $i=0,1.$}
\]
On account of $\vartheta(\nu_0)>1$ and $\partial_rf_i\in\F_{\delta_2}^\infty(\nu_0)$ with $\delta_2>0,$ we can assert that $A_0(t_1,t_2;\nu)$ tends to zero as $(t_1,t_2,\nu)\to (0^+,0^+,\nu_0)$  and that $A_1(t_1,t_2;\nu)>0$ on $W_\varepsilon$ for $\varepsilon>0$ small enough. Since $\nu_5>0$, from \refc{thCeq5} we conclude that $(t_2-t_1)D(t_1,t_2,t_3;\nu)>0$ for all $(t_1,t_2,t_3,\nu)\in W_\varepsilon$ with $\varepsilon>0$ small enough. This contradicts $D(t_1,t_2,t_3;\nu)=0$ and so $\mr{Cycl}\big((\{\Gamma_u,\Gamma_\ell\},X_{\mu_0}),X_\mu\big)\leqslant 2.$ This concludes the proof of the result. 
\end{prooftext}

\section{Proof of Theorem~\ref{ThmD}}\label{provateoD}

In this section we shall demonstrate \teoc{ThmD}. However, prior to that, we shall give two general results regarding  the 
different notions of cyclicity considered in this paper. Thus, otherwise explicitly stated, we consider a germ $\{X_\mu\}_{\mu\approx\mu_0}$ of an arbitrary analytic family of vector fields on $\Sc^2$. Given $U\subseteq\Sc^2,$ we denote by $\mc C(U)$ the set of compact subsets $K\subset U$ and, as usual, $N_\varepsilon(U)$ stands for the open $\varepsilon$-neighbourhood of~$U.$ We also denote the set of limit periodic sets of the germ $\{X_\mu\}_{\mu\approx\mu_0}$ by $\mc L$, so that $\mc L\subset\mc C(\mb S^2)$, see \defic{lps}.

\begin{lem}\label{lema_inside}
If $U$ is an open subset of~$\Sc^2$ then $\underline{\mr{Cycl}}^{\,U}_{\,G}\big((\partial U,X_{\mu_0}),X_\mu\big)\leqslant\mr{Cycl}_G\big((\partial U,X_{\mu_0}),X_\mu\big).$
\end{lem} 
 
 \begin{prova}
Fix a natural number $c\leqslant\underline{\mr{Cycl}}^{\,U}_{\,G}\big((\partial U,X_{\mu_0}),X_\mu\big)$. For any $\rho>0$ the set $K_\rho\!:=\overline{U}\setminus N_{\rho}(\partial U)\in\mc C(U)$ verifies $U\setminus K_\rho\subset N_{\rho}(\partial U)$ and, on the other hand (see \defic{inside}), $\mr{Cycl}_G\big((U\setminus K_\rho,X_{\mu_0}),X_\mu\big)\geqslant {c}$. This means, recall \defic{cycl-G}, that there exists $L_\rho\in\mc C(U\setminus K_\rho)$ for which $\mr{Cycl}_G\big((L_\rho,X_{\mu_0}),X_\mu\big)\geqslant {c}$, i.e., for all $\varepsilon,\delta>0$ there exists $\mu\in B_\delta(\mu_0)$ such that $X_{\mu}$ has at least ${c}$ limit cycles inside $N_{\varepsilon}(L_\rho)\subset N_{\varepsilon+\rho}(\partial U)$. According to \defic{cycl-G} again, we conclude that $\mr{Cycl}_G\big((\partial U,X_{\mu_0}),X_\mu\big)\geqslant {c}$. If $\underline{\mr{Cycl}}^{\,U}_{\,G}\big((\partial U,X_{\mu_0}),X_\mu\big)$ is finite we can take $ c=\underline{\mr{Cycl}}^{\,U}_{\,G}\big((\partial U,X_{\mu_0}),X_\mu\big)$ to obtain that $\mr{Cycl}_G\big((\partial U,X_{\mu_0}),X_\mu\big)\geqslant\underline{\mr{Cycl}}^{\,U}_{\,G}\big((\partial U,X_{\mu_0}),X_\mu\big)$, otherwise we easily deduce $\underline{\mr{Cycl}}^{\,U}_{\,G}\big((\partial U,X_{\mu_0}),X_\mu\big)=\infty=\mr{Cycl}_G\big((\partial U,X_{\mu_0}),X_\mu\big)$.
\end{prova}

\begin{lem}\label{L}
If  $K\in\mc C(\mb S^2)$ then
$
 \mr{Cycl}_G\big((K,X_{\mu_0}),X_\mu\big)=\mr{Cycl}\big((\mc L(K),X_{\mu_0}),X_\mu\big),
$ 
where $\mc L(K)=\mc L\cap \mc C(K)$.
\end{lem}

\begin{prova}
By \obsc{compact}, the set $\{\gamma \text{ limit cycle of }X_\mu \text{ contained in } N_\varepsilon(K)\}$ contains
\[
\{\gamma \text{ limit cycle of }X_\mu \text{ with }d_H(\gamma,\Gamma)<\varepsilon \text{ for some }\Gamma\in\mc L(K)\}.
\]
Accordingly, on account of Definitions~\ref{defcic} and~\ref{cycl-G},
it follows that 
\[
 \mr{Cycl}_G\big((K,X_{\mu_0}),X_\mu\big)\geqslant\mr{Cycl}\big((\mc L(K),X_{\mu_0}),X_\mu\big).
\]  
Fix a natural number $c\leqslant\mr{Cycl}_G\big((K,X_{\mu_0}),X_\mu\big)$. Then, see \defic{cycl-G} again, for any $n\in\mb N$ there exists $\mu_n\in B_{1/n}(\mu_0)$ such that $X_{\mu_n}$ has at least $c$ limit cycles $\gamma_n^1,\ldots,\gamma_n^c$ contained in $ N_{1/n}(K)$. Since $\big(\mc C(\mb S^2),d_H\big)$ is compact (see \obsc{compact} again), by taking a subsequence we can assume that $\gamma_n^j\to \Gamma^j\in\mc L(K)$ as $n\to\infty$. Consequently, for each $\varepsilon,\delta>0$, there exists $n\in\mb N$ such that $\mu_n\in B_\delta(\mu_0)$ and $d_H(\gamma_n^j,\Gamma^j)<\varepsilon$. Therefore, see \defic{defcic}, $\mr{Cycl}\big((\mc L(K),X_{\mu_0}),X_\mu\big)\geqslant c$.
If $\mr{Cycl}_G\big((K,X_{\mu_0}),X_\mu\big)$ is finite then we can take $c=\mr{Cycl}_G\big((K,X_{\mu_0}),X_\mu\big)$ to conclude that $\mr{Cycl}\big((\mc L(K),X_{\mu_0}),X_\mu\big)\geqslant \mr{Cycl}_G\big((K,X_{\mu_0}),X_\mu\big)$, and the result follows. Otherwise one can easily show that 
\[
\mr{Cycl}_G\big((K,X_{\mu_0}),X_\mu\big)=\infty=\mr{Cycl}\big((\mc L(K),X_{\mu_0}),X_\mu\big),
\]
and so the result follows as well.
\end{prova}

\begin{prooftext}{Proof of \teoc{ThmD}.} 
Let $\{X_\mu\}_{\mu\in\Lambda}$ be the whole quadratic family of vector fields and $X_{\mu_0}$ the vector field \refc{DSq} with $(a_0,b_0)\in\{(-1,1),(-\frac{1}{2},\frac{1}{2}),(-\frac{1}{2},\frac{3}{2})\}$. Setting $U=\R^2\setminus\{y=0\}$, so that $\partial U=\Gamma_u\cup\Gamma_\ell$, we get 
\begin{align*}
\underline{\mr{Cycl}}_{\,G}^{\,U}\big((\partial U,X_{\mu_0}),X_\mu\big)&\stackrel{(1)}{\leqslant}\mr{Cycl}_G\big((U,X_{\mu_0}),X_\mu\big)\stackrel{(2)}{=}2<3\stackrel{(3)}{\leqslant}\mr{Cycl}\big((\{\Gamma_u,\Gamma_\ell\},X_{\mu_0}),X_\mu\big)\\[4pt]
&\stackrel{(4)}{=}\mr{Cycl}\big((\mc L(\partial U),X_{\mu_0}),X_\mu\big)\stackrel{(5)}{=}\mr{Cycl}_G\big((\partial U,X_{\mu_0}),X_\mu\big).
\end{align*}
The inequality $(1)$ follows from \defic{inside} taking $K=\emptyset$.
The equality $(2)$ for $(a_0,b_0)=(-1,1)$ follows from \cite[Theorem 11]{FG}, for $(a_0,b_0)=(-\frac{1}{2},\frac{1}{2})$ follows from~\cite[Theorem 1.2]{PFL14} and for $(a_0,b_0)=(-\frac{1}{2},\frac{3}{2})$ is a consequence of the latter by applying Lemmas~\ref{5param} and~\ref{involucio}. The inequality $(3)$ follows from Theorem~\ref{thmC}.
The equality $(4)$ is due to the fact that the only limit periodic sets inside $\partial U=\Gamma_u\cup\Gamma_\ell$ are $\Gamma_u$ and $\Gamma_\ell.$ Finally, the equality $(5)$ follows from \lemc{L}. This proves the result. 
\end{prooftext}

We conclude the present section by resuming the remark that we made in the paragraph just after \defic{alien_defi}. The following is the intrinsic notion of alien limit cycle for an unfolding of a polycycle that we propose: 

\begin{defi}\label{robert1}
Let $\{X_\mu\}_{\mu\approx\mu_0}$ be a germ of an analytic family of vector fields on $\Sc^2$ such that $X_{\mu_0}$
has a polycycle $\Gamma$ with only a well-defined return map on one side, which is the identity. Assume moreover that~$\Gamma$ does not contain any proper subset being a limit periodic set of the unfolding. Let $U$ be the connected component of $\Sc^2\setminus\Gamma$ containing the side of $\Gamma$ where the return map is defined. Then, if
\begin{equation*}
\underline{\mr{Cycl}}^{\,U}_{\,G}\big((\partial U,X_{\mu_0}),X_\mu\big)<\mr{Cycl}_G\big((\partial U,X_{\mu_0}),X_\mu\big),
\end{equation*}
we say that an \emph{alien limit cycle bifurcation} occurs at $\partial U=\Gamma$ from inside $U$ for $\{X_\mu\}_{\mu\approx\mu_0}$.
\end{defi}

In the above definition, the hypothesis that $\Gamma$ has a well-defined return map only on one side, together with the requirement that $\Gamma$ does not contain any proper subset being a limit periodic set of the unfolding, guarantee that $\mr{Cycl}_G\big((\partial U,X_{\mu_0}),X_\mu\big)$ accounts for the limit cycles coming from $U$ only. On the other hand, if the return map is not the identity then $\underline{\mr{Cycl}}^{\,U}_{\,G}\big((\partial U,X_{\mu_0}),X_\mu\big)=0$. 

Next we particularize \defic{robert1} to the case of a 2-saddle cycle and show its relation with Melnikov functions. 
With this aim, let $\{X_\mu\}_{\mu\approx\mu_0}$ be a germ of an analytic family of vector fields on $\Sc^2$ such that~$X_{\mu_0}$
has a hyperbolic $2$-saddle cycle $\Gamma$ homeomorphic to $\Sc^1$ and with the return map being the identity. We assume moreover that
at most one saddle connection in $\Gamma$ breaks when $\mu\approx\mu_0$. Similarly as we do in \figc{dib5} we take a transversal section $\Sigma_1$ in the unbroken connection and a transversal section $\Sigma_2$ in the other one, and we consider the difference map $\mathscr D(s;\mu)$ between the corresponding Dulac maps, which is defined on $(0,s_0)$. If $f(s)$ is a smooth function on $(0,s_0)$ and $I$ is an interval inside 
 $(0,s_0)$, we denote by $Z_I(f)$ (respectively, $Z^m_I(f)$) the number of zeros of $f$ in $I$ (respectively, counted with multiplicities). Then, following this notation, we have that 
\begin{equation}\label{alien_eq2}
 \mr{Cycl}_G((\partial U,X_{\mu_0}),X_\mu)= \mr{Cycl}((\Gamma,X_{\mu_0}),X_\mu)=\inf\limits_{\varepsilon,\delta>0}\sup\limits_{\mu\in B_\delta(\mu_0)}Z_{(0,\varepsilon)}\big(\mathscr D(\cdot;\mu)\big)=:\mc Z,
 \end{equation}
where the first equality follows by using \lemc{L} and the assumption that~$\Gamma$ does not contain any proper subset being a limit periodic set of the unfolding. In the second equality we use that the limit cycles of $X_\mu$ which are Hausdorff close to $\Gamma$ correspond to small isolated zeros of the displacement function $\mathscr D(s;\mu).$ 
On the other hand, for each analytic arc $\mu=\xi(\epsilon)$ with $\xi(0)=\mu_0$ such that $\mathscr D\big(s;\xi(\epsilon)\big)\not\equiv 0$, 
we can take the Taylor's expansion at $\epsilon=0$ and write $\mathscr D\big(s;\xi(\epsilon)\big)=\epsilon^{k_\xi}(M_\xi(s)+O(\epsilon))$, where $M_\xi(s)$ is the first non-identically zero Melnikov function associated to the one-parameter unfolding $\{X_{\xi(\epsilon)}\}_{\epsilon\approx 0}$. We then define 
\[
 \mc M\!:=\inf\limits_{\varepsilon>0}\sup\limits_{\xi(0)=\mu_0} Z_{(0,\varepsilon)}^m(M_\xi),
\]
where the supremum ranges over all the analytic arcs $\mu=\xi(\epsilon)$ with $\xi(0)=\mu_0$ such that $\mathscr D\big(s;\xi(\epsilon)\big)\not\equiv 0$. (We point out that here $\varepsilon$ and $\epsilon$ play different roles.)

\begin{lem}\label{robert2}
Under the previous assumptions and notation, let $U$ be the connected component of $\Sc^2\setminus\Gamma$ where the return map of $X_{\mu_0}$ is well defined and suppose that the boundary cyclicity of $U$ from inside is finite. If $\mc M<\mc Z$ then an alien limit cycle bifurcation occurs at $\partial U=\Gamma$ from inside $U$ for $\{X_\mu\}_{\mu\approx\mu_0}$.
\end{lem}

\begin{prova}
To show the result we note that
\begin{align*}
\underline{\mr{Cycl}}_{\,G}^{\,U}((\partial U,&X_{\mu_0}),X_\mu)\stackrel{(1)}{=}\inf_{K\in \mc C(U)}\mr{Cycl}_{\,G}\big((U\setminus K,X_{\mu_0}),X_\mu\big)\stackrel{(2)}{=}\inf_{K\in\mc C(U)}\mr{Cycl}_{\,G}\big((U\setminus K,X_{\mu_0}),X_{\xi_K(\epsilon)}\big)\\[5pt]
&\stackrel{(3)}{\leqslant}\inf_{K\in\mc C(U)}Z^m_{(0,\varepsilon_K)}\big(M_{\xi_K}\big)\stackrel{(4)}{\leqslant} 
\inf\limits_{\varepsilon>0}\sup\limits_{\xi(0)=\mu_0} Z_{(0,\varepsilon)}^m(M_\xi)\stackrel{(5)}{=}\mc M<\mc Z
\stackrel{(6)}{=}\mr{Cycl}_G((\partial U,X_{\mu_0}),X_\mu).
\end{align*}
Here the equalities $(1)$ and $(5)$ follow by definition. The equality $(2)$ follows by the assumption that the boundary cyclicity of $U$ from inside is finite and applying~\cite[Theorem 1]{G08} to the period annulus $U\setminus K$ for each fixed $K\in\mc C(U)$. (It is clear that we can take $K$ to be an invariant closed disc of $X_{\mu_0}$ without loss of generality.) The inequality $(3)$ is a consequence of the Weierstrass Preparation Theorem and $(4)$ is obvious because we take the supremum over all the analytic arcs instead of the ones given by the realization theorem of Gavrilov. Finally the equality $(6)$ follows from \refc{alien_eq2}. 
\end{prova}

This lemma is related with the approach made by Dumortier, Roussarie and collaborators in \cite{alien1,alien2,alien3,alien4}. Indeed, following our notation, they say that an alien limit cycle bifurcation occurs in case that $\mc M_1<\mc Z$, where~$\mc M_1$ is defined as $\mc M$ but taking the supremum only over all the radial arcs $\xi$ for which $k_\xi=1$. For other results related with alien limit cycles the reader is referred to the contributions of Han and collaborators in~\cite{alien5,alien6,alien7} and references therein.

\appendix

\section{The asymptotic expansion of the Dulac map and related results}\label{A}

In order to prove Theorems~\ref{thmA} and~\ref{thmB} we will appeal to some previous results from \cite{MV19,MV20,MV21} about the asymptotic expansion of the Dulac map. For reader's convenience we gather these results in \propc{3punts}. To this end it is first necessary to introduce some new notation and definitions. For simplicity in the exposition, we use $\varpi\in\{\infty,\omega\}$ as a wild card in $\cc^\varpi$ for the smooth class~$\cc^\infty$ and the analytic class $\cc^\omega$. 

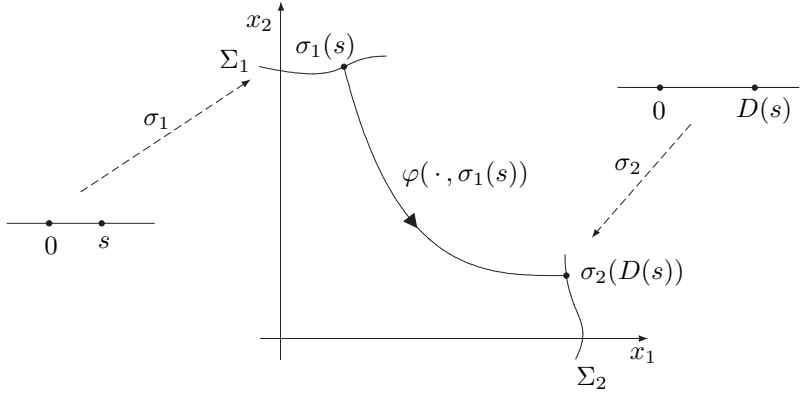
\begin{figure}[t]
   \centering
  \begin{lpic}[l(0mm),r(0mm),t(0mm),b(5mm)]{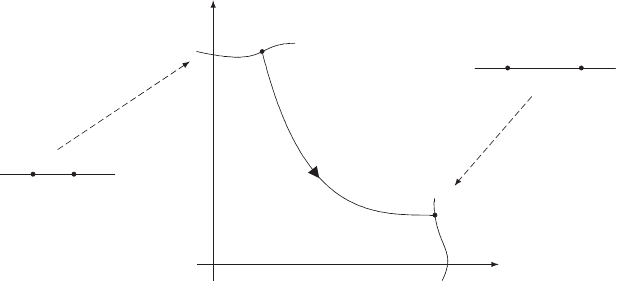}
   \lbl[l]{5,15.5;$0$}   
   \lbl[l]{12,16;$s$}   
   \lbl[l]{18,32;$\sigma_1$}
   \lbl[l]{28,40;$\Sigma_1$}      
   \lbl[l]{31.5,45;$x_2$}   
   \lbl[l]{38,42;$\sigma_1(s)$}   
   \lbl[l]{52,25;$\varphi(\,\cdot\,,\sigma_1(s))$}
   \lbl[l]{75.5,12;$\sigma_2(D(s))$}  
   \lbl[l]{75,-2;$\Sigma_2$}   
   \lbl[l]{82,1;$x_1$}
   \lbl[l]{80,26;$\sigma_2$}
   \lbl[l]{85,33.5;$0$}
   \lbl[l]{96,33.5;$D(s)$}                   
   \end{lpic}
  \caption{Definition of the Dulac map $D(\,\cdot\,;{\hat\nu})$, where $\varphi(t,p;{\hat\nu})$ is the solution of~$X_{\hat\nu}$ passing through the point $p\in\mathscr U$ at time $t=0.$}
  \label{DefTyR}
 \end{figure}

Setting $\hat\nu\!:=(\lambda,\nu)\in\hat W\!:=(0,+\infty)\times W$ with $W$ an open set of $\R^N,$ we consider the family of vector fields $\{X_{\hat\nu}\}_{{\hat\nu}\in\hat W}$ with
\begin{equation}\label{X}
 X_{\hat\nu}({x_1},{x_2})={x_1}P_1({x_1},{x_2};{\hat\nu})\partial_{x_1}+{x_2}P_2({x_1},{x_2};{\hat\nu})\partial_{x_2}
\end{equation}
where
\begin{itemize}
\item $P_1$ and $P_2$ belong to $\mathscr C^{\varpi}(\mathscr U\!\times\!\hat W)$ for some open set $\mathscr U$ of $\R^2$ containing the origin, 
\item $P_1({x_1},0;{\hat\nu})>0$ and $P_2(0,{x_2};{\hat\nu})<0$ for all $({x_1},0),(0,{x_2})\in\mathscr U$ and ${\hat\nu}\in \hat W,$
\item $\lambda=-\frac{P_2(0,0;\hat\nu)}{P_1(0,0;\hat\nu)}$.
\end{itemize}
Thus, for all $\hat\nu\in\hat W$, the origin is a hyperbolic saddle of $X_{\hat\nu}$ with the separatrices lying in the axis. We point out that here the hyperbolicity ratio of the saddle is an independent parameter, although in the applications we will have $\lambda=\lambda(\nu)$. The reason for this is that the hyperbolicity ratio turns out to be the ruling parameter in our results and, besides, having it uncoupled from the rest of parameters simplifies the notation in the statements.
Moreover, for $i=1,2,$ we consider a $\mathscr C^{\varpi}$~transverse section \map{\sigma_i}{(-\varepsilon,\varepsilon)\times \hat W}{\Sigma_i} to~$X_{{\hat\nu}}$ at $x_i=0$ defined by
 \[
  \sigma_i(s;{\hat\nu})=\bigl(\sigma_{i1}(s;{\hat\nu}),\sigma_{i2}(s;{\hat\nu})\bigr)
 \]
such that $\sigma_1(0,{\hat\nu})\in\{(0,x_2);x_2>0\}$ and $\sigma_2(0,{\hat\nu})\in\{(x_1,0);x_1>0\}$
for all ${\hat\nu}\in \hat W.$ We denote the Dulac map of~$X_{\hat\nu}$ from $\Sigma_1$ to $\Sigma_2$ by $D(\,\cdot\,;{\hat\nu})$, see \figc{DefTyR}. 
 
The asymptotic expansion of $D(s;\hat\nu)$ at $s=0$ consists of a remainder and a principal part. The principal part is given in a monomial scale that contains a deformation of the logarithm, the so-called Ecalle-Roussarie compensator, whereas the remainder has good flatness properties with respect to the parameters. We next give precise definitions of these key notions. 

\begin{defi}\label{defi_comp}
The function defined for $s>0$ and $\alpha\in\R$ by means of
 \[
  \omega(s;\alpha)\:=
  \left\{
   \begin{array}{ll}
    \frac{s^{-\alpha}-1}{\alpha} & \text{if $\alpha\neq 0,$}\\[2pt]
    -\log s & \text{if $\alpha=0,$}
   \end{array}
  \right.
 \]
is called the \emph{Ecalle-Roussarie compensator}. 
\end{defi}

\begin{defi}
Consider $K\in\Z_{\geq 0}\cup\{\infty\}$ and an open subset $U\subset\hat W\subset \R^{N+1}.$ We say that a function $\psi(s;{\hat\nu})$ belongs to the class $\mathscr C^K_{s>0}(U)$, respectively $\mathscr C^K_{s=0}(U)$, if there 
exist an open neighbourhood $\Omega$ of 
\[ 
 \{(s,{\hat\nu})\in\R^{N+2};s=0,{\hat\nu}\in U\}=\{0\}\times U
\]  
in $\R^{N+2}$ such that $(s,{\hat\nu})\mapsto \psi(s;{\hat\nu})$ is $\mathscr C^K$ on $\Omega\cap\big((0,+\infty)\times U\big)$, respectively $\Omega.$ 
\end{defi}

\begin{defi}\label{defi2} 
Consider $K\in\Z_{\geq 0}\cup\{\infty\}$ and an open subset $U\subset\hat W\subset\R^{N+1}.$ Given $L\in\R$ and ${\hat\nu}_0\in U$, we say that a function $\psi(s;{\hat\nu})\in\mathscr C^K_{s>0}(U)$ is \emph{$(L,K)$-flat with respect to $s$ at ${\hat\nu}_0$}, and we write $\psi\in\F_L^K({\hat\nu}_0)$, if for each $\ell=(\ell_0,\ldots,\ell_{N+1})\in\Z_{\geq 0}^{N+2}$ with $|\ell|=\ell_0+\ldots+\ell_{N+1}\leqslant K$  there exist a neighbourhood~$V$ of ${\hat\nu}_0$ and $C,s_0>0$ such that
\begin{equation*}
 \left|\frac{\partial^{|\ell|}\psi(s;{\hat\nu})}{\partial s^{\ell_0}
 \partial{\hat\nu}_1^{\ell_1}\cdots\partial{\hat\nu}_{N+1}^{\ell_{N+1}}}\right|\leqslant C s^{L-\ell_0}
 \text{ for all $s\in(0,s_0)$ and ${\hat\nu}\in V$.} 
\end{equation*}
If $W$ is a (not necessarily open) subset of $U$ then define $\F_L^K(W)\!:=\bigcap_{{\hat\nu}_0\in W}\F_L^K({\hat\nu}_0).$
\end{defi}

Apart from the remainder and the monomial order, the most important ingredient for our purposes is the explicit expression of the coefficients in the asymptotic expansion. In order to give them we introduce next some additional notation, where for the sake of shortness the dependence on ${\hat\nu}=(\lambda,\nu)$ is omitted. We define the functions:
\begin{equation}\label{def_fun}
\begin{array}{ll}
\dsp L_1(u)\!:=\exp\int_0^u\left(\frac{P_1(0,z)}{P_2(0,z)}+\frac{1}{\lambda}\right)\frac{dz}{z} & 
\dsp L_2(u)\!:=\exp\int_0^u\left(\frac{P_2(z,0)}{P_1(z,0)}+{\lambda}\right)\frac{dz}{z} \\[15pt]
\dsp M_1(u)\!:=L_1(u)\partial_1\!\left(\frac{P_1}{P_2}\right)(0,u)&  
\dsp M_2(u)\!:=L_2(u)\partial_2\!\left(\frac{P_2}{P_1}\right)(u,0)\\[15pt]
\end{array}
\end{equation}
On the other hand, for shortness as well, we use the compact notation $\sigma_{ijk}$ for the $k$th derivative at $s=0$ of the $j$th component of $\sigma_i(s;{\hat\nu})$, i.e., 
\[
 \sigma_{ijk}({\hat\nu})\!:=\partial^k_s\sigma_{ij}(0;{\hat\nu}).
\]
Taking this notation into account we also introduce the following real values, where once again we omit the dependence on ${\hat\nu}$:
\begin{equation}\label{def_S}
\begin{array}{l}
\dsp S_1\!:
=\frac{\sigma_{112}}{2\sigma_{111}}-\frac{\sigma_{121}}{\sigma_{120}}\left(\frac{P_1}{P_2}\right)\!(0,\sigma_{120})-\frac{\sigma_{111}}{L_1(\sigma_{120})}\gorro{M}_1(1/\lambda,\sigma_{120})\\[20pt]
\dsp S_2\!:
=\frac{\sigma_{222}}{2\sigma_{221}}-\frac{\sigma_{211}}{\sigma_{210}}\left(\frac{P_2}{P_1}\right)\!(\sigma_{210},0)-{\frac{\sigma_{221}}{L_2(\sigma_{210})}}\gorro{M}_2(\lambda,\sigma_{210})
\\[20pt]
\dsp S_3\!:
=\frac{\sigma_{221}\sigma_{210}}{L_2(\sigma_{210})}M_2'(0).
\end{array}
\end{equation}
Here $\hat M_i$ stands for a sort of incomplete Mellin transform of $M_i$ that will be defined by \propc{L8} below. 
The next proposition gathers the essential results in \cite{MV21} that we shall need to prove the first main result in the present paper.

\begin{prop}\label{3punts}
Let $D(s;{\hat\nu})$ be the Dulac map of the hyperbolic saddle \refc{X} from $\Sigma_1$ and $\Sigma_2$ and consider any $\lambda_0>0.$ Then
$D(s;{\hat\nu})=\Delta_{0}({\hat\nu})s^\lambda +\F_{\ell}^\infty(\{\lambda_0\}\times W)$ for any ${\ell}\in  \big[\lambda_0,\min(2\lambda_0,\lambda_0+1)\big)$ where $\Delta_{0}$ is a strictly positive $\cc^\varpi$ function on $\hat W$ and 
\[
\Delta_{0}(\hat\nu)=\frac{\sigma_{111}^\lambda\sigma_{120}}{L_1^\lambda(\sigma_{120})}\frac{L_2(\sigma_{210})}{\sigma_{221}\sigma_{210}^\lambda}.
\]
Moreover,
\begin{enumerate}[$(1)$]

\item If $\lambda_0>1$ then $D(s;{\hat\nu})=\Delta_{0}({\hat\nu})s^\lambda +\Delta_{1}({\hat\nu})s^{\lambda+1} +\F_{\ell}^\infty(\{\lambda_0\}\times W)$ for any ${\ell}\in \big[\lambda_0+1,\min(\lambda_0+2,2\lambda_0)\big)$ where $\Delta_1$ is a $\cc^\varpi$ function in a neighbourhood of $\{\lambda_0\}\times W$ and $\Delta_{1}(\hat\nu)=\Delta_{0}\lambda S_1.$

\item If $\lambda_0<1$ then $D(s;{\hat\nu})=\Delta_{0}({\hat\nu})s^\lambda +\Delta_{2}({\hat\nu})s^{2\lambda}+\F_{\ell}^\infty(\{\lambda_0\}\times W)$ for any ${\ell}\in  \big[2\lambda_0,\min(3\lambda_0,\lambda_0+1)\big)$ where $\Delta_2$ is a $\cc^\varpi$ function in a neighbourhood of $\{\lambda_0\}\times W$ and $\Delta_{2}(\hat\nu) =-\Delta_{0}^2S_2$.

\item If $\lambda_0=1$ then $D(s;{\hat\nu})=\Delta_{0}({\hat\nu})s^\lambda +\Delta_{3}({\hat\nu})s^{\lambda+1}\omega(s;1-\lambda)+\Delta_4(\hat\nu)s^{\lambda+1}+\F_{\ell}^\infty(\{1\}\times W)$
for any ${\ell}\in [2,3)$ where $\Delta_3$ and $\Delta_4$ are $\cc^\varpi$ functions in a neighbourhood of $\{1\}\times W$ and $\Delta_3(\hat\nu)|_{\lambda=1}=\Delta_0^2S_3|_{\lambda=1}.$
\end{enumerate}
\end{prop}

For the ease of the reader, let us explain regarding this result that the structure of the asymptotic expansion follows from  \cite[Theorem 4.1]{MV21}, whereas the properties (i.e., regularity and explicit expression) of the coefficients follow by applying Theorem A, Corollary B and Proposition 3.2 of the same paper. 
Furthermore, the flatness $\ell$ of the remainder can range in a certain interval depending on $\lambda_0.$ The left endpoint of this interval is only given for completeness to guarantee that all the monomials in the principal part are relevant (i.e., they cannot be included in the remainder). The important information about the flatness is given by the right endpoint.
A key tool in order to give a closed expression of the coefficients $\Delta_i$ is the use of a sort of incomplete Mellin transform, which is accurately defined in the next result. For a proof of this result the reader is referred to \cite[Appendix B]{MV21}.

\begin{prop}\label{L8}
Let us consider an open interval $I$ of $\R$ containing $x=0$ and an open subset $U$ of $\R^M$.
\begin{enumerate}[$(a)$]

\item Given $f(x;{\upsilon})\in\mathscr C^{\infty}(I\times U)$, there exits a unique $\hat f(\alpha,x;{\upsilon})\in\mathscr C^{\infty}((\R\setminus\Z_{\geqslant 0})\times I\times U)$ such that 
\begin{equation*}
 x\partial_x\hat f({\alpha},x;{\upsilon})-\alpha\hat f({\alpha},x;{\upsilon})=f(x;{\upsilon}).
\end{equation*}

\item If $x\in I\setminus\{0\}$ then $\partial_x(\hat f({\alpha},x;{\upsilon})|x|^{-\alpha})=f(x;{\upsilon})\frac{|x|^{-\alpha}}{x}$ and, taking any $k\in\Z_{\ge0}$ with $k>\alpha$,
\begin{equation*}
\hat f(\alpha,x;{\upsilon})=
\sum_{i=0}^{k-1}\frac{\partial_x^if(0;{\upsilon})}{i!(i-\alpha)}x^i+|x|^{\alpha}\int_0^x\!\left(f(s;{\upsilon})-T_0^{k-1}f(s;{\upsilon})\right)|s|^{-\alpha}\frac{ds}{s},
\end{equation*}
where $T_0^kf(x;{\upsilon})=\sum_{i=0}^{k}\frac{1}{i!}\partial_x^if(0;{\upsilon})x^i$ is the $k$-th degree Taylor polynomial of $f(x;{\upsilon})$ at $x=0$.

\item 
For each $(i_0,x_0,{\upsilon}_0)\in\Z_{\geqslant 0}\times I\times W$ the function $(\alpha,x,{\upsilon})\mapsto(i_0-\alpha)\hat f(\alpha,x;{\upsilon})$ extends $\mathscr C^\infty$ at $(i_0,x_0,{\upsilon}_0)$ and, moreover, it tends to $\frac{1}{i_0!}\partial_x^{i_0}f(0;{\upsilon}_0)x_0^{i_0}$ as $(\alpha,x,{\upsilon})\to (i_0,x_0,{\upsilon}_0).$

\item If $f(x;{\upsilon})$ is analytic on $I\times U$ then $\hat f(\alpha,x;{\upsilon})$ is analytic on $(\R\setminus\Z_{\geqslant 0})\times I\times U$. Finally,
for each $(\alpha_0,x_0,{\upsilon}_0)\in\Z_{\geqslant 0}\times I\times U$ the function
$(\alpha,x,{\upsilon})\mapsto(\alpha_0-\alpha)\hat f(\alpha,x;{\upsilon})$ extends analytically to $(\alpha_0,x_0,{\upsilon}_0)$.
\end{enumerate}
\end{prop}

On account of this result for each $M_i(u;\hat\nu)$ in \refc{def_fun} we have that $(\alpha,u;\hat\nu)\mapsto\hat M_i(\alpha,u;\hat\nu)$ is a well defined meromorphic function with poles only at $\alpha\in\Z_{\geq 0}$. Accordingly, see \refc{def_S}, $\hat M_1(1/\lambda,\sigma_{120})$ and $\hat M_2(\lambda,\sigma_{210})$ are the values (depending on $\hat\nu$) that we obtain by taking $\hat M_1(\alpha,u;\hat\nu)$ with $\alpha=1/\lambda$ and $u=\sigma_{120}(\hat\nu)$ and by taking $\hat M_2(\alpha,u;\hat\nu)$ with $\alpha=\lambda$ and $u=\sigma_{210}(\hat\nu),$ respectively.

The next result (see \cite[Lemma 4.3]{MV20}) is addressed to study the case in which the separatrices depicted in \figc{DefTyR} are not straight lines. 

\begin{lem}\label{rectes}
Consider a $\cc^\infty$ family $\{X_\mu\}_{\mu\in\R^N}$ of planar vector fields defined in some open set $W$ of $\R^2$. Let us fix some $\mu_0\in\R^N$ and assume that, for all $\mu\approx\mu_0,$ $X_\mu$ has a hyperbolic saddle point at $p_\mu\in W$ with $($global$)$ stable and unstable separatrices $S_\mu^+$ and $S_\mu^-$, respectively. Consider two closed connected arcs $\ell^{\pm}\subset S^\pm_{\mu_0}$, having both an endpoint at $p_{\mu_0}$. In case of a homoclinic connection $($i.e., $S^+_{\mu_0}=S^-_{\mu_0})$ we require additionally that $\ell^                                                    +\cap\ell^{-}=\{p_{\mu_0}\}$. Then there exists a neighbourhood $V$ of $(\ell^+\cup\ell^-)\times\{\mu_0\}$ in $\R^2\times\R^N$ and a $\cc^\infty$ diffeomorphism $\Phi:V\rightarrow\Phi(V)\subset\R^2\times\R^N$ with
$\Phi(x,y,\mu)=(\phi_\mu(x,y),\mu)$ such that
\[
 \Phi((S_\mu^+\times\{\mu\})\cap V)\subset\{x=0\}\times\{\mu\}\text{ and }
 \Phi((S_\mu^-\times\{\mu\})\cap V)\subset\{y=0\}\times\{\mu\}.
\]
In other words, $(\phi_\mu)_\star(X_\mu)=\hat X_\mu$ where $\hat X_\mu(x,y)=xP(x,y;\mu)\partial_x+yQ(x,y;\mu)\partial_y,$ with $P,Q\in\cc^\infty(\Phi(V)).$
\end{lem}

Next result gathers some general properties (see \cite[Lemma A.3]{MV19})  with regard to operations between functions in $\F_L^K(W)$ with  $L\in\R$.

\begin{lem}\label{FLK} 
Let $U$ and $U'$ be open sets of $\R^N$ and $\R^{N'}$ respectively and consider $W\subset U$ and $W'\subset U'.$ 
Then the following holds:
\begin{enumerate}[$(a)$]
\item $\F_L^{K}(W)\subset\F_L^{K}(\hat W)$ for any $\hat W\subset W$ and $\bigcap_n\F_L^{K}(W_n)=\F_L^{K}\left(\bigcup_n W_n\right)$.
\item $\F_L^{K}(W)\subset\F_L^{K}(W\times W')$.
\item $\mathscr C^{K}(U)\subset\mathscr C^{K}_{s=0}(U)\subset\F_0^{K}(W)$.
\item If $K\geqslant K'$ and $L\geqslant L'$ then $\F_L^{K}(W)\subset\F_{L'}^{K'}(W)$.
\item $\F_L^{K}(W)$ is closed under addition.
\item If $f\in\F_L^{K}(W)$ and $\nu\in\Z_{\ge0}^{N+1}$ with $|\nu|\leqslant K$ then 
         $\partial^\nu f\in\F_{L-\nu_0}^{K-|\nu|}(W)$.
\item $\F_L^{K}(W)\cdot\F_{L'}^{K}(W)\subset\F_{L+L'}^{K}(W)$.
\item Assume that \map{\phi}{U'}{U} is a $\mathscr C^{K}$ function with $\phi(W')\subset W$ and let us take 
        $g\in\F_{L'}^{K}(W')$ with $L'>0$ and verifying $g(s;\eta)>0$ for all $\eta\in W'$ and $s>0$ small enough. 
        Consider also any $f\in\F_L^{K}(W)$. Then $h(s;\eta)\!:=f(g(s;\eta);\phi(\eta))$ is a well-defined function 
        that belongs to $\F_{LL'}^{K}(W')$. 
\end{enumerate}
\end{lem}

\begin{obs}\label{converse}
From \defic{defi2} it follows easily that if $\partial^\nu f\in\F_{L-\nu_0}^K(W)$ for all $\nu\in\Z^{N+1}_{\geq 0}$ with $|\nu|\leqslant 1$ then $f\in \F_{L}^{K+1}(W).$ This is a sort of converse for assertion $(f)$ in \lemc{FLK}. 
\end{obs}

\begin{lem}\label{difeoF}
Let us consider $f(s;\mu)\in\F_\delta^\infty(\mu_0)$ with $\delta>0$ and define $\psi(s,\mu)=\big(s(1+f(s;\mu)),\mu\big)$ for $0<s\ll 1$ and $\mu\approx\mu_0.$ Then $\psi$ extends to a local $\mathscr C^1$ diffeomorphism  on a neighbourhood of $(0,\mu_0)$. Moreover its inverse, for $0<s\ll 1$ and $\mu\approx\mu_0,$ writes as 
$\psi^{-1}(s,\mu)=\big(s(1+g(s;\mu)),\mu\big)$ with $g\in\F_\delta^\infty(\mu_0)$.

\end{lem}
\begin{prova}
Since $f(s;\mu)\in\F_\delta^\infty(\mu_0)$ with $\delta>0$ then $sf(s;\mu)\in\F_{1+\delta}^\infty(\mu_0)$ extends to a $\mathscr C^1$ function on  some neighbourhood of $(0,\mu_0)$ by \cite[Lemma A.1]{MV19}. Thus $F(s,u,\mu)\!:=s(1+f(s;\mu))-u$ is $\mathscr C^1$ at $(0,0,\mu_0)$, $F(0,0,\mu_0)=0$ and $\partial_sF(0,0,\mu_0)=1$,
and by applying the Implicit Function Theorem there exists a unique $\mathscr C^1$ function $\sigma(u,\mu)$ on a neighbourhood $(-\varepsilon,\varepsilon)\times U$ of $(0,\mu_0)$ such that $\sigma(0,\mu_0)=0$ and $F(\sigma(u,\mu),u,\mu)\equiv 0,$ i.e., $\sigma(u,\mu)(1+f(\sigma(u,\mu);\mu))=u$. Moreover the uniqueness implies that $\sigma(0,\mu)=0$ for all $\mu\in U.$ 

We claim that $\sigma\in\bigcap_{K= 0}^\infty\F_1^K(\mu_0)=\F_1^\infty(\mu_0)$. The proof follows by induction on $K\in\Z_{\geq 0}.$ Indeed, 
due to $\sigma(0,\mu)=0$ for all $\mu\in U,$ we can write 
\[
 \sigma(u,\mu)=u\int_0^1\partial_u\sigma(tu;\mu)dt\in u\mathscr C^0_{u=0}(U)\subset\F_1^0(\mu_0),
\] 
where the inclusion follows by $(c)$ in \lemc{FLK}. Since $f\in\cc_{s>0}^\infty(U)$, by applying the Implicit Function Theorem to the equality $F(s,u,\mu)=0$ at the points $(s,u,\mu)=(\sigma(u_\star,\mu_\star),u_\star,\mu_\star)$ with $(u_\star,\mu_\star)\in (0,\varepsilon)\times U$ and taking the uniqueness of $\sigma$ into account, we deduce that $\sigma\in\cc_{s>0}^\infty(U)$. Furthermore
\begin{align*}
&\partial_u\sigma(u,\mu)=-\left(\frac{\partial_uF}{\partial_sF}\right)\!\big(\sigma(u,\mu),u,\mu\big)=\frac{1}{1+f_1(\sigma(u,\mu);\mu)}
\intertext{and}
&\partial_{\mu_i}\sigma(u,\mu)=-\left(\frac{\partial_{\mu_i}F}{\partial_s F}\right)\!\big(\sigma(u,\mu),u,\mu\big)=-\frac{\sigma(u,\mu)\partial_{\mu_i}f(\sigma(u,\mu);\mu)}{1+f_1(\sigma(u,\mu);\mu)},
\end{align*}
where $f_1\!:=f+s\partial_sf\in\F_\delta^\infty(\mu_0)$ by $(f)$ in \lemc{FLK}. On account of these two expressions, and by applying \lemc{FLK} once again, we can assert that if 
$\sigma\in\F_1^{K}(\mu_0)$ then $\partial_u\sigma\in\F_0^{K}(\mu_0)$ and $\partial_{\mu_i}\sigma\in\F_1^{K}(\mu_0),$ and consequently, see \obsc{converse}, $\sigma\in\F_1^{K+1}(\mu_0)$. Accordingly, since we already proved that $\sigma\in\F_1^0(\mu_0)$, we conclude that  
$\sigma\in\F_1^\infty(\mu_0)$ and $f(\sigma(u,\mu);\mu)\in\F_\delta^\infty(\mu_0)$ by induction. Hence 
\[
\sigma(u,\mu)=\frac{u}{1+f(\sigma(u,\mu);\mu)}=\frac{u}{1+\F_\delta^\infty(\mu_0)}=u(1+g(u,\mu))
\] 
with $g\in\F_\delta^\infty(\mu_0)$, thanks to Lemma~\ref{FLK} again. This concludes the proof of the result.
\end{prova}

The following result is a kind of division theorem among the class of flat functions and its proof can be found in \cite[Lemma 4.1]{
MV22}.

\begin{lem}\label{Fdiv}
Let us fix $L\geqslant 0$ and $n\in\N.$ If $f(s;\mu_1,\ldots,\mu_n)\in\F^\infty_L(0_n)$ verifies that 
\[
 f(s;\mu_1,\ldots,\mu_{k-1},0,\ldots,0)\equiv 0\text{, for some $k\in\{1,2,\ldots,n\}$,}
\]
then there exist $f_{k},\ldots,f_n\in\F_L^{\infty}(0_n)$ such that $f=\sum_{i=k}^n\mu_if_i$.
\end{lem}

We give at this point the precise definition of independence of functions that we use in this paper and a subsequent result addressed to obtain lower bounds for the number of bifurcating zeros. 

\begin{defi}\label{indepe}
Let $W$ be a subset of $\R^N$ (not necessarily open) and consider the functions \map{g_i}{W}{\R} for $i=1,2,\ldots,k$. The \emph{real variety} $V(g_1,g_2,\ldots,g_k)\subset W$ is defined to be the set of $\mu\in W$ such that $g_i(\mu)=0$ for $i=1,2,\ldots,k.$ We say that $g_1,g_2,\ldots,g_k$ are \emph{independent} at $\mu_{\star}\in V(g_1,g_2,\ldots,g_k)$ if the following conditions are satisfied:\begin{enumerate}[$(1)$]
\item Every neighbourhood of $\mu_{\star}$ contains two points $\mu_1,\mu_2\in V(g_1,\ldots,g_{k-1})$
         such that $g_k(\mu_1)g_k(\mu_2)<0$ (if $k=1$ then we set $V(g_1,\ldots,g_{k-1})=V(0)=W$ 
         for this to hold).
\item The varieties $V(g_1,\ldots,g_i)$, $2\leqslant i\leqslant k-1,$ are such that if $\mu_0\in V(g_1,\ldots,g_i)$ and 
        $g_{i+1}(\mu_0)\neq 0$, then every neighbourhood of $\mu_0$ contains a point $\mu\in V(g_1,\ldots,g_{i-1})$ such 
        that $g_{i}(\mu)g_{i+1}(\mu_0)<0.$
\item If $\mu_0\in V(g_1)$ and $g_2(\mu_0)\neq 0$, then every open neighbourhood of $\mu_0$ contains a point 
         $\mu\in W$ such that $g_{1}(\mu)g_{2}(\mu_0)<0.$        
\end{enumerate}
It is clear that if $W$ is an open subset of $\R^N$ and $g_i\in\mathscr C^1(W)$ for $i=1,2,\ldots,k$ then a sufficient condition for $g_1,g_2,\ldots,g_k$ to be independent at $\mu_{\star}$ is that the gradients $\nabla g_1(\mu_{\star}),\nabla g_2(\mu_{\star})\dots,\nabla g_k(\mu_{\star})$ are linearly independent vectors of $\R^{N}.$ 
\end{defi}
 
\begin{prop}\label{cota_inf}
Let $W$ be a subset of $\R^N$ $($not necessarily open$)$ and consider 
\[
 F(s;\mu)=\sum_{i=1}^n\delta_i(\mu)f_i(s;\mu)+f_{n+1}(s;\mu),
\]
where \map{f_i}{(0,\varepsilon)\!\times\! W}{\R} and \map{\delta_i}{W}{\R} are continuous functions $($with respect to the induced topology$)$. If $\mu_{\star}\in V(\delta_1,\delta_2,\ldots,\delta_{n})\subset W$ satisfies 
\begin{enumerate}[$(a)$]
\item $F(s;\mu_{\star})$ is not identically zero on $(0,\rho)$ for every $\rho\in (0,\varepsilon),$ 
\item $f_i(s;\mu)>0$, $1\leqslant i\leqslant n,$ for all $(s,\mu)$ in a neighbourhood of $(0,\mu_{\star})$ 
          in $(0,\varepsilon)\!\times\! W$,
\item $\lim_{s\to 0}\frac{f_{i+1}(s;\mu)}{f_{i}(s;\mu)}=0$, $1\leqslant i\leqslant n,$ 
         for every $\mu$ in a neighbourhood of $\mu_{\star}$ in $W$, and
\item $\delta_1,\delta_2,\ldots,\delta_n$ are independent at $\mu_{\star},$ 
\end{enumerate}
then for every neighbourhood $V$ of $\mu_\star$ in $W$ and $\rho>0$ there exists $\mu_0\in V$ such that $F(s;\mu_0)$ has at least~$n$ different zeros inside the interval $(0,\rho).$ %If, additionally, \map{F(\,\cdot\,;\mu)}{(0,\varepsilon)}{\R} is a smooth function for each $\mu\in W$ then $\mathcal Z_{0}(F(\,\cdot\,;\mu),\mu_{\star})\geqslant n$.
\end{prop}

\begin{prova}
Fix any $\rho>0$ and any neighbourhood $U$ of $\np_{\star}$ in $W$. Then, by the assumption $(a)$, there exists  $s_1\in (0,\rho)$ such that $F(s_1;\np_{\star})=f_{n+1}(s_1;\np_{\star})\neq 0.$ Suppose for instance that $F(s_1;\np_{\star})>0.$ Then, on account of~$(1)$ in \defic{indepe}, we can take $\np_1\in U\cap V(\delta_1,\delta_2,\ldots,\delta_{n-1})$ such that $\delta_{n}(\np_1)<0$ and close enough to~$\np_{\star}$ so that, by continuity, $F(s_1;\np_{1})>0.$ Observe that 
\[
 F(s;\np_1)=\delta_{n}(\np_1)f_{n}(s;\np_1)+f_{n+1}(s;\np_1).
 \]
Thus, by $(b)$ and $(c)$, $\lim_{s\to 0}\frac{F(s;\np_1)}{f_{n}(s;\np_1)}=\delta_{n}(\np_1)<0$ and we can take $s_2\in (0,s_1)$ such that $F(s_2;\np_1)<0$. 
Next, thanks to $(2)$ in \defic{indepe}, we can choose $\np_2\in U\cap V(\delta_1,\delta_2,\ldots,\delta_{n-2})$ with $\delta_{n-1}(\np_2)>0$ and close enough to $\np_1$ so that $F(s_1;\np_2)>0$ and $F(s_2;\np_2)<0$. Note that $$F(s;\np_2)=\delta_{n-1}(\np_2)f_{n-1}(s;\np_2)+\delta_{n}(\np_2)f_{n}(s;\np_2)+f_{n+1}(s;\np_2).$$ Consequently, by $(b)$ and $(c)$, $\lim_{s\to 0}\frac{F(s;\np_2)}{f_{n-1}(s;\np_2)}=\delta_{n-1}(\np_2)>0$ and we can choose $s_3\in (0,s_2)$ such that $F(s_3;\np_2)>0.$ Next we take $\np_3\in U\cap V(\delta_1,\delta_2,\ldots,\delta_{n-3})$ with $\delta_{n-2}(\np_3)<0$ and close enough to $\np_2$ so that $F(s_1;\np_3)>0$, $F(s_2;\np_3)<0$ and $F(s_3;\np_3)>0$. We repeat this process $n-2$ times after which we find a parameter $\np_{n+1}\in U$ and $0<s_{n+1}<s_{n}<\ldots<s_2<s_1<\rho$, such that $(-1)^{i+1}F(s_i;\np_{n+1})>0$ for all $i=1,2,\ldots,n+1.$ By applying Bolzano's theorem we can assert the existence of at least $n$ different zeros of $F(\,\cdot\,;\np_{n+1})$ inside the interval $(0,\rho).$ %Accordingly $\mathcal Z_0(F(\;\cdot\,;\np),\np_{\star})\geqslant n$ and this concludes the proof.
This concludes the proof of the result.
\end{prova}

\section{Deferred proofs}\label{B}

In this section, we collect the longest and most technical proofs.

\subsection{Proof of \teoc{teo1}}

\begin{prooftext}{Proof of \teoc{teo1}.}
We shall study first the Dulac map $D_+(\,\cdot\,;\mu)$ of $X_\mu$ from $\Sigma_1$ to $\Sigma_2$. For convenience 
we introduce auxiliary transverse sections $\Sigma_1^\eta$ and $\Sigma_2^\eta$ parametrized by $\sigma_1^\eta(s)=(\frac{\eta}{s},\frac{1}{s})$ and $\sigma_2^\eta(s)=(\eta,s)$ with $\eta\approx 0$, respectively. On the other hand, setting $\ell_\alpha\!:=x+\alpha( y+1)$, we perform 
the projective change of coordinates $(x_1,x_2)=\phi(x,y;\alpha)\!:=(\frac{1}{\ell_\alpha},\frac{y}{\ell_\alpha})$
to the vector field $X_\mu$, that recall is given by
\[
\sist{yf(x,y;\mu)+g(x;\mu),}{yq(x,y;\mu).}
\]
In doing so we obtain that
\[
 {\phi(\,\cdot\,;\alpha)}_\star X_\mu=\frac{1}{x_1^n}\left(
 x_1\bar P_1(x_1,x_2;\mu,\alpha)\partial_{x_1}+x_2\bar P_2(x_1,x_2;\mu,\alpha)\partial_{x_2}
 \right),
\]
where one can verify that
\begin{align}
 \textstyle P_1(x_1,x_2;\mu)&\!:=\bar P_1(x_1,x_2;\mu,0)=\textstyle -x_2x_1^n f\left(\frac{1}{x_1},\frac{x_2}{x_1}\right)-x_1^{n+1}g\left(\frac{1}{x_1}\right)\label{defP1}
 \intertext{and}
 \textstyle P_2(x_1,x_2;\mu)&\!:=\bar P_2(x_1,x_2;\mu,0)=\textstyle -x_2x_1^n f\left(\frac{1}{x_1},\frac{x_2}{x_1}\right)-x_1^{n+1}g\left(\frac{1}{x_1}\right)+x_1^n q\left(\frac{1}{x_1},\frac{x_2}{x_1}\right).
 \label{defP2}
\end{align}
Let us note at this point, see \refc{K}, that 
\begin{equation}\label{K_teo1}
 \frac{P_2(x_1,x_2;\mu)}{P_1(x_1,x_2;\mu)}=\left. 1-\frac{xq(x,y)}{yf(x,y)+g(x)}\right|_{(x,y)=\left(\frac{1}{x_1},\frac{x_2}{x_1}\right)}=K(x_1,x_2).
\end{equation}
The origin $(x_1,x_2)=(0,0)$ is a hyperbolic saddle of $x_1\phi_\star(X_\mu;\alpha)$ 
with hyperbolicity ratio equal to 
\[
 \lambda(\mu)=-K(0,0;\mu)=-1+\frac{q_n(1,0)}{g_{n+1}}.
\]
By introducing $\alpha$ and $\eta$ (that will make easier the forthcoming computations) we shall work in an extended parameter space $\bar\mu\!:=(\mu,\alpha,\eta)$ with the admissibility conditions $\Sigma_i^\eta\subset\{\ell_\alpha>0\}$ for $i=1,2$. Let $\bar D(\,\cdot\,;\mu,\alpha,\eta)$ be the Dulac map of $\bar X_{\bar\mu}\!:=x_1{\phi(\,\cdot\,;\alpha)}_\star X_\mu$ from $\Sigma_1^\eta$ to $\Sigma_2^\eta$. The key point is that, by construction, $\bar D(\,\cdot\,;\mu,\alpha,\eta)$ does not depend on~$\alpha$ and that $\bar D(\,\cdot\,;\mu,\alpha,0)=D_+(\,\cdot\,;\mu).$ 

Let us fix any admissible $\alpha_0$ and $\eta_0$. By applying \propc{3punts} to the analytic family of vector fields 
\[
 \bar X_{\bar\mu}=x_1\bar P_1(x_1,x_2;\mu,\alpha)\partial_{x_1}+x_2\bar P_2(x_1,x_2;\mu,\alpha)\partial_{x_2}
\]
at $\bar\mu_0=(\mu_0,\alpha_0,\eta_0)$ we can assert that
\[
\bar D(s;{\bar\mu})=\bar\Delta_{0}({\bar\mu})s^\lambda+
\left\{
\begin{array}{ll}
\bar\Delta_{1}({\bar\mu})s^{\lambda+1}+ \F_{\ell_1}^\infty(\bar\mu_0)&\text{ if $\lambda_0>1,$}\\[10pt]
\bar\Delta_{2}({\bar\mu})s^{2\lambda}+ \F_{\ell_2}^\infty(\bar\mu_0)&\text{ if $\lambda_0<1,$}\\[10pt]
\bar\Delta_{3}({\bar\mu})s^{\lambda+1}\omega(s;1-\lambda)+\bar\Delta_4(\bar\mu)s^{\lambda+1}+ \F_{\ell_3}^\infty(\bar\mu_0)&\text{ if $\lambda_0=1,$}
\end{array}
\right.
\]
for any ${\ell_1}\in \big[\lambda_0+1,\min(2\lambda_0,\lambda_0+2)\big)$, ${\ell_2}\in  \big[2\lambda_0,\min(3\lambda_0,\lambda_0+1)\big)$ and $\ell_3\in [2,3)$, respectively.

We remark that $\lambda_0=\lambda(\mu_0)=-K(0,0;\mu)$ because, although the new vector field $\bar X_{\bar\mu}$ depends on $\alpha$, the hyperbolicity ratio of the saddle does not. We only need to compute the coefficients of the asymptotic development for $\eta=0$ and to this aim notice that
\[
 \Delta_i^+(\mu)\!:=\bar\Delta_i(\mu,\alpha,0)=\lim_{\eta\to 0^+}\bar\Delta_i(\mu,\alpha,\eta)=\lim_{\eta\to 0^+}\bar\Delta_i(\mu,0,\eta),
\]
where in the third equality we use that the coefficients do not depend on $\alpha.$ So it suffices to perform all the computations with $\alpha=0$. The parametrisations of the auxiliary transverse sections 
$\Sigma_1^\eta$ and $\Sigma_2^\eta$ in coordinates $(x_1,x_2)$ for $\alpha=0$ are $\sigma_1(s)=(\frac{s}{\eta},\frac{1}{\eta})$ and $\sigma_2(s)=(\frac{1}{\eta},\frac{s}{\eta})$ respectively, so that $\sigma_{ijk}=\frac{1}{\eta}$ for $(i,j,k)\in\{(1,1,1),(1,2,0),(2,1,0),(2,2,1)\}$. Taking this into account, by applying \propc{3punts},
\[
 \bar\Delta_0(\mu,0,\eta)=\exp\left(\int_0^{1/\eta}\left(\frac{P_2(z,0)}{P_1(z,0)}+{\lambda}-\lambda\frac{P_1(0,z)}{P_2(0,z)}-1\right)\frac{dz}{z}\right),
\]
where
\[
\frac{P_2(z,0)}{P_1(z,0)}=1-\frac{q(1/z,0)}{zg(1/z)}\;\text{ and }\;
\frac{P_1(0,z)}{P_2(0,z)}=1+\frac{q_n(1,z)}{zf_n(1,z)+g_{n+1}-q_n(1,z)}
=1+\frac{q_n(1,z)}{\h(1,z)}.
\]
Consequently
\begin{align}\notag
 \Delta_0^+(\mu)=\bar\Delta_0(\mu,\alpha,0)=\lim_{\eta\to 0^+}\bar\Delta_0(\mu,0,\eta)
 &=\exp\left(-\int_0^{+\infty}\left(
 \frac{q(1/z,0)}{zg(1/z)}+\lambda\frac{q_n(1,z)}{\h(1,z)}
 \right)\frac{dz}{z}\right)\\[8pt]\label{teo1eq10}
 &=\exp\left(-\int_0^{+\infty}\left(
 \frac{q(w,0)}{g(w)}+\lambda\frac{q_n(w,1)}{\h(w,1)}
 \right)dw\right).
\end{align}
In the third equality we apply the Dominated Convergence Theorem \cite[Theorem 11.30]{Rudin} taking into account that the integrand does not grow faster than $z^{-2}$ at infinity, which follows by the assumptions {\bf H1} and {\bf H2}. Moreover, in the last equality, we perform the change of coordinates $w=1/z$ and take advantage of the homogeneity of the functions $q_n$ and $\h.$

Next, we compute $\bar\Delta_2(\mu,\alpha,0)$ under the assumption $\lambda_0<1$. By \propc{3punts}, $\bar\Delta_{2}=-(\bar\Delta_0)^2S_2$ with
\[
 S_2=\frac{\sigma_{222}}{2\sigma_{221}}-\frac{\sigma_{211}}{\sigma_{210}}\frac{P_2}{P_1}(\sigma_{210},0)-\frac{\sigma_{221}}{L_2(\sigma_{210})}\hat M_2(\lambda,\sigma_{210})=-\frac{1/\eta}{L_2(1/\eta)}\hat M_2(\lambda,1/\eta),
\]
where, see \refc{K_teo1} and \refc{def_fun}, $M_2(u)=L_2(u)\partial_2K(u,0)$ with
\[
 L_2(u)=\exp\int_0^u\big(K(z,0)+{\lambda}\big)\frac{dz}{z}
\]
and we take $\sigma_2(s)=(\frac{1}{\eta},\frac{s}{\eta})$ into account.
To perform the limit of $S_2$ as $\eta\to 0^+$ we need to study the growth of the functions that are involved. With this aim observe that, 
since $\lambda(\mu)<1$ for $\mu\approx\mu_0$, we can take $k=1$ in $(b)$ of \propc{L8} to get
\begin{equation}\label{teo1eq1}
 \hat M_2(\lambda,1/\eta)=\frac{M_2(0)}{-\lambda}
    + \eta^{-\lambda}\int_0^{1/\eta}(M_2(u)-M_2(0))u^{-\lambda}\frac{du}{u}.
\end{equation}
Setting $\tilde f(x_1,x_2)=x_1^n f\big(\frac{1}{x_1},\frac{x_2}{x_1}\big)$, $\tilde q(x_1,x_2)=x_1^n q\big(\frac{1}{x_1},\frac{x_2}{x_1}\big)$ and $\tilde g(x_1)=x_1^{n+1}g\big(\frac{1}{x_1}\big)$, from \refc{defP1} and \refc{defP2}, 
\[
\partial_2K(u,0)=\partial_2\left(\frac{P_2}{P_1}\right)(u,0)
%=\frac{(\tilde f(u,0)-\partial_2 \tilde q(u,0))\tilde g(u)-(\tilde g(u)-\tilde q(u,0))\tilde f(u,0)}{\tilde g(u)^2}
=\frac{\tilde q(u,0)\tilde f(u,0)-\partial_2\tilde q(u,0)\tilde g(u)}{\tilde g(u)^2}.
\]
Hence, using that $\deg(\tilde g)=n+1$ due to $g(0)\neq 0$ (see {\bf H1}) it follows that $\partial_2K(u,0)$ does not grow faster than $u^{-2}$ at $u=+\infty$. We write this assertion as $\partial_2K(u,0)\prec u^{-2}$ and in what follows we shall use this notation for shortness.
Since $(\lambda+1)\int_1^{1/\eta}\frac{dz}{z}=-\log\eta^{1+\lambda}$, an easy computation yields
\begin{align*}
\log L_2(1/\eta)&=\int_0^{1/\eta}\big(K(z,0)+\lambda\big)\frac{dz}{z}
=\int_0^{1/\eta} \left(1-\frac{q(1/z,0)}{zg(1/z)}+\lambda\right)\frac{dz}{z}\\[8pt]
&=\int_0^1 \left(1-\frac{q(1/z,0)}{zg(1/z)}+\lambda\right)\frac{dz}{z}-\int_\eta^1\frac{q(w,0)}{g(w)}dw-\log\eta^{\lambda+1}.
\end{align*}
Accordingly, due to $g(0)\neq 0$ (see {\bf H1}), setting 
\[
 G_{2}^+\!:=\int_0^1 \left(\lambda+1-\frac{q(1/z,0)}{zg(1/z)}-\frac{zq(z,0)}{g(z)}\right)\frac{dz}{z},
\]
the Dominated Convergence Theorem shows the validity of the limit
\begin{equation}\label{teo1eq2}
 \lim_{\eta\to 0^+}\eta^{\lambda+1}L_2(1/\eta)=\exp(G_{2}^+).
 \end{equation}
In particular, $L_2(u)\prec u^{\lambda+1}$. Therefore $M_2(u)=L_2(u)\partial_2K(u,0)\prec u^{\lambda-1}$. Hence, due to $\lambda<1$, we can assert that 
$(M_2(u)-M_2(0))u^{-\lambda-1}\prec u^{-\lambda-1}\prec u^{-2}.$ Accordingly, from \refc{teo1eq1},
\[
 \lim_{\eta\to 0^+}\eta^{\lambda}\hat M_2(\lambda,1/\eta)=\int_0^{+\infty}\big(M_2(u)-M_2(0)\big)\frac{du}{u^{\lambda+1}}.
\]
Finally, the combination of this with \refc{teo1eq2} yields
\begin{align}\notag
\Delta_2^+(\mu)&=\bar\Delta_2(\mu,\alpha,0)=\lim_{\eta\to 0^+}\bar\Delta_2(\mu,0,\eta)
=-\lim_{\eta\to 0^+}\big((\bar\Delta_0)^2S_2\big)(\mu,0,\eta)\\\label{teo1eq6}
&=(\Delta_0^+)^2\exp(-G_{2}^+)\int_0^{+\infty}\big(M_2(u)-M_2(0)\big)\frac{du}{u^{\lambda+1}}.
\end{align}

Our next task is to compute $\bar\Delta_1(\mu,\alpha,0)$ under the assumption $\lambda_0>1$, which is given by $\bar\Delta_1=\lambda\bar\Delta_0S_1$ thanks to the first assertion in \propc{3punts}. Taking the derivatives of
$\sigma_1(s)=(\frac{s}{\eta},\frac{1}{\eta})$ at $s=0$ into account we get that
\[
S_1(\mu,0,\eta)=\frac{\sigma_{112}}{2\sigma_{111}}-\frac{\sigma_{121}}{\sigma_{120}}\frac{P_1}{P_2}(0,\sigma_{120})-\frac{\sigma_{111}}{L_1(\sigma_{120})}\hat M_1(1/\lambda,\sigma_{120})=-\frac{1}{\eta L_1(1/\eta)}\hat M_1(1/\lambda,1/\eta),
\]
where, see \refc{K_teo1} and \refc{def_fun}, $M_1(u)=L_1(u)\partial_1\!\left(\frac{1}{K}\right)(0,u)$ with
\[
L_1(u)=\exp\int_0^u\left(\frac{1}{K(0,z)}+\frac{1}{\lambda}\right)\frac{dz}{z}.
\]
Moreover
 \begin{equation*}
\partial_1\!\left(\frac{1}{K}\right)(0,u)=\partial_1\left(1+\frac{\tilde q(x_1,x_2)}{x_2\tilde f(x_1,x_2)+\tilde g(x_1)-\tilde q(x_1,x_2)}\right)\Big|_{(x_1,x_2)=(0,u)}
\prec u^{-2}.
\end{equation*}
Here the assertion with regard to the growth at infinity is a consequence of $f_n(0,1)\neq 0$ (see {\bf H2}), which implies that $\tilde f(0,u)$ has degree exactly $n$. On the other hand, by applying $(b)$ in \propc{L8} and taking $1/\lambda<1$ into account, we get
\begin{equation}\label{teo1eq3}
\hat M_1(1/\lambda,1/\eta)=-\lambda M_1(0)+\eta^{-1/\lambda}\int_0^{1/\eta}(M_1(u)-M_1(0))u^{-1/\lambda}\frac{du}{u}.
\end{equation}
Moreover
\begin{align*}
\log L_1(1/\eta)=&\int_0^{1/\eta}\left(\frac{1}{K(0,z)}+\frac{1}{\lambda}\right)\frac{dz}{z}
=\int_0^{1/\eta}\left(\frac{q_n(1,z)}{\h(1,z)}+1+\frac{1}{\lambda}\right)\frac{dz}{z}\\[10pt]
=&\int_0^1\left(\frac{q_n(1,z)}{\h(1,z)}+1+\frac{1}{\lambda}\right)\frac{dz}{z}
+
\int_\eta^1\frac{q_n(w,1)}{\h(w,1)}dw-\left(1+\frac{1}{\lambda}\right)\log\eta,
\end{align*}
where in the last equality we use the coordinate change $z=1/w.$ Consequently, by applying the Dominated Convergence Theorem using that $f_n(0,1)\neq 0,$
\begin{equation}\label{teo1eq4}
\lim_{\eta\to 0^+}\eta^{1+1/\lambda}L_1(1/\eta)=\exp(G_{1}^+),
\end{equation}
where
\[
G_{1}^+\!:=
\int_0^1\left(\frac{q_n(1,z)}{\h(1,z)}+1+\frac{1}{\lambda}+\frac{zq_n(z,1)}{\h(z,1)}\right)\frac{dz}{z}.
\]
This implies in particular that $L_1(u)\prec u^{1+1/\lambda}$ and, accordingly, $M_1(u)\prec u^{-1+1/\lambda}$. The combination of this, together with \refc{teo1eq3} and \refc{teo1eq4}, yields
\[
\lim_{\eta\to 0^+}S_1(\mu,0,\eta)
%=-\lim_{\eta\to 0^+}\frac{\hat M_1(1/\lambda,1/\eta)}{\eta L_1(1/\eta)}
=-\lim_{\eta\to 0^+}\frac{\eta^{1/\lambda}\hat M_1(1/\lambda,1/\eta)}{\eta^{1+1/\lambda}L_1(1/\eta)}
=-\exp(-G_{1}^+)
\int_0^{+\infty}(M_1(u)-M_1(0))\frac{du}{u^{1+1/\lambda}}.
\]
Therefore
\begin{align}\notag
\Delta_1^+(\mu)=&\,\bar\Delta_1(\mu,\alpha,0)=\big(\lambda\bar\Delta_0S_1\big)(\mu,\alpha,0)\\\label{teo1eq8}
=&-\lambda\Delta_0^+\exp(-G_{1}^+)\int_0^{+\infty}\big(M_1(u)-M_1(0)\big)\frac{du}{u^{1+1/\lambda}}.
\end{align}

Now we turn to the computation of the coefficient $\bar\Delta_3(\mu,\alpha,0)$ in case that $\lambda(\mu)=1.$ By the third assertion in \propc{3punts} we have that $\left.\bar\Delta_3\right|_{\lambda=1}=\left.(\bar\Delta_0)^2S_3\right|_{\lambda=1}$ with 
\[
 S_3=\frac{\sigma_{221}\sigma_{210}}{L_2(\sigma_{210})}M_2'(0).
\] 
Note that if $\lambda=1$ then the quotient $\frac{\sigma_{221}\sigma_{210}}{L_2(\sigma_{210})}=\frac{1}{\eta^2L_2(1/\eta)}$ tends to $\exp(-G_2^+)$ as $\eta\to 0^+$ thanks to \refc{teo1eq2}, which is true for any $\lambda>0.$ Consequently, if $\lambda=1$ then
\begin{equation}\label{teo1eq4.5}
 \Delta_3^+(\mu)=\bar\Delta_3(\mu,\alpha,0)=\lim_{\eta\to 0^+}\bar\Delta_3(\mu,0,\eta)=
 (\Delta_0^+)^2\exp(-G_2^+)M'_2(0).
\end{equation}

So far we have proved that
\[
 D_+(s;{\mu})=\Delta_{0}^+({\mu})s^\lambda+
\left\{
\begin{array}{ll}
\Delta_{1}^+({\mu})s^{\lambda+1}+ \F_{\ell_1}^\infty(\mu_0)&\text{ if $\lambda_0>1,$}\\[10pt]
\Delta_{2}^+({\mu})s^{2\lambda}+ \F_{\ell_2}^\infty(\mu_0)&\text{ if $\lambda_0<1$}\\[10pt]
\Delta_{3}^+({\mu})s^{\lambda+1}\omega(s;1-\lambda)+\Delta_4^+(\mu)s^{\lambda+1}+ \F_{\ell_3}^\infty(\mu_0)&\text{ if $\lambda_0=1.$}
\end{array}
\right.
\]
We turn next to the study of the Dulac map $D_-(\,\cdot\,;\mu)$ of $-X_\mu$ from $\Sigma_1$ to $\Sigma_2$. To this aim the idea is to take advantage of the previous results for $D_+(\,\cdot\,;\mu)$ using the fact that 
$(x,y)\mapsto (-x,y)$ sends $-X_\mu$ to
\[
 \tilde X_{\mu}\!:=\big(y\bar f(x,y;\mu)+\bar g(x;\mu)\big)\partial_x+y\bar q(x,y;\mu)\partial_y
\]
with $\bar f(x,y)=f(-x,y)$, $\bar g(x)=g(-x)$ and $\bar q(x,y)=-q(-x,y).$ In particular, following the obvious notation one can check that $\bar\ell_{n+1}(x,y)=\h(-x,y),$ together with
\begin{equation}\label{teo1eq5}
 \bar L_i(u)=L_i(-u)\text{ and }\bar M_i(u)=-M_i(-u)\text{ for $i=1,2$}
\end{equation}
is verified. By applying the above assertions to the Dulac map of $\tilde X_{\mu}$ from $\Sigma_1$ to $\Sigma_2$ we get that
\[
 D_-(s;{\mu})=\Delta_{0}^-({\mu})s^\lambda+
\left\{
\begin{array}{ll}
\Delta_{1}^-({\mu})s^{\lambda+1}+  \F_{\ell_1}^\infty(\mu_0)&\text{ if $\lambda_0>1,$}\\[10pt]
\Delta_{2}^-({\mu})s^{2\lambda}+ \F_{\ell_2}^\infty(\mu_0)&\text{ if $\lambda_0<1,$}\\[10pt]
\Delta_{3}^-({\mu})s^{\lambda+1}\omega(s;1-\lambda)+\Delta_4^-(\mu)s^{\lambda+1}+ \F_{\ell_3}^\infty(\mu_0)&\text{ if $\lambda_0=1,$}
\end{array}
\right.
\]
where each coefficient $\Delta_{i}^-$ is the counterpart for $\tilde X_\mu$ of the coefficient $\Delta_i^+$ that we have obtained previously for $X_\mu.$
We can thus assert that
\begin{align*}
 \mathscr D(s;\mu)&=D_+(s;{\mu})-D_-(s;\mu)\\[10pt]
 &=\Delta_{0}({\mu})s^\lambda+
\left\{
\begin{array}{ll}
\Delta_{1}({\mu})s^{\lambda+1}+ \F_{\ell_1}^\infty(\mu_0)&\text{ if $\lambda_0>1,$}\\[10pt]
\Delta_{2}({\mu})s^{2\lambda}+ \F_{\ell_2}^\infty(\mu_0)&\text{ if $\lambda_0<1,$}\\[10pt]
\Delta_{3}({\mu})s^{\lambda+1}\omega(s;1-\lambda)+\Delta_4(\mu)s^{\lambda+1}+ \F_{\ell_3}^\infty(\mu_0)&\text{ if $\lambda_0=1,$}
\end{array}
\right.
\end{align*}
where $\Delta_i\!:=\Delta_i^+-\Delta_i ^-$ for $i=0,1,2,3,4.$ Our next task is to compute each coefficient. Note that, from \refc{teo1eq10},
\begin{align*}
 \Delta_0^-(\mu)&
= \exp\left(\int_0^{+\infty}\left(
 \frac{q(-w,0)}{g(-w)}+\lambda \frac{q_n(-w,1)}{\h(-w,1)}
 \right)dw\right)\\[5pt]
&= \exp\left(\int^0_{-\infty}\left(
 \frac{q(z,0)}{g(z)}+\lambda \frac{q_n(z,1)}{\h(z,1)}
 \right)dz\right).
\end{align*}
It is clear now that
\[
\log(\Delta_0^+)-\log(\Delta_0^-)
=-\int^{+\infty}_{-\infty}\left(
 \frac{q(z,0)}{g(z)}+\lambda \frac{q_n(z,1)}{\h(z,1)}
 \right)dz
=:\!d_{0}
\]
On account of this, and the fact that $x\mapsto\log x$ is strictly increasing, the application of the mean value theorem shows that $\Delta_0=\Delta_0^+-\Delta_0^-=\kappa_0d_{0}$ for some analytic function $\kappa_0$ with $\kappa_0(\mu_0)>0.$

We turn next to the computation of $\Delta_2^-$. To this end we again take advantage of the expression of $\Delta_2^+$ thanks to the fact that $(x,y)\mapsto (-x,y)$ sends $-X_\mu$ to $\tilde X_\mu$. In doing so, recall \refc{teo1eq5}, from \refc{teo1eq6} we get
\begin{equation}\label{teo1eq7}
 \Delta_2^-=-(\Delta_0^-)^2\exp(-G_{2}^-)\int_0^{+\infty}\big(M_2(-u)-M_2(0)\big)\frac{du}{u^{\lambda+1}}.
\end{equation}
where
\[
 G_{2}^-\!:=\int_0^1 \left(\lambda+1+\frac{q(-1/z,0)}{zg(-1/z)}+\frac{zq(-z,0)}{g(-z)}\right)\frac{dz}{z}.
\]
In order to study $\Delta_2=\Delta_2^+-\Delta_2^-$ we first observe that
\begin{align*}
 G_{2}^--G_{2}^+&=
 \int_0^1\left(\frac{q(1/z,0)}{zg(1/z)}+\frac{q(-1/z,0)}{zg(-1/z)}\right)\frac{dz}{z}
 +\int_0^1\left(\frac{q(z,0)}{g(z)}+\frac{q(-z,0)}{g(-z)}\right)dz\\
 &=-\int_{+\infty}^1\left(\frac{q(u,0)}{g(u)}+\frac{q(-u,0)}{g(-u)}\right)du
 +\int_0^1\left(\frac{q(z,0)}{g(z)}+\frac{q(-z,0)}{g(-z)}\right)dz\\
 &=\int_0^{+\infty}\left(\frac{q(z,0)}{g(z)}+\frac{q(-z,0)}{g(-z)}\right)dz=:\!G_{2},
\end{align*}
where in the second equality we perform the change of coordinates $u=1/z$. Then, from \refc{teo1eq6} and \refc{teo1eq7},
\begin{align*}
\Delta_2=\Delta_2^+-\Delta_2^-&=\exp(-G_{2}^-)
\left(
(\Delta_0^-)^2\int_0^{+\infty}\big(M_2(-u)-M_2(0)\big)\frac{du}{u^{\lambda+1}}\right.\\
&\qquad\left.+(\Delta_0^+)^2\exp(G_{2})\int_0^{+\infty}\big(M_2(u)-M_2(0)\big)\frac{du}{u^{\lambda+1}}
\right)\\
&=\bar\kappa_2\Delta_0+\exp(-G_{2}^-)(\Delta_0^-)^2
\left(
\int_0^{+\infty}\big(M_2(-u)-M_2(0)\big)\frac{du}{u^{\lambda+1}}\right.\\
&\qquad\left.+\exp(G_{2})\int_0^{+\infty}\big(M_2(u)-M_2(0)\big)\frac{du}{u^{\lambda+1}}
\right)\\
&=\bar\kappa_2\Delta_0+\kappa_2F_{2},
\end{align*}
where in the second equality we use that $G_{2}^--G_{2}^+=G_{2}$, in the third one we plug $\Delta_0^+=\Delta_0+\Delta_0^-$ to get an analytic function $\bar\kappa_2=\bar\kappa_2(\mu)$ multiplying $\Delta_0$ and in the last one we set $\kappa_2=\exp(-G_{2}^-)(\Delta_0^-)^2$. Accordingly $\Delta_2=\kappa_2F_{2}+\bar\kappa_2\Delta_0$ with $\kappa_2(\mu_0)>0,$ so the assertion $(2)$ in the statement is true.

In order to obtain the expression for $\Delta_1=\Delta_1^+-\Delta_1^-$ we follow the same strategy as before. First we take advantage of the expression of $\Delta_1^+$ in \refc{teo1eq8} and the equalities in \refc{teo1eq5} to get that
\begin{equation}\label{teo1eq9}
\Delta_1^-=\lambda\Delta_0^-\exp(-G_{1}^-)\int_0^{+\infty}\big(M_1(-u)-M_1(0)\big)\frac{du}{u^{1+1/\lambda}},
\end{equation}
where
\begin{align*}
 G_{1}^-\!:=&\int_0^1\left(-\frac{q_n(-1,z)}{\h(-1,z)}+1+\frac{1}{\lambda}-\frac{zq_n(-z,1)}{\h(-z,1)}\right)\frac{dz}{z}\\[4pt]
   =&-\int_{-1}^0\left(\frac{q_n(1,u)}{\h(1,u)}+1+\frac{1}{\lambda}+\frac{uq_n(u,1)}{\h(u,1)}\right)\frac{du}{u}.
\end{align*}
Here we use first the homogeneity of $q_n$ and $\h$ and then we perform the change of coordinates $u=-z.$
Consequently
\[
 G_{1}^+-G_{1}^-=\int_{-1}^1\left(\frac{q_n(1,z)}{\h(1,z)}+1+\frac{1}{\lambda}+\frac{zq_n(z,1)}{\h(z,1)}\right)\frac{dz}{z}=:\!G_{1}
\]
On account of this, the combination of \refc{teo1eq8} and \refc{teo1eq9} yields
\begin{align*}
\Delta_1=\Delta_1^+-\Delta_1^-&=-\lambda\exp(-G_{1}^+)
\left(
\Delta_0^+\int_0^{+\infty}\big(M_1(u)-M_1(0)\big)\frac{du}{u^{1+1/\lambda}}\right.\\
&\qquad\left.+\Delta_0^-\exp(G_{1})\int_0^{+\infty}\big(M_1(-u)-M_1(0)\big)\frac{du}{u^{1+1/\lambda}}
\right)\\
&=\bar\kappa_1\Delta_0-\lambda\Delta_0^+\exp(-G_{2}^+)
\left(
\int_0^{+\infty}\big(M_1(u)-M_1(0)\big)\frac{du}{u^{1+1/\lambda}}\right.\\
&\qquad\left.+\exp(G_{1})\int_0^{+\infty}\big(M_1(-u)-M_1(0)\big)\frac{du}{u^{1+1/\lambda}}
\right)\\
&=\bar\kappa_1\Delta_0+\kappa_1F_{1},
\end{align*}
where in the second equality we use that $G_{1}^+-G_{1}^-=G_{1}$, in the third one we replace $\Delta_0^+$ by $\Delta_0+\Delta_0^-$ to obtain a function $\bar\kappa_1$ multiplying $\Delta_0$ and in the last one we set $\kappa_1=\lambda\Delta_0^+\exp(-G_{2}^+)$. Therefore $\Delta_1=\kappa_1F_{1}+\bar\kappa_1\Delta_0$ with $\kappa_1(\mu_0)>0$. Since one can easily verify that $\bar\kappa_1$ is analytic at $\mu_0$ with $\lambda(\mu_0)>1,$
this concludes the proof of assertion $(1)$.

It only remains to compute $\Delta_3=\Delta_3^+-\Delta_3^-$ in case that $\lambda(\mu)=1.$ Exactly as before, since $(x,y)\mapsto (-x,y)$ sends $-X_\mu$ to $\tilde X_\mu$, from the expression of $\Delta_3^+$ in \refc{teo1eq4.5} and taking \refc{teo1eq5} into account we get
\[
\left.\Delta_3^-\right|_{\lambda=1}=(\Delta_0^-)^2\exp(-G_2^-)\bar M_2'(0)= (\Delta_0^-)^2\exp(-G_2^-)M_2'(0).
\] 
Hence some straightforward computations show that 
\begin{align*}
\left.\Delta_3\right|_{\lambda=1}&=\left((\Delta_0^+)^2\exp(-G_2^+)-(\Delta_0^-)^2\exp(-G_2^-)\right)M_2'(0)\\[5pt]
&=\left((\Delta_0+\Delta_0^-)^2\exp(-G_2-G_2^-)-(\Delta_0^-)^2\exp(-G_2^-)\right)M_2'(0)\\[5pt]
&=-\kappa_3G_2M_2'(0)+\bar\kappa_3\Delta_0,
\end{align*}
where  
\[
 \kappa_3\!:=(\Delta_0^-)^2\exp(-G_2^-)\frac{1-\exp(-G_2)}{G_2}\text{ and }
 \bar\kappa_3\!:=(\Delta_0+2\Delta_0^-)\exp(-G_2^+)M_2'(0),
\] 
which are analytic functions at $\mu_0$ and $\kappa_3(\mu_0)>0.$ Finally, due to
$M_2(u)=L_2(u)\partial_2K(u,0)$ with
\[
 M_2(u)=L_2(u)\partial_2K(u,0)\text{ and }L_2(u)=\exp\int_0^u\big(K(z,0)+{\lambda}\big)\frac{dz}{z},
\]
one can easily show that 
$M_2'(0)=L_2'(0)\partial_2K(0,0)+L_2(0)\partial_{12}K(0,0)=\partial_1K(0,0)\partial_2K(0,0)+\partial_{12}K(0,0)$. We thus obtain that $\Delta_3|_{\lambda=1}=\kappa_3F_3+\bar\kappa_3\Delta_0$ with $F_3=-G_2\big(\partial_1K\partial_2K+\partial_{12}K\big)(0,0),$ as desired. This proves the validity of the third assertion in the statement and concludes the proof of the result.
\end{prooftext}

\subsection{Proof of \propc{prop2}}\label{ApB2}

In this section we prove \propc{prop2}, which gives the asymptotic development of the difference map 
\[
 \mathscr D_u(s;\mu)\!:=D^u_+(s;\mu)-D^u_-(s;\mu),
\]
see \figc{dib5}, together with some properties of its coefficients. To this end we need first two auxiliary results. 

\begin{lem}\label{lem1} 
Fix any $\mu_0=(a_0,b_0,\varepsilon_0,\varepsilon_1,\varepsilon_2)$ with $(a_0,b_0)\in (-2,0)\times (0,2)$ and $\varepsilon_i\approx 0$ for $i=0,1,2.$ Then 
\[
 D^u_\pm(s;{\mu})=\delta_\pm+ \Delta_{0}^\pm s^\lambda+\F_{\ell}^\infty(\mu_0),
 \text{ for any ${\ell}\in \big[\lambda_0,\min(2\lambda_0,\lambda_0+1)\big)$,}
\]
where $\lambda$, $\delta_\pm$ and $\Delta_0^\pm$ are $\cc^\infty$ functions on $\mu\approx\mu_0$ and $\lambda_0\!:=\lambda(\mu_0)=-\frac{a_0+2}{a_0}.$ Moreover, for $a_0\neq -1,$

\begin{enumerate}[$(1)$]

\item If $a_0>-1$ then $D^u_\pm(s;{\mu})=\delta_\pm+ \Delta_{0}^\pm s^\lambda +\Delta_{1}^\pm s^{\lambda+1} + \F_{\ell}^\infty(\mu_0)$ for any ${\ell}\in \big[\lambda_0+1,\min(2\lambda_0,\lambda_0+2)\big)$, where
$\Delta_1^\pm$ is a $\cc^\infty$ function on $\mu\approx\mu_0.$ 

\item If $a_0<-1$ then $D^u_\pm(s;{\mu})=\delta_\pm+\Delta_{0}^\pm s^\lambda +\Delta_{2}^\pm s^{2\lambda}+ \F_{\ell}^\infty(\mu_0)$ for any ${\ell}\in\big[2\lambda_0,\min(3\lambda_0,\lambda_0+1)\big)$, where
$\Delta_2^\pm$ is a $\cc^\infty$ function on $\mu\approx\mu_0.$ 

\end{enumerate}
\end{lem}

\begin{prova}
For the sake of simplicity in the exposition we omit the superscript in $D^u_\pm$. That being said, let us prove the result for the Dulac map $D_+(\,\cdot\,;\mu)$, the proof for $D_-(\,\cdot\,;\mu)$ follows verbatim. We denote the $y$-coordinate of the intersection point with $x=0$ of the unstable separatrix of the saddle at~$s_1$ by $\delta_+(\mu)$. 
The function $\mu\mapsto \delta_+(\mu)$ is $\cc^\infty$ in a neighbourhood of $\mu=\mu_0$.
Indeed, this follows by first applying the local center-stable manifold theorem (see \cite[Theorem 1]{Kelley} for instance) to $s_1$ and then appealing to the smooth dependence of the solutions of $X_\mu$ on initial conditions and parameters.
It is clear moreover that $\delta_+|_{\varepsilon_0=0}\equiv 0.$ For convenience we change the parametrisation on $\Sigma_2$ by $\hat s\mapsto \big(0,\hat s+\delta_+(\mu)\big)$ for $\hat s>0$ small enough and we denote by $\hat D_+(s;\mu)$ the Dulac map of $X_\mu$ from $\Sigma_1$ to $\Sigma_2$ with this new parametrisation in the arrival section. It is then clear that $D_+(s;\mu)=\delta_+(\mu)+\hat D_+(s;\mu)$ for $s>0.$
To study $\hat D_+(\,\cdot\,;\mu)$ we first compactify the vector field $X_\mu$ by using the projective coordinates $(u,v)=\phi_1(x,y)\!:=(\frac{1}{x+y+1},\frac{y}{x+y+1}).$ The key point here is that the trajectories of $X_\mu$ from $\Sigma_1$ to~$\Sigma_2$ do not intersect $x+y+1=0$. In doing so we obtain an analytic family of vector fields which is orbitally equivalent to a polynomial one, say $Y_\mu,$ that has a finite hyperbolic saddle at the origin. By construction its stable separatrix is at $u=0$ for all $\mu,$ whereas its unstable one is at $v=0$ only when $\varepsilon_0=0.$ In order to straighten both separatrices for all $\mu$ we apply \lemc{rectes}, that gives a $\cc^\infty$ family of diffeomorphisms $\phi_2(u,v;\mu)$ such that the push-forward $(\phi_2)_\star\big(Y_\mu\big)$ writes as in~\refc{X} with $\varpi=\infty.$ By construction, setting $\phi\:=\phi_2\circ\phi_1$, its Dulac map from $\phi(\Sigma_1)$ to $\phi(\Sigma_2)$, parametrised, respectively, by $\sigma_1(s;\mu)=\phi(0,1/s;\mu)$ and $\sigma_2(s;\mu)=\phi(0,s+\delta_+(\mu);\mu)$, is precisely $\hat D(s;\mu).$ Observe in this regard that the parametrisations of the transverse sections are $\cc^\infty.$ Accordingly, by applying \propc{3punts},
\[
 \hat D_+(s;{\mu})=\Delta_{0}^+({\mu})s^\lambda+
\left\{
\begin{array}{ll}
\Delta_{1}^+({\mu})s^{\lambda+1}+  \F_{\ell_1}^\infty(\mu_0)&\text{ if $\lambda_0>1,$}\\[10pt]
\Delta_{2}^+({\mu})s^{2\lambda}+ \F_{\ell_2}^\infty(\mu_0)&\text{ if $\lambda_0<1,$}
\end{array}
\right.
\]
for any ${\ell_2}\in  \big[2\lambda_0,\min(3\lambda_0,\lambda_0+1)\big)$ and ${\ell_1}\in \big[\lambda_0+1,\min(2\lambda_0,\lambda_0+2)\big)$. Here $\lambda=\lambda(\mu)$ is the hyperbolicity ratio of the saddle of $X_\mu$ at $s_1$ and
$\lambda_0=\lambda(\mu_0)=-\frac{a_0+2}{a_0}$. Moreover the coefficient $\Delta_0^+$ is~$\cc^\infty$ at $\mu_0$ and, on the other hand, the coefficient $\Delta_1^+$  (respectively, $\Delta_2^+$) is $\cc^\infty$ at $\mu_0$ provided that $\lambda_0>1$ (respectively, $\lambda_0<1$). On account of $D_+(s;\mu)=\delta_+(\mu)+\hat D_+(s;\mu)$ this concludes the proof of the result.
\end{prova}

\begin{lem}\label{delta}
$\partial_{\varepsilon_0}\big(\delta_+-\delta_-\big)(\mu)> 0$ for all $\mu=(a,b,0,0,0)$ with $a\in (-2,0)\setminus\{-1\}$ and $b\in (0,2).$ 
\end{lem}

\begin{prova}
The differential form associated to system \refc{pert} is given by
\[
 \textstyle\Omega\!:=\big(2xy-\varepsilon_0\big)dx+\big(\frac{b-2}{4}+\varepsilon_1x+(1-b)y+ax^2+\varepsilon_2xy+by^2\big)dy.
\]
We know on the other hand that 
\begin{equation}\label{lem0eq1}
 H(x,y)\!:=y(x^2+\ell y^2+my+n)^{\frac{1}{a}},
\end{equation}
with $\ell=\frac{b}{a+2},$ $m=-\frac{b-1}{a+1}$ and $n=\frac{b-2}{4a}$,
is a first integral of \refc{pert} for $\varepsilon_0=\varepsilon_1=\varepsilon_2=0$. We observe in this regard that
\[
 a\frac{dH}{H}=a\frac{dy}{y}+\frac{2xdx+(2\ell y+m)dy}{x^2+\ell y^2+my+n},
\]
which yields
\[
ay^{1-a}H^{a-1}dH=2xy\,dx+
\big(an+(a+1)my+ax^2+(a+2)\ell y^2\big)dy=\Omega|_{\varepsilon_1=\varepsilon_2=0}+\varepsilon_0dx,
\]
where in the second equality we use the expression of $\ell$, $m$ and $n$ in terms of $a$ and $b.$ This shows that $\Omega|_{\varepsilon_1=\varepsilon_2=0}$ is proportional to $\Omega_0\!:=dH-\frac{\varepsilon_0}{a}H^{1-a}y^{a-1}dx.$ On account of this, if we take any $\mu_0=(a,b,\varepsilon_0,0,0)$ and denote by $\Gamma_{s,\varepsilon_0}$ the oriented arc of orbit of $X_{\mu_0}$ that joins the points $\big(0,D^u_+(s;\mu_0)\big)$ and $\big(0,D^u_-(s;\mu_0)\big)$ then we have that
\[
 0=\int_{\Gamma_{s,\varepsilon_0}}\Omega_0
 =H\big(0,D^u_+(s;\mu_0)\big)-H\big(0,D^u_-(s;\mu_0)\big)-\frac{\varepsilon_0}{a}\int_{\Gamma_{s,\varepsilon_0}}H(x,y)^{1-a}y^{a-1}dx,
\]
where $D^u_\pm(s;\mu)$ is Dulac map in \lemc{lem1}. Consequently 
\[
H\big(0,D^u_+(s;\mu_0)\big)-H\big(0,D^u_-(s;\mu_0)\big)=\frac{\varepsilon_0}{a}\int_{\Gamma_{s,\varepsilon_0}}H(x,y)^{1-a}y^{a-1}dx\text{ for all $\varepsilon_0\approx 0.$}
\]
The derivative of this expression with respect to $\varepsilon_0$ evaluated at $\bar\mu_0\!:=(a,b,0,0,0)$ yields
\[
\partial_yH\big(0,D^u_+(s;\bar\mu_0)\big)\partial_{\varepsilon_0}D^u_+(s;\bar\mu_0)
-\partial_yH\big(0,D^u_-(s;\bar\mu_0)\big)\partial_{\varepsilon_0}D^u_-(s;\bar\mu_0)=\frac{1}{a}\int_{\Gamma_{s,0}}H(x,y)^{1-a}y^{a-1}dx.
\]
Our next goal will be to make $s\to 0^+$ in this equality. With this aim in view note that, by the first assertion in \lemc{lem1}, $D^u_\pm(s;\mu)=\delta_\pm(\mu)+\F_\rho^\infty(\mu_0)$ for any $\rho>0$ small enough. Consequently, since $\delta_\pm(\bar\mu_0)=0,$ we get that
\[
 \lim_{s\to 0^+}\partial_yH\big(0,D^u_\pm(s;\bar\mu_0)\big)\partial_{\varepsilon_0}D^u_\pm(s;\bar\mu_0)=
 \partial_yH(0,0)\partial_{\varepsilon_0}\delta_\pm(\bar\mu_0)=n^{1/a}\partial_{\varepsilon_0}\delta_\pm(\bar\mu_0),
\]
where in the first equality we use the good properties of the remainder with respect to the derivation of the parameters, see \defic{defi2}, and in the second one the expression in \refc{lem0eq1}. Therefore 
\begin{equation}\label{lem0eq2}
 an^{1/a}\big(\partial_{\varepsilon_0}\delta_+(\bar\mu_0)-\partial_{\varepsilon_0}\delta_-(\bar\mu_0)\big)=
 \lim_{s\to 0^+}\int_{\Gamma_{s,0}}H(x,y)^{1-a}y^{a-1}dx.
\end{equation}
Note at this point that $\Gamma_{s,0}$ is a periodic orbit of $X_{\bar\mu_0}.$ Thus it is contained inside the level set $H(x,y)=h$ where $h=h(s)$ verifies
\[
 h=H(0,1/s)=s^{-1-2/a}(\ell+ms+ns^2)^{\frac{1}{a}}. 
\]
Here we use \refc{lem0eq1} once again and that the parametrization of $\Sigma_1$ is given by $s\mapsto (0,1/s).$ Since $a\in (-2,0)$ by assumption, this shows that $\lim_{s\to 0^+}h(s)=0.$ Accordingly, if we denote by $\gamma_h$ the periodic orbit of $X_{\bar\mu_0}$ inside the level curve $H=h$, from \refc{lem0eq2} we get that
\[
an^{1/a}\partial_{\varepsilon_0}\big(\delta_+-\delta_-\big)(\bar\mu_0)=
 \lim_{h\to 0}h^{1-a}\int_{\gamma_h}y^{a-1}dx.
\]
It is clear then that the result will follow once we prove that the above limit exists and is different from zero. To this end, setting $\gamma_h^+\!:=\gamma_h\cap\{x\geqslant 0\}$, we first observe that
\[
 \int_{\gamma_h}y^{a-1}dx=2\int_{\gamma^+_h}y^{a-1}dx
\]
since~$X_{\bar\mu_0}$ is symmetric with respect to $x=0$.
To compute this Abelian integral we perform the projective change of coordinates $(u,v)=(\frac{1}{x},\frac{y}{x})$ and in these new variables, see \refc{lem0eq1}, we have that
\[
\gamma_h^+\subset\{\hat H(u,v)=h^a\},\text{ where $\hat H(u,v)\!:=u^{-a-2}v^{a}({1+\ell v^2+muv+nu^2}).$}
\]
A computation shows that
\begin{equation}\label{lem0eq3}
 \frac{\partial_u\hat H(u,v)}{\partial_v\hat H(u,v)}=-\frac{u}{v}\frac{(u-2v)\big((b-2)u-2bv\big)+4a}{(u-2v)\big((b-2)u-2bv\big)+4(a+2)},
\end{equation}
which gives, up to a unity, the expression of the partial derivatives of $\hat H.$ Then, taking $(a,b)\in (-2,0)\times (0,2)$ into account, it follows that $\partial_v\hat H(u,v)\neq 0$ on $0<u\leqslant 2v$ and $\partial_u\hat H(u,v)\neq 0$ on $0<2v\leqslant u.$ Observe also that, for each $h>0$, the arc $\gamma_h^+$ has exactly one intersection point with the straight line $u=2v$ because $\hat H(u,u)=h^a$ if, and only if, $u=\pm c(h)$ where $c(h)\!:=(2^{a+2}h^a-(\ell+2m+4n))^{-1/2}.$ Therefore, by applying (twice) the Implicit Function Theorem to $\hat H(u,v)=h^a$ we can split $\gamma_h^+$ as
\[ 
 \gamma_h^+=\big\{u=u(v;h), v\in[c(h),+\infty)\big\}\cup\big\{v=v(u;h), u\in[c(h),+\infty)\big\}.
\]
Accordingly, from \refc{lem0eq1} once again,
\begin{align*}
 \lim_{h\to 0}h^{1-a}\int_{\gamma_h^+}y^{a-1}dx&=\lim_{h\to 0}\int_{\gamma_h^+}(x^2+\ell y^2+my+n)^{\frac{1}{a}-1}dx\\
&=- \lim_{h\to 0}\left(\int_{c(h)}^{+\infty} \left.(1+\ell v^2+mu v+nu^2)^{\frac{1-a}{a}}\right|_{v=v(u;h)}u^{-\frac{2}{a}}du\right.
\\
&\hspace{1.3truecm}\left.-\int_{c(h)}^{+\infty} \left.(1+\ell v^2+mvu +nu^2)^{\frac{1-a}{a}}u^{-\frac{2}{a}}\right|_{u=u(v;h)}\partial_v u(v;h)dv\right).
\end{align*}
In order to make this limit let us first observe that $\lim_{h\to 0}c(h)= 0$ due to $a<0.$ On the other hand, $\lim_{h\to 0} u(v;h)=0$, uniformly in~$v$, and $\lim_{h\to 0} v(u;h)=0$, uniformly in~$u$, because the oval $\gamma_h$ tends to the polycycle (in  Hausdorff sense) as $h\to 0.$ Furthermore, due to
\[
\partial_vu(v;h)=\frac{du}{dv}=-\left.\frac{\partial_v\hat H(u,v)}{\partial_u\hat H(u,v)}\right|_{u=u(v;h)},
\]
from the expression in \refc{lem0eq3} we deduce that $|\partial_vu(v;h)|$ is uniformly bounded since $0<u(v;h)\leqslant 2v$ for any $v\in[c(h),+\infty)$. Taking these facts into account, together with the assumption $a\in(-2,0)$, by applying the Dominated Convergence Theorem we conclude that 
\[
\lim_{h\to 0}h^{1-a}\int_{\gamma_h^+}y^{a-1}dx=-\int_0^{+\infty}(1+nu^2)^{\frac{1-a}{a}}u^{-\frac{2}{a}}du=:\!p\in\R_{<0}.
\]
Hence$\partial_{\varepsilon_0}\big(\delta_+-\delta_-\big)(\bar\mu_0)=\frac{2pn^{-1/a}}{a}>0$ and this finishes the proof of the result.
\end{prova}

\begin{prooftext}{Proof of \propc{prop2}.}
The three assertions with regard to structure of the asymptotic development follow from \lemc{lem1} setting $\Delta^{u}_i\!:=\Delta_i^+-\Delta_i^-$ for $i=0,1,2$ and $\delta_u\!:=\delta_+-\delta_-$ because then
\begin{equation}\label{lem2eq0}
 \mathscr D_{u}(s;\mu)=D_+(s;\mu)-D_-(s;\mu)=\delta_{u}(\mu)+\Delta^{u}_{0}({\mu})s^\lambda+
\left\{
\begin{array}{ll}
\Delta^{u}_{1}({\mu})s^{\lambda+1}+ \F_{\ell_1}^\infty(\mu_0)&\text{ if $a_0>-1,$}\\[10pt]
\Delta^{u}_{2}({\mu})s^{2\lambda}+ \F_{\ell_2}^\infty(\mu_0)&\text{ if $a_0<-1,$}
\end{array}
\right.
\end{equation}
for any ${\ell_2}\in  \big[2\lambda_0,\min(3\lambda_0,\lambda_0+1)\big)$ and ${\ell_1}\in \big[\lambda_0+1,\min(2\lambda_0,\lambda_0+2)\big)$.
Since we will deal with the ``upper case'' only, for simplicity in the exposition we shall omit any subscript and superscript $u$ from now on. 

It is clear that $\mathscr D(s;{\mu_{0}})\equiv 0$ because $X_\mu$ is inside the center variety when $\mu=\mu_{0}$. On the other hand, by \lemc{delta}, $\partial_{\varepsilon_0}\delta(\mu_{0})>0$. Note also that the straight line $y=0$ is invariant in case that $\varepsilon_0=0.$ Hence $\delta(\mu)|_{\varepsilon_0=0}\equiv 0$ by definition and, consequently,
$\partial_{\varepsilon_1}\delta(\mu_{0})=\partial_{\varepsilon_2}\delta(\mu_{0})=0.$ That being stablished, our main task is to compute the partial derivatives $\partial_{\varepsilon_1}\Delta_k$ and $\partial_{\varepsilon_2}\Delta_k$ evaluated at $\mu_{0}=(a_0,b_0,0,0,0)$ for each $k=0,1,2.$ To this end the key point is that we can perform the computations setting $\varepsilon_0=0$ and that in this case $X_\mu$ is a D-system, more concretely, with $f(x,y)=1-b+\varepsilon_2 x+by,$ $g(x)=\frac{b-2}{4}+\varepsilon_1x+ax^2$, $q(x,y)=-2x$ and $n=1,$ so that
\[
 \ell_2(x,y)=(a+2)x^2+\varepsilon_2xy+by^2.
\]
Let us remark that it is only for $\varepsilon_0=0$ that $X_{\mu}$ becomes a D-system. Thus, for the sake of consistency we shall denote $\bar\mu=(a,b,\varepsilon_1,\varepsilon_2)$ and $\bar\mu_0=(a,b,0,0)$. That being said, following the notation in \teoc{teo1}, that we stress it is addressed to D-systems, from \refc{K} we have
\[
K(x_1,x_2;\bar\mu)=\left. 1-\frac{xq(x,y)}{yf(x,y)+g(x)}\right|_{(x,y)=\left(\frac{1}{x_1},\frac{x_2}{x_1}\right)}
=1+\frac{2}{a+\varepsilon_1x_1+\varepsilon_2x_2+\frac{b-2}{4}x_1^2+(1-b)x_1x_2+bx_2^2}.
\]
Hence $\lambda(\bar\mu)=-\frac{a+2}{a}.$ From \refc{d0} we get that
\begin{align*}
 d_0(\bar\mu)&=2\int_{-\infty}^{+\infty}\left(\frac{z}{\frac{b-2}{4}+\varepsilon_1z+az^2}+\lambda\frac{z}{(a+2)z^2+\varepsilon_2z+b}\right)dz\\[4pt]
 &=\frac{2\pi}{a}\left(\frac{\varepsilon_1}{\sqrt{(b-2)a-\varepsilon_1^2}}+\frac{\varepsilon_2}{\sqrt{4b(a+2)-\varepsilon_2^2}}\right).
\end{align*}
On account of this one can verify that $d_0(\bar\mu)=-\rho_0(\bar\mu)\Big(2\frac{\sqrt{b(a+2)}}{\sqrt{a(b-2)}}\,\varepsilon_1+\varepsilon_2\Big)$ where $\rho_0$ is a smooth function with $\rho_0(\bar\mu_{0})>0$ since $a_0\in (-2,0).$ Hence, from \refc{lem2eq0} and applying \teoc{teo1}, 
\begin{equation*}
\left.\Delta_0(\mu)\right|_{\varepsilon_0=0}=-\kappa_{01}(\bar\mu)\left(2\frac{\sqrt{b(a+2)}}{\sqrt{a(b-2)}}\,\varepsilon_1+\varepsilon_2\right)\text{ with $\kappa_{01}(\bar\mu_{0})> 0.$}
\end{equation*}
Consequently, there exists a smooth function $\rho_1=\rho_1(\mu)$ such that
\[
 \Delta_0(\mu)=-\kappa_{01}(\bar\mu)\left(2\frac{\sqrt{b(a+2)}}{\sqrt{a(b-2)}}\,\varepsilon_1+\varepsilon_2\right)
 +\varepsilon_0\rho_1(\mu)
 =-\kappa_{01}(\bar\mu)\left(2\frac{\sqrt{b(a+2)}}{\sqrt{a(b-2)}}\,\varepsilon_1+\varepsilon_2\right)
 +\kappa_{02}(\mu)\delta(\mu),
\]
where in the second equality we use that we can write $\delta(\mu)=\varepsilon_0\rho_2(\mu)$ with $\rho_2(\mu_0)\neq 0$ due to $\delta(\mu)|_{\varepsilon_0=0}\equiv0 $ and $\partial_{\varepsilon_0}\delta(\mu_{0})\neq 0$. Since $\kappa_{01}(\bar\mu)$ is a smooth function on $\mu$, this proves the assertion with regard to $\Delta_0(\mu).$

Let us assume now that $a_0\in (-2,-1)$ and turn to the study of $\Delta_2$. This, on account of \teoc{teo1}, leads to the computation of $F_{2}$. According to \refc{F2} its expression is given by
\begin{equation}\label{lem2eq0bis}
F_{2}(\bar\mu)=\int_0^{+\infty} \Big(M_2(-z)-M_2(0)+\exp(G_{2})\big(M_2(z)-M_2(0)\big)\Big)\frac{dz}{z^{1+\lambda}}
\end{equation}
where $M_2(u)=L_2(u)\partial_2K(u,0)$ with $L_2(u)\!:=\exp\left(\int_0^u\big(K(z,0)+\lambda\big)\frac{dz}{z}\right).$
After some lengthy computations we obtain that
\[
L_2(u)=\left(1+\frac{\varepsilon_1}{a}u+{\eta_2} u^2\right)^{-\frac{1}{a}}B_2(u),
\]
where ${\eta_2}\!:=\frac{b-2}{4a}>0$ for all $(a,b)\in(-2,0)\times(0,2)$ and
\[
B_2(u)\!:=\exp\left(\frac{-2\varepsilon_1}{a\sqrt{a(b-2)-\varepsilon_1^2}}
\left(
\arctan\left(\frac{\frac{b-2}{2}u+\varepsilon_1}{\sqrt{a \left(b -2\right)-\varepsilon_{1}^{2}}}\right)
-\arctan\left(\frac{\varepsilon_1}{\sqrt{a(b-2)-\varepsilon_1^2}}\right)
\right)
\right).
\]
The explicit computation of $F_2(\bar\mu)$ for arbitrary $\bar\mu$ requires a primitive of 
$u\mapsto (M_2(u)-M_2(0))\,u^{-1-\lambda},$ which is not feasible because 
$M_2(u)=L_2(u)\partial_2K(u,0)$
where
\[
\partial_2K(u,0)=\frac{2}{a^2}\frac{(b-1)u-\varepsilon_2}{\big(1+\frac{\varepsilon_1}{a}u+{\eta_2} u^2\big)^2}.
\]
To bypass this problem the strategy is to compute only the first order Taylor's expansion of this function at $(\varepsilon_1,\varepsilon_2)=(0,0).$ In doing so we get
\begin{align*}
M_2(u)-M_2(0)=&\,\frac{2(b-1)}{a^2}\,u(1+{\eta_2} u^2)^{-2-\frac{1}{a}}\\
&-\frac{2(b-1)}{a^4\sqrt{{\eta_2}}}\,u(1+{\eta_2} u^2)^{-3-\frac{1}{a}}\big((1+{\eta_2} u^2)\arctan(\sqrt{{\eta_2}}u)+\sqrt{{\eta_2}}(1+2a)u\big)\varepsilon_1\\ 
&-\frac{2}{a^2}\big((1+{\eta_2} u^2)^{-2-\frac{1}{a}}-1\big)\varepsilon_2+\op(\|(\varepsilon_1,\varepsilon_2)\|).
\end{align*}
Thus, on account of the parity of each coefficient with respect to $u,$ if we write
\begin{equation}\label{lem2eq4}
 \int_0^{+\infty} \big(M_2(\pm u)-M_2(0)\big)\frac{du}{u^{1+\lambda}}
 =m_0^\pm+m_1^\pm\varepsilon_1+m_2^\pm\varepsilon_2
 +\op(\|(\varepsilon_1,\varepsilon_2)\|)
\end{equation}
then it turns out that $m_0^-=-m_0^+$, $m_1^-=m_1^+$ and $m_2^-=m_2^+$. Of course to obtain the above equality we must prove that the higher order terms can also be neglected after integration. To show this let us note first that, as a matter of fact, the higher order terms do not depend on $\varepsilon_2$ because $M_2(u;\bar\mu)$ is linear in this parameter. Therefore to get $m_0^\pm$ we need a result to pass the limit $\varepsilon_1\to 0$ under the integral sign, and to get $m_1^\pm$ a similar result for the derivation with respect to $\varepsilon_1$. With this aim we appeal to the results in \cite[\S 17.2]{Zorich} about improper integrals depending on a parameter. More concretely, Proposition~2, which is a sort of Weierstrass test for the uniform convergence of an improper integral depending on a parameter, and Proposition 6, that gives sufficient conditions for the differentiation of an improper integral with respect to a parameter. 
To this end the key points are that, on one hand, $\lambda=\lambda(\bar\mu)=-\frac{a+2}{a}\in (0,1)$ for $\bar\mu\approx\bar\mu_{0}$ due to $a_0\in (-2,-1)$ and, on the other hand, that $B_2(u;\bar\mu)$ and $\partial_{\varepsilon_1}B_2(u;\bar\mu)$ are bounded for $u\in (0,+\infty)$ and $\varepsilon_1\approx 0$ by a constant. That being said, some computations show that
\begin{align*}
 m_0^+=&
% \frac{2(b-1)}{a^2}\int_0^\infty u(1+{\eta_2} u^2)^{-2-1/a}u^{-\lambda-1}du=
 \frac{2(b-1)}{a^2}\int_0^{+\infty}(1+{\eta_2} u^2)^{-2-\frac{1}{a}} u^{\frac{2}{a}+1}du=\frac{(b-1){\eta_2}^{-1-\frac{1}{a}}}{a(a+1)}
 \intertext{and}
 m_2^+=&-\frac{2}{a^2}\int_0^{+\infty}\big((1+{\eta_2} u^2)^{-2-\frac{1}{a}}-1\big)u^{\frac{2}{a}}du
 =-\frac{\sqrt{\pi}}{2a^2}\frac{\Gamma\left(\frac{a+2}{2a}\right)}{\Gamma\left(\frac{2a+1}{a}\right)}{\eta_2}^{-\frac{a+2}{2a}}.
\end{align*}
One can readily check in particular that $m_2^+>0$ for all $(a,b)\in (-2,-1)\times (0,2).$
Computing $m_1^+$ is a little more involved. In this case
\begin{align*}
-\frac{a^4}{2(b-1)}m_1^+=&\frac{1}{\sqrt{{\eta_2}}}\int_0^{+\infty}(1+{\eta_2} u^2)^{-2-\frac{1}{a}}\arctan(\sqrt{{\eta_2}}u)u^{1+\frac{2}{a}}du+(1+2a)\int_0^{+\infty}(1+{\eta_2} u^2)^{-3-\frac{1}{a}}u^{2+\frac{2}{a}}du\\
=&
\frac{a\pi{\eta_2}^{-\frac{3a+2}{2a}}}{4(1+a)}-\frac{a\sqrt{\pi}}{4(a+1)}\frac{\Gamma\left(\frac{3a+2}{2a}\right)}{\Gamma\left(\frac{2a+1}{a}\right)}{\eta_2}^{-\frac{3a+2}{2a}}
+(1+2a)\frac{\sqrt{\pi}}{4}\frac{\Gamma\left(\frac{3a+2}{2a}\right)}{\Gamma\left(\frac{3a+1}{a}\right)}{\eta_2}^{-\frac{3a+2}{2a}},\\
%=&\frac{a\sqrt{\pi}}{4}\left(-\frac{{\eta_2}}{a+1}+1\right)\frac{\Gamma((3a+2)/(2a))}{\Gamma((2a+1)/a)}{\eta_2}^{-3/2-1/a}
\end{align*}
and after some simplifications we get that
\[
 m_1^+=-\frac{\sqrt{\pi}(b-1)}{2a^2(a+1)}\left(\frac{\Gamma\left(\frac{3a+2}{2a}\right)}{\Gamma\left(\frac{2a+1}{a}\right)}+\frac{\sqrt{\pi}}{a}\right){\eta_2}^{-\frac{3a+2}{2a}}.
\]
On the other hand, 
 \[
  G_2=\int_0^{+\infty}\left(\frac{q(u,0)}{g(u)}+\frac{q(-u,0)}{g(-u)}\right)du=-\frac{2\pi\varepsilon_1}{a\sqrt{(b-2)a-\varepsilon_1^2}},
\]
so that $\exp(G_2)=1+\frac{2\pi}{\sqrt{a^3(b-2)}}\varepsilon_1+\op(\varepsilon_1)$ due to $a<0.$ Accordingly the substitution of  \refc{lem2eq4} in \refc{lem2eq0bis} yields
\begin{align*}
F_2&= m_0^-+m_1^-\varepsilon_1+m_2^-\varepsilon_2+\left(1+\frac{2\pi}{\sqrt{a^3(b-2)}}\varepsilon_1\right)\left(m_0^++m_1^+\varepsilon_1+m_2^+\varepsilon_2\right)+\op(\|(\varepsilon_1,\varepsilon_2)\|)\\
&=m_0^-+m_0^++\left(m_1^-+m_1^++\frac{2\pi}{\sqrt{a^3(b-2)}}m_0^+\right)\varepsilon_1+\left(m_2^-+m_2^+\right)\varepsilon_2+\op(\|(\varepsilon_1,\varepsilon_2)\|)\\
&=2\left(m_1^++\frac{\pi}{\sqrt{a^3(b-2)}}m_0^+\right)\varepsilon_1+2m_2^+\varepsilon_2+\op(\|(\varepsilon_1,\varepsilon_2)\|)
\end{align*}
and let us note that
\[
 m_1^++\frac{\pi}{\sqrt{a^3(b-2)}} m_0^+=
%\frac{\sqrt{\pi}(b-1){\eta_2}^{-1-\frac{1}{a}}}{a^2}\frac{\Gamma\left(\frac{3}{2}+\frac{2}{a}\right)}{\Gamma\left(2+\frac{1}{a}\right)}\left[\frac{1}{(a+1)\sqrt{a(b-2)}}-\frac{{\eta_2}^{-1/2-1/a}}{2a}\right]
-\frac{\sqrt{\pi}}{2a^2}\frac{b-1}{a+1}\frac{\Gamma\left(\frac{3a+2}{2a}\right)}{\Gamma\left(\frac{2a+1}{a}\right)}{\eta_2}^{-\frac{3a+2}{2a}}.
\]
Hence $m_1^++\frac{\pi}{\sqrt{a^3(b-2)}} m_0^+=\frac{2(a+2)(b-1)}{(a+1)(b-2)}m_2^+,$ so that 
\[
 F_2(\bar\mu)=\rho_3(\bar\mu)\left(\frac{2(a+2)(b-1)}{(a+1)(b-2)}\varepsilon_1+\varepsilon_2+\op(\|(\varepsilon_1,\varepsilon_2)\|)\right)\text{ with $\rho_3(\bar\mu_0)>0.$}
 \]
Finally, from \refc{lem2eq0} and the last assertion in $(2)$ of \teoc{teo1} we get that
\[
 \left.\Delta_2(\mu)\right|_{\varepsilon_0=0}=\kappa_{21}(\bar\mu)\left( \frac{2(a+2)(b-1)}{(a+1)(b-2)}\varepsilon_1+\varepsilon_2+\op(\|(\varepsilon_1,\varepsilon_2)\|)\right)+\kappa_{22}(\bar\mu)\left.\Delta_0(\mu)\right|_{\varepsilon_0=0} \text{ with $\kappa_{21}(\bar\mu_{0})> 0.$}
\] 
Let us stress here that $\kappa_{21}$ and $\kappa_{22}$ are smooth functions in a neighbourhood of $\bar\mu_{0}$ provided that $a_0\in (-2,-1).$ Consequently, since $\delta(\mu)|_{\varepsilon_0=0}\equiv 0$ and $\partial_{\varepsilon_0}\delta(\mu_{0})\neq 0$, we get that
\[
 \Delta_2(\mu) =\kappa_{21}(\bar\mu)\left( \frac{2(a+2)(b-1)}{(a+1)(b-2)}\varepsilon_1+\varepsilon_2+\op(\|(\varepsilon_1,\varepsilon_2)\|)\right)+\kappa_{22}(\bar\mu)\Delta_0(\mu)+\kappa_{23}(\mu)\delta(\mu)
\] 
for some smooth function $\kappa_{23}$ and this proves the assertion in $(2).$ 

So far we have studied the coefficient $\Delta_2$ assuming $a_0\in (-2,-1),$ i.e., $\lambda_0<1.$ Our next task is to do the same for the coefficient $\Delta_1$ assuming $a_0\in (-1,0),$ i.e., $\lambda_0>1.$ In this case, see \refc{F1}, we have to compute 
\begin{equation}\label{lem2eq1bis}
 F_{1}(\mu)=-\int_0^{+\infty} \Big(M_1(z)-M_1(0)
   +\exp(G_{1})\big(M_1(-z)-M_1(0)\big)\Big)\frac{dz}{z^{1+1/\lambda}},
\end{equation}
where $M_1(u)=L_1(u)\partial_1\big(\frac{1}{K}\big)(0,u)$ with $L_1(u)\!:=\exp\left(\int_0^u\left(\frac{1}{K(0,z)}+\frac{1}{\lambda}\right)\frac{dz}{z}\right).$ In doing so exactly as before we obtain that 
\[
 L_1(u)=\left(1+\frac{\varepsilon_2}{a+2}u+\eta_1u^2\right)^{\frac{1}{a+2}}B_1(u),
\]
where $\eta_1\!:=\frac{b}{a+2}>0$ for all $(a,b)\in(-2,0)\times(0,2)$ and
\[
B_1(u)\!:=\exp\left(\frac{2\varepsilon_2}{(a+2)\sqrt{4b(a+2)-\varepsilon_2^2}}
\left(
\arctan\left(
\frac{2bu+\varepsilon_2}{\sqrt{4b(a+2)-\varepsilon_2^2}}
\right)
-
\arctan\left(
\frac{\varepsilon_2}{\sqrt{4b(a+2)-\varepsilon_2^2}}
\right)
\right)
\right).
\]
Since one can also verify that
\[
\partial_1\!\left(\frac{1}{K}\right)(0,u)=\frac{2}{(a+2)^2}\frac{(1-b)u+\varepsilon_1}{\big(1+\frac{\varepsilon_2}{a+2}u+\eta_1u^2\big)^2},
\]
it turns out that the function $M_1(u)=L_1(u)\partial_1\!\left(\frac{1}{K}\right)(0,u)$ is linear in $\varepsilon_1$. That being said, some computations show that
\begin{align*}
M_1(u)-M_1(0)&=\frac{2(1-b)}{(a+2)^2}u(1+{\eta_1} u^2)^{-\frac{2a+3}{a+2}}\\
&+\frac{2}{(a+2)^2}\left((1+{\eta_1} u^2)^{-\frac{2a+3}{a+2}}-1\right)\varepsilon_1
\\ 
&+\frac{2(1-b)u}{(a+2)^3\sqrt{b(a+2)}}(1+{\eta_1} u^2)^{-\frac{2a+3}{a+2}}\left(\arctan(\sqrt{{\eta_1}}u)-\frac{(2a+3)\sqrt{{\eta_1}}u}{1+{\eta_1} u^2}\right)\varepsilon_2
+\op(\|(\varepsilon_1,\varepsilon_2)\|).
\end{align*}
Following the obvious notation, if we write 
\[
 \int_0^{+\infty}\big(M_1(\pm u)-M_1(0)\big)u^{-1-1/\lambda} du=n_0^\pm+n_1^\pm\varepsilon_1+n_2^\pm\varepsilon_2+\op(\|(\varepsilon_1,\varepsilon_2)\|),
\]
then $n_0^-=-n_0^+$, $n_1^-=n_1^+$ and $n_2^-=n_2^+$ due to the parity of each coefficient with respect to $u.$ Here we follow exactly the same strategy as before, by using the results from \cite[\S 17.2]{Zorich} about improper integrals depending on a parameter, to show that the higher order terms can be neglected after integration. 
Moreover
\[
G_{1}=\int_{-1}^1\left(\frac{q_n(1,z)}{\h(1,z)}+1+\frac{1}{\lambda}
        +\frac{zq_n(z,1)}{\h(z,1)}\right)\frac{dz}{z}=\frac{2\pi\varepsilon_2}{(a+2)\sqrt{4b(a+2)-\varepsilon_2^2}},
\]        
so that $\exp(G_1)=1+\frac{\pi}{(a+2)\sqrt{b(a+2)}}\varepsilon_2+\op(\varepsilon_2).$ Accordingly, from \refc{lem2eq1bis} we can assert that
\begin{align*}
F_1&= -(n_0^++n_1^+\varepsilon_1+n_2^+\varepsilon_2)-\left(1+\frac{\pi}{\sqrt{b(a+2)^3}}\varepsilon_2\right)\left(n_0^-+n_1^-\varepsilon_1+n_2^-\varepsilon_2\right)+\op(\|(\varepsilon_1,\varepsilon_2)\|)\\
&=-n_0^--n_0^+-\left(n_1^-+n_1^+\right)\varepsilon_1-\left(n_2^-+n_2^++\frac{\pi}{\sqrt{b(a+2)^3}}n_0^-\right)\varepsilon_2+\op(\|(\varepsilon_1,\varepsilon_2)\|)\\
&=-2n_1^+\varepsilon_1-\left(2n_2^+-\frac{\pi}{\sqrt{b(a+2)^3}}n_0^+\right)\varepsilon_2
+\op(\|(\varepsilon_1,\varepsilon_2)\|).
\end{align*}
In order to compute this coefficients let us note that
\begin{align*}
 n_0^+&=\frac{2(1-b)}{(a+2)^2}\int_0^{+\infty}(1+{\eta_1} u^2)^{-\frac{2a+3}{a+2}}
 \frac{du}{u^{1/\lambda}}=\frac{1-b}{(a+1)(a+2)}\eta_1^{-\frac{a+1}{a+2}}\\
 \intertext{and}
 n_1^+&=\frac{2}{(a+2)^2}\int_0^{+\infty}\left((1+{\eta_1} u^2)^{-\frac{2a+3}{a+2}}-1\right)\frac{du}{u^{1+1/\lambda}}=\frac{\sqrt{\pi}}{2(a+2)^2}\frac{\Gamma\left(\frac{a}{2(a+2)}\right)}{\Gamma\left(\frac{2a+3}{a+2}\right)}\eta_1^{-\frac{a}{2(a+2)}}.
\end{align*}
The computations of $n_2^+$ is a little more involved. In this case
\begin{align*}
\frac{(a+2)^3\sqrt{b(a+2)}}{2(1-b)}\,n_2^+=&
\int_0^{+\infty}(1+{\eta_1} u^2)^{-\frac{2a+3}{a+2}}\arctan(\sqrt{{\eta_1}}u)u^{-1/\lambda}du\\
&-(2a+3)\sqrt{{\eta_1}}\int_0^{+\infty}(1+{\eta_1} u^2)^{-\frac{3a+5}{a+2}}
u^{1-1/\lambda}du\\
=&\,
\frac{\sqrt{\pi}(a+2)^2}{4(a+1)}\left(\frac{\Gamma\left(\frac{3a+4}{2(a+2)}\right)}{\Gamma\left(\frac{2a+3}{a+2}\right)}-\frac{\sqrt{\pi}}{a+2}\right)\eta_1^{-\frac{a+1}{a+2}},
\end{align*}
where to obtain the expression of the first integral we perform integration by parts.
From here some additional computations show that
\[
2n_2^+-\frac{\pi}{\sqrt{b(a+2)^3}}n_0^+=
\frac{\sqrt{\pi}(b-1)\eta_1^{-\frac{a+1}{a+2}}}{(a+1)\sqrt{b(a+2)^3}}\frac{\Gamma\left(\frac{3a+4}{2(a+2)}\right)}{\Gamma\left(\frac{2a+3}{a+2}\right)}
\]
and, on account of this,
\[
 \frac{2n_2^+-\frac{\pi}{\sqrt{b(a+2)^3}}n_0^+}{2n_1^+}=\frac{a(b-1)}{2(a+1)b}.
\]
Since $n_1^+<0$ for all $a\in (-1,0)$ and $b\in (0,2),$ we have that $F_1(\bar\mu)=\rho_4(\bar\mu)\left(\varepsilon_1+\frac{a(b-1)}{2(a+1)b}\varepsilon_2
+\op(\|(\varepsilon_1,\varepsilon_2)\|)\right)$ with $\rho_4(\bar\mu_0)>0$.
Accordingly, the combination of \refc{lem2eq0} and the last assertion in $(1)$ of \teoc{teo1} yields
\[
 \left.\Delta_1(\mu)\right|_{\varepsilon_0=0}=\kappa_{11}(\bar\mu)\left(\varepsilon_1+\frac{a(b-1)}{2(a+1)b}\varepsilon_2
+\op(\|(\varepsilon_1,\varepsilon_2)\|)\right)+\kappa_{12}(\bar\mu)\left.\Delta_0(\mu)\right|_{\varepsilon_0=0}
\]
with $\kappa_{11}(\bar\mu_0)>0.$ 
Finally, once again thanks to $\partial_{\varepsilon_0}\delta(\mu_{0})\neq 0$, we can write
\[
 \Delta_1(\mu)=\kappa_{11}(\bar\mu)\left(\varepsilon_1+\frac{a(b-1)}{2(a+1)b}\varepsilon_2
+\op(\|(\varepsilon_1,\varepsilon_2)\|)\right)+\kappa_{12}(\bar\mu)\Delta_0(\mu)+\kappa_{13}(\mu)\delta(\mu)
\]
for some smooth function $\kappa_{13}$ in a neighbourhood of $\mu_0$. This proves the last assertion in $(1)$ and completes the proof of the result.
\end{prooftext}

\bibliographystyle{plain}

\begin{thebibliography}{99}


\bibitem{ADL} J.C. Art\'es, F. Dumortier and J. Llibre, ``Qualitative theory of planar differential systems'', Universitext, Springer-Verlag, Berlin, 2006.

\bibitem{alien1} M. Caubergh, F. Dumortier and R. Roussarie, {\it Alien limit cycles in rigid unfoldings of a Hamiltonian 2-saddle cycle,} Commun. Pure Appl. Anal. {\bf 6} (2007) 1--21.

\bibitem{CLP09} B. Coll, C. Li, R. Prohens, {\it Quadratic perturbations of a class of quadratic reversible systems with two centers}, Discrete and Continuous Dynamical Systems {\bf 24} (2009) 699--729.

\bibitem{alien3} B. Coll, F. Dumortier, and R. Prohens, {\it Alien limit cycles in Li\'enard equations,} J. Differerential Equations {\bf 254} (2013) 1582--1600.

\bibitem{DGR02} F. Dumortier, A. Guzm\'an and C. Rousseau, 
{\it Finite cyclicity of elementary graphics surrounding a focus or center in quadratic systems,}
Qual. Theory Dyn. Syst. {\bf 3} (2002) 123--154.

\bibitem{alien2} F. Dumortier and R. Roussarie, {\it Abelian integrals and limit cycles}, J. Differential Equations {\bf 224} (2006) 296--313. 


\bibitem{DRR2} F. Dumortier, R. Roussarie and C. Rousseau,  {\it Hilbert's 16th problem for quadratic
vector fields,} J.~Differential Equations {\bf 110} (1994) 86--133.


\bibitem{FG} J.-P. Fran\c{c}oise and L. Gavrilov, {\it Perturbation theory of the quadratic Lotka-Volterra double center}, 
Commun. Contemp. Math. {\bf 24} (2022), no. 5, paper no. 2150064, 38 pp.
%\href{https://arxiv.org/abs/2011.08316}{arXiv:2011.08316v1}.

\bibitem{GMM02} A. Gasull, V. Ma\~nosa and F. Ma\~nosas, {\it Stability of certain planar unbounded polycycles,} 
J. Math. Anal. Appl. {\bf 269} (2002) 332--351. 


\bibitem{G08}  L. Gavrilov, {\it Cyclicity of period annuli and principalization of Bautin ideals}, Ergod. Th. \& Dynam. Sys. \textbf{28} (2008) 1497--1507.

\bibitem{alien8} L. Gavrilov and I. Illiev, {\it Perturbations of quadratic Hamiltonian two-saddle cycles,}
Ann. Inst. H. Poincar\'e C Anal. Non Lin\'eaire {\bf 32} (2015) 307--324.

\bibitem{Iliev} I.D. Iliev, {\it Perturbations of quadratic centers}, Bull. Sci. Math. \textbf{122} (1998) 107--161.

\bibitem{Ily} Y. Ilyashenko, {\it Centennial history of Hilbert's 16th problem,} Bull. Amer. Math. Soc. {\bf 39} (2002) 301--354.

\bibitem{Kelley} A. Kelley, {\it The stable, center-stable, center, center-unstable, unstable manifolds,} J. Differential 
Equations {\bf 3} (1967) 546--570.

\bibitem{Li} J. Li, {\it Hilbert's 16th problem and bifurcations of planar polynomial vector fields,} Internat. J. Bifur. Chaos Appl. Sci. Engrg. {\bf 13} (2003) 47--106.

\bibitem{alien4} S. Luca, F. Dumortier, M. Caubergh and R. Roussarie, {\it Detecting alien limit cycles near a Hamiltonian 2-saddle cycle,} Discrete Contin. Dyn. Syst. {\bf 4} (2009) 723--781.

\bibitem{Liu12} C. Liu, {\it The cyclicity of period annuli of a class of quadratic reversible systems with two centers}, J.~Differential Equations {\bf 252} (2012) 5260--5273.

\bibitem{Mou5} A. Mourtada,  {\it  Action de derivations irreductibles sur les algebres quasi-regulieres d'Hilbert}, preprint (2009), \href{https://arxiv.org/abs/0912.1560v1}{arXiv:0912.1560v1}.

\bibitem{MV19} D. Mar\'{\i}n and J. Villadelprat, \emph{Asymptotic expansion of the Dulac map and time for unfoldings of hyperbolic saddles: local setting,} J. Differential Equations {\bf 269} (2020) 8425--8467. 

\bibitem{MV20} D. Mar\'{\i}n and J. Villadelprat, {\it Asymptotic expansion of the Dulac map and time for unfoldings of hyperbolic saddles: general setting}, J. Differential Equations {\bf 275} (2021) 684--732.

\bibitem{MV21} D. Mar\'{\i}n and J. Villadelprat, {\it Asymptotic expansion of the Dulac map and time for unfoldings of hyperbolic saddles: coefficient properties},  J. Differential Equations {\bf 404} (2024) 43--107.

\bibitem{MV22} D. Mar\'{\i}n and J. Villadelprat, {\it The criticality of reversible quadratic centers at the outer boundary of its period annulus}, J. Differential Equations {\bf 332} (2022), 123--201.

\bibitem{Munkres} J.R. Munkres, ``Topology: a first course'', Prentice-Hall, Inc., Englewood Cliffs, NJ, 1975

\bibitem{PFL14}  L. Peng, Z. Feng and C. Liu, {\it Quadratic perturbations of a quadratic reversible Lotka-Volterra system with two centers}, Discrete and Continuous Dynamical Systems {\bf 34} (2014) 4807--4826.

\bibitem{Roussarie} R. Roussarie, ``Bifurcations of planar vector fields and Hilbert's sixteenth problem'' [2013] reprint of the 1998 edition. Modern Birkh\"auser Classics. Birkh\"auser/Springer, Basel, 1998.

\bibitem{RouRou08} R. Roussarie and C. Rousseau, {\it Finite cyclicity of nilpotent graphics of pp-type surrounding a center,}  Bull. Belg. Math. Soc. Simon Stevin {\bf 15} (2008) 889--920.

\bibitem{Rou} C. Rousseau, {\it Normal forms, bifurcations and finiteness properties of vector fields},
NATO Sci. Ser. II Math. Phys. Chem. {\bf 137}, Kluwer Academic Publishers, Dordrecht (2004) 431--470.

\bibitem{Rudin} W. Rudin, ``Real and complex analysis'' McGraw-Hill Book Co., New York-Toronto, Ont.-London 1966.

\bibitem{ShaZeg94} D. S. Shafer and A. Zegeling, {\it Bifurcation of limit cycles from quadratic centers}, J. Differential Equations {\bf 122} (1995) 48--70.

\bibitem{alien5} L. Sheng and M. Han, {\it Bifurcation of limit cycles from a compound loop with five saddles,} J. Appl. Anal. Comput. {\bf 9} (2019) 2482--2495.

\bibitem{alien6} L. Sheng, M. Han and Y. Tian, {\it On the number of limit cycles bifurcating from a compound polycycle,} Int. J. Bifur. Chaos Appl. Sci. Eng. {\bf 30}  (2020)
no. 7, paper no. 2050099, 16 pp.

\bibitem{Swirszcz} G. Swirszcz, {\it Cyclicity of Infinite Contour around Certain Reversible Quadratic Center}, J. Differential Equations {\bf 154} (1999) 239--266.

\bibitem{alien7} J. Yang, Y. Xiong and M. Han, {\it Limit cycle bifurcations near a 2-polycycle or
double 2-polycycle of planar systems,} Nonlinear Anal. {\bf 95}  (2014) 756--773.

\bibitem{Zorich} V. A. Zorich, ``Mathematical analysis II'' Translated from the 2002 fourth Russian edition by Roger Cooke. Universitext. Springer-Verlag, Berlin, 2004.


\end{thebibliography}

\end{document}